# Directional Deduction

Roderick Batchelor

Department of Philosophy, University of São Paulo

*Abstract*.  We present some new methods for logical deduction, based on ideas from ground theory. Roughly speaking, in our calculi a typical deduction will proceed as follows: we first analyse the premiss down to its ultimate grounds; then we discard information irrelevant to the conclusion; and then we synthesize the conclusion up from *its* ultimate grounds. – We give a series of calculi for: classical propositional logic (Chapter 1); classical predicate logic (Chapter 2); modal propositional logic and modal predicate logic (Chapter 3); and some relevantistic fragments of these various systems (Chapters 4 and 5). In connection with these fragments we develop also some new semantic constructions, of 'truthmaker semantics' type.

# CONTENTS

# CHAPTER 1

## DIRECTIONAL DEDUCTION FOR CLASSICAL PROPOSITIONAL LOGIC

### 1. Introduction

(The present introductory section expounds a line of reasoning which provides natural motivation for the construction of our Data Calculus and our other, related calculi. This is by no means the only line of reasoning which might lead to this construction; but it is perhaps a particularly natural one, and corresponds approximately to the way in which I myself originally hit upon the idea of the Data Calculus.)

p ∧ q implies p. Also, p ∧ q implies (p ∧ q) ∨ r. But there is an important difference between these two cases of implication. (The present remarks can be understood in at least two different ways: in terms of formulas and a semantical notion of implication; or in terms of a fine-grained structural notion of proposition [situation, state of affairs] and say the modal notion of strict implication.) In passing from p ∧ q to p we move to what seems to be in some reasonable sense a '*part of the content*' of the original proposition – something obtainable from it by a kind of *analysis*. On the other hand in passing from p ∧ q to (p ∧ q) ∨ r (or more generally from A to A ∨ B), far from moving to a 'part of the content' by a kind of 'analysis', on the contrary we seem to be in some sense *expanding* on the content and effecting a kind of *synthesis*. It is not that there isn't *implication* just as much in the second case as in the first; but the *direction* (as we move from implying proposition to implied one) seems to be different: *downward* in the first case and *upward* in the second. It is natural to represent this difference graphically thus:

$$(p \wedge q) \vee r$$
$$\uparrow$$
$$p \wedge q$$
$$\downarrow$$
$$p$$

Or in the general case for conjunctions and disjunctions:

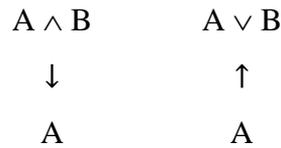

Double negation gives another typical example of this phenomenon. ¬¬p 'analytically implies' p – p is a kind of 'part of the content' of ¬¬p, obtainable from it by 'analysis' –; whereas p 'synthetically implies' ¬¬p – the latter is obtainable from the former by a kind of 'expansion' or 'synthesis':

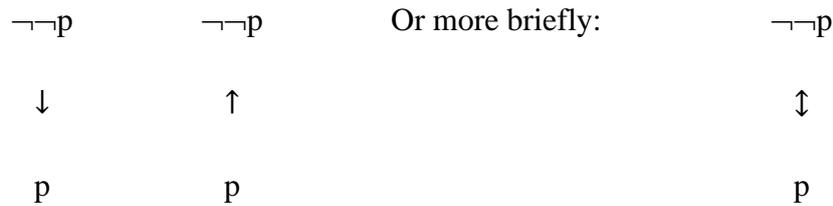

It is not that all cases of implication should be classified either as analytic or as synthetic. E.g. p ∧ q implies q ∨ r but this is neither purely analytic nor purely synthetic implication. However a natural way of deducing q ∨ r from p ∧ q – viz. inferring q from p ∧ q, then q ∨ r from q – consists of one purely analytic step and then one purely synthetic step. (It may not even be too hasty to suggest at once that perhaps the '*atoms*' of 'inferential connection' are always either purely analytic or purely synthetic.)

So these intuitive ideas of analytic and synthetic implication seem to make some sense, and if legitimate are obviously of considerable interest. So how might we try to *define* (in precise terms) concepts of analytic and synthetic implication?

An idea which very naturally occurs here is to try to use the idea of *grounding*, of some things being the case *in virtue of* or *because of* other things being the case. The definitions then would be something like this: *analytic* implication is implication from groundeds to grounds; *synthetic* implication is implication from grounds to groundeds. Thus e.g. from p ∧ q to p we move from grounded to ground; from p to ¬¬p from ground to grounded.

There is however an immediate general difficulty with this idea. Analytic and synthetic implication should be relations correctly applicable to propositions in general whether they be true or false propositions: e.g. p ∧ q should analytically imply p for *any* propositions p and q; p synthetically imply ¬¬p for *any* proposition p. But *grounding* by contrast is usually considered as a '*factive*' relation, i.e. a relation which holds only between *facts* or *true* propositions or the like. It is *not* the case e.g. that $2 + 2 = 5$ is a ground for $\neg\neg(2 + 2 = 5)$, or in other words that $\neg\neg(2 + 2 = 5)$ because (or at least in part because) $2 + 2 = 5$.

Still, perhaps this difficulty is not unsurmountable. Perhaps we can make sense of a kind of '*neutral*' notion of grounding – something *like* the usual notion of grounding but with the facticity requirement abstracted away. It sounds plausible to say that if p (whence also $\neg\neg$p) is the case, then: $\neg\neg$p is the case because p is the case. But it also sounds plausible to say that if p is or were (maybe even *per impossibile*) the case, then: $\neg\neg$p is or would be the case because p is or would be the case. I.e. it seems plausible to say that, for any proposition p (true or false), p is 'neutral ground' of $\neg\neg$p. Or again, that for any propositions p, q we have: p, q neutrally ground $p \wedge q$.

Indeed these two cases (double negation and conjunction) give almost the whole of a case-by-case definition of 'immediate neutral grounding' (for 'truth-functional compounds') similar to the case-by-case definition of immediate (factive) grounding in my paper Batchelor 2010. (And in the factive case the corresponding notion of 'general' [not necessarily immediate] ground can be defined as the transitive closure of the given notion of immediate ground.) The immediate grounding connections

$$\frac{\neg\neg\text{p}}{\text{p}} \qquad \frac{\text{p} \wedge \text{q}}{\text{p, q}}$$

are as reasonable for the 'neutral' notion as for the 'factive' one.

However, in the remaining case of *negations of conjunctions*, the difference between factive and neutral context becomes relevant. In my earlier paper I took the (factive) immediate grounds of a fact $\neg(p \wedge q)$ to be: the facts among $\neg$p, $\neg$q. Or if we take all of negation, conjunction and disjunction as primitive connections, we have also (in addition to the other clauses already mentioned and a clause for negations of disjunctions directly similar to the clause for conjunctions and also equally suitable in neutral as in factive context): the immediate grounds of a fact $p \vee q$ are: the facts among p, q. Here of course facticity plays a significant role, and it is not immediately clear what if anything might be done in the neutral context.

Some examples may be useful here. Suppose we ask: Why is it the case that

$$(\text{Prawitz is Swedish}) \lor (\text{Prawitz is Norwegian})?$$

(In the sense of 'immediate' reason.) It can then sensibly be answered: Because Prawitz is Swedish. Or again, we may ask: Why is it that

$$(\text{Skolem is Swedish}) \lor (\text{Skolem is Norwegian})?$$

And it can sensibly be answered: Because Skolem is Norwegian. (And for someone with the relevant double nationality both disjuncts would constitute the immediate grounds.) But now suppose we ask: What neutrally grounds

$$(\text{Quine is Swedish}) \lor (\text{Quine is Norwegian})?$$

Here it is less clear what can sensibly be answered. We can of course say what would be the factive grounds in the different circumstances under which the disjunction would be true. But what we would like to say is not what *would be* the *factive* grounds in different circumstances but what *are* the *neutral* grounds! Indeed even in the case of *true* disjunctions, like the Prawitz and Skolem examples above, it is not clear what might be said to *be* the *neutral* grounds.

One might think that these considerations show that facticity is essential to grounding and that the whole idea of neutral grounding, at least as a fully general notion, is misguided and should be abandoned. But I don't think one should be so pessimistic. Let us consider again the case of conjunction:

$$\frac{p \land q}{p, q}$$

This can be read as saying: the proposition p and the proposition q are the immediate neutral grounds of the proposition $p \land q$. Or alternatively (and more germanely to present purposes) it can be read as saying: the immediate neutral ground for it being the case that $p \land q$ is: it being the case that p and it being the case that q. Thus the (putative) 'datum'

It is the case that $p \land q$

is taken to be immediately neutrally grounded in the (putative) datum

It is the case that p and it is the case that q.

Here we have as it were transferred the complexity from the internal structure of the proposition (viz. p ∧ q) to the external structure of the datum; and plausibly this constitutes a grounding analysis. – But then, once things are put in this way, it is clear that a similar clause can be used for *dis*junction. We can put

$$\frac{p \vee q}{p \mathbin{/} q}$$

meaning: the immediate neutral ground for it being the case that p ∨ q is: it being the case that p or it being the case that q. – Of course 'p / q' here does not represent some particular proposition, but rather the 'alternational datum': it is the case that p or it is the case that q. But then also before 'p, q' did not represent a particular proposition but rather the 'list-like datum': it is the case that p and it is the case that q.

– This then is our idea: to use here a notion of '*datum*' including both 'list-like' data which give immediate (neutral) ground for conjunctions and 'alternational' data which give immediate (neutral) ground for disjunctions. Now in giving a precise implementation of this idea, one might consider permitting free iterations of list-like and alternational constructions. But that does not seem to give a very useful notion. The better course, and the one which I will adopt, is to allow a single layer of lists (of propositions, or in linguistic context formulas) inside a single alternation. I.e. the datum is a set of sets of propositions, the inner set-formation being understood as 'conjunction-like' and the outer as 'disjunction-like'. (We will also consider a notion of 'dual datum', with 'conjunction-like' and 'disjunction-like' inverted.)

With this notion of datum in place, one quickly sees that the simple (equivalential) grounding connections indicated above, each corresponding to both an analytic (top to bottom) and a synthetic (bottom to top) inference rule – one sees that this is already a very large chunk of a sufficient (complete) logical calculus. One needs only a few strategic additions. Obviously one must be able to weaken data (since of course sometimes a weaker datum follows from a stronger one) – but our 'structural devices', the comma and the slash, are naturally associated with corresponding weakening rules: an analytic Deletion rule to reduce a 'list', and a synthetic Expansion rule to prolong an 'alternation'. And with only a couple more strategic additions related to negation, we are done.

(We are now firmly in the context of 'ordinary logic', and this will be the context of the whole of the present work. But it should be clear that some analogous developments are possible in the analogous 'metaphysical' context of structured propositions – though this raises some delicate issues. For example: here it seems clear that [neutral] grounding cannot adequately be defined as the transitive closure of [neutral] immediate grounding. For take e.g. the infinite alternational datum

$$p / \neg\neg p / \neg\neg\neg\neg p / \neg\neg\neg\neg\neg\neg p / \ldots$$

Surely p should be neutral ground [indeed *the* ultimate neutral ground] of this datum, but no finite number of double negation eliminations will take us to p.)

A crucial characteristic of this 'Data Calculus' is the possibility of reducing an arbitrary datum, through the analytic grounding inference rules, to an 'elementary' datum where no connectives occur at all except perhaps negations immediately preceding atomic formulas – i.e. a sort a datum-correlate of a disjunctive normal form formula. – Reflection on the analogy between disjunctive normal forms in propositional logic and Skolem normal forms in predicate logic then led us to the construction of our 'Quantificational Data Calculus' QDC (and cognate calculi) presented in Chapter 2 below. And that in turn, together with the familiar correspondence of necessity and possibility with the quantificational notions of truth in *every* possible world and truth in *some* possible world, led us to our Modal Data Calculus S5DC and its quantificational version S5QDC (and cognate calculi), presented in Chapter 3. (I will not here speculate on the question of possible 'metaphysical significance' of *these* constructions [quantificational and modal cases].)

The relevantistic calculi in Chapters 4 and 5 were a kind of afterthought. I noticed that dropping a couple of rules from the Data Calculus yields a natural calculus for the Belnap–Dunn logic of First-Degree Entailment (FDE). This logic is one of an interesting family of relevantistic logics treatable by the method of 'truthmaker semantics' – a family which includes also e.g. R. B. Angell's equivalential logic and Fabrice Correia's 'logic of factual equivalence'. So I was lead to calculi also for these other logics (and modal and quantificational extensions). (E.g. the Data Calculus stripped down to just the grounding rules corresponds to Correia's system!)

*Remarks* (on literature etc.).  (1) The idea of 'partial content', i.e. (the converse of) what I have called here *analytic* implication, is I suppose a generally familiar idea. However it has usually not been considered in conjunction with the cognate notion of 'expansive' or *synthetic* implication. An exception (i.e. someone who *has* considered something like these two ideas in conjunction) is Fine: see Fine 2015 § 4 on 'constitutively necessary conditions' and 'constitutively sufficient conditions', and Fine 2017a § 5 on 'containment' and 'entailment'.

(2) I should perhaps say that analytic implication in my sense of the term here is quite different from analytic implication in the sense of Parry (on which see e.g. Fine 1986). There the idea is something like this: p strictly implies q, and all basic non-logical constituents of q are also basic non-logical constituents of p. So for instance with p and ¬¬p there will always be mutual 'analytic implication' (or as we may say 'analytic equivalence') in this sense. – No doubt analytic implication in my sense implies analytic implication in this Parryan sense, though not of course vice versa. – A notion along somewhat similar lines and closer to analytic implication in my sense is: p strictly implies q, and q is a constituent of p. But this is not quite what we want. An obvious counterexample is ¬(p ∨ q) and ¬p: we would like to say that the former analytically implies the latter, but ¬p is not quite a constituent of ¬(p ∨ q). In fact however *such* counterexamples are not the most serious: just as in some similar contexts one extends the notion of 'subformula' so that negations of subformulas-in-the-literal-sense also count as subformulas, so also it is reasonable enough here to consider the modification of the above condition where 'q is constituent of p' is replaced by 'either q is constituent of p or ¬q is constituent of p' – which disposes of the counterexamples in question. But consider now the case of ¬p ∧ (q ∧ ¬q) and p: here there is strict implication, and even the stronger (and of course a fortiori the weaker) constituency condition holds; but we would not want to say that this is a case of analytical implication in our sense – p is not a 'part of what is said' in ¬p ∧ (q ∧ ¬q).

(3) Just as strict implication is the necessary truth of the material implication (□(p → q)), so also one might reasonably speak of analytic implication in the sense of the analytic truth of the material implication (An(p → q)) (supposing analytic truth is a legitimate notion, which is of course debatable). (So this would be at least as strong as, and perhaps stronger than, strict implication.) This too is of course very different from analytic implication in the sense in which we are using the term here. E.g. An(p → (p ∨ q)) is no less always-true than An((p ∧ q) → p).

(4) Neutral grounding has been occasionally discussed in the recent literature on ground. See e.g. Fine 2012 § 5 pp. 48–50, Fine 2022, and Correia 2010. (Fine uses the term '*non-factive* grounding' for this notion. I prefer 'neutral grounding' because (i) 'non-factive grounding' carries the – of course unintended – suggestion that the terms are *not* factive, and (ii) perhaps the notion is important enough to deserve a 'positive' name.)

(5) The general idea of using some kind of 'structural device' corresponding to disjunction is of course not new here. It goes at least as far back as Gentzen's sequent calculus; it is present of course in Smullyan's tableaux and its ancestors from Beth and Hintikka; and it appears also in the method of 'multiple-conclusion logic' pioneered by Carnap (1943) and Kneale (1956) and later investigated by some others – see esp. Shoesmith and Smiley 1980.

(6) The notion of what we call *data* in classical propositional logic was considered in Smullyan 1968 Ch. 17 under the name *configurations*, and earlier in Smullyan 1963 § 4 under the name *statements*. But Smullyan, in accordance with his 'purely analytic' tendencies as in his tableaux, considers only what we will call the 'downward G-rules' and notes their completeness for refutation of unsatisfiable 'configurations' in the sense of reduction to 'explicitly contradictory' configurations (i.e. with each component containing some pair of formulas A, ¬A). So *this* is hardly more than a reformulation of tableaux in a slightly different style. Also Smullyan does not consider anything like my notions of data for classical predicate logic and modal logic.

(7) The method of 'coupled trees' for classical propositional logic in Jeffrey 1967 pp. 93 ff. (by the way, unfortunately apparently omitted from later editions of the book) has some definite similarities with our Data Calculus. But it is less flexible in that there analysis needs always to be carried out to the end (as in our *canonical* Data Calculus deductions). Also Jeffrey gives no quantificational or modal version of the method.

(8) The calculus for classical predicate logic in the paper Craig 1957 should also be mentioned here as having significant similarities to our methods, in particular in the following points: linearity; only rules (no logical axioms); extensive use of equivalential rules; all rules 'correct' (i.e. premiss implying conclusion); rules with close control of the structure of formulas, and relatedly possibility of 'normal' deductions, this among other things permitting proof of interpolation theorem. (Indeed this is the paper which gave the original proof of that theorem for classical predicate logic.)

(9) For references concerning the relevantistic systems see the Historical Notes at the end of the introductory section in Chapter 4 below.

## 2. The Data Calculus (DC)

We use a standard language for Classical Propositional Logic, with propositional variables p, q, r, … and with ¬, ∧, ∨ as the primitive connectives; *formulas* are then defined in the usual way.

Datum $=_{df}$ finite set of finite sets of formulas.

We use α, β, γ, … as variables for data.

By the *components* of a datum we mean its elements.

We call a datum *elementary* if all its formulas (i.e. all elements of all components) are literals. We say that an elementary datum is *purified* if none of its components includes a variable together with its negation. And we say that a purified elementary datum α is a *canonical elementary datum* if L(Γ) = L(Δ) for all components Γ and Δ of α. (We use here L for 'the set of propositional variables occurring in'.)

We usually represent data using slashes to separate the components, and commas to separate the formulas inside each component. Thus e.g. instead of {{p, q}, {p, ¬q}, {¬p}, {r ∨ ¬r}} we write

$$p, q \ / \ p, \neg q \ / \ \neg p \ / \ r \vee \neg r.$$

We use also some conventions of a familiar kind to lighten notation, such as writing a component as Γ, A instead of Γ ∪ {A}. We use also sometimes a Greek letter such as ρ to represent the 'rest' of a datum: thus e.g. p / ¬p / ρ represents data of the form {{p}, {¬p}, …}. Note that this ρ need not be *one* component: it can also be more, or less – thus e.g. just p / ¬p is of the form p / ¬p / ρ.

We say that an interpretation (i.e. assignment of truth-values to the propositional variables) σ *verifies* a datum α if there exists a component Γ of α s.t. σ verifies every formula in Γ. Thus we understand a datum as having the 'effect' of the corresponding disjunction of conjunctions. We then define in obvious way *valid* datum, *satisfiable* datum, etc. etc.

We use ⊥ for { }, and ⊤ for {{ }}. Clearly ⊥ is unsatisfiable, and ⊤ valid.

– We list now the rules of our *Data Calculus* (DC). They are divided into three groups: *Grounding*-rules or for short G-rules; *Negation*-rules or N-rules (rules deriving from the general character of negation); and *Structural*-rules or S-rules. Rules represented with ↓ we call *analytic*, and those represented with ↑ *synthetic*. (The ↕-figures in the group of G-rules represent each a pair of rules, one analytic and the other synthetic.) Notice that the synthetic rules are depicted with premiss at the bottom and conclusion at the top. – The abbreviated names NC, EM, Del, and Exp stand for respectively Non-Contradiction, Excluded Middle, Deletion, and Expansion.

<u>G-rules:</u>

|  | Γ, A ∧ B / ρ |  | Γ, A ∨ B / ρ |  | Γ, ¬¬A / ρ |
|---|---|---|---|---|---|
| (↕∧) | ↕ | (↕∨) | ↕ | (↕¬¬) | ↕ |
|  | Γ, A, B / ρ |  | Γ, A / Γ, B / ρ |  | Γ, A / ρ |

|  | Γ, ¬(A ∨ B) / ρ |  | Γ, ¬(A ∧ B) / ρ |
|---|---|---|---|
| (↕¬∨) | ↕ | (↕¬∧) | ↕ |
|  | Γ, ¬A, ¬B / ρ |  | Γ, ¬A / Γ, ¬B / ρ |

<u>N-rules:</u>

|  | Γ, A, ¬A / ρ |  | Γ, A / Γ, ¬A / ρ |
|---|---|---|---|
| (NC) | ↓ | (EM) | ↑ |
|  | ρ |  | Γ / ρ |

<u>S-rules:</u>

|  | Γ, Δ / ρ |  | Γ / ρ |
|---|---|---|---|
| (Del) | ↓ | (Exp) | ↑ |
|  | Γ / ρ |  | ρ |

Clearly each of these rules is *correct* in the sense that, for any premiss-datum and conclusion-datum instantiating the indicated form, if an interpretation verifies the premiss-datum then it also verifies the conclusion-datum. Each of the G-rules and N-rules, though not either of the S-rules, is also *equivalential* (or *equivalentially correct*) in the sense that both it and also its 'inverse' are correct. Indeed in the case of a G-rule its inverse is itself a (correct) primitive rule of our calculus. The inverses of NC and EM are readily derivable using respectively Exp and Del.

A *deduction* (in DC) of datum β from datum α is a finite sequence of data, starting with α and ending with β, s.t. each datum in the sequence after the first is immediately inferable from the immediately preceding one by one of the rules listed above.

(Note that no significant gain would result from liberalizing the notion of deduction by replacing 'from *the immediately preceding* one' by 'from a previous one' in this definition. For, given any deduction of β from α in this liberalized sense, a deduction in our official sense can be obtained by simply successively erasing, from the bottom up, all the blocks of intermediate lines between such 'jumping points'. – Note also that this is not just a point about the Data Calculus but applies in general to any calculus of 'linear reasoning', i.e. with only single-premiss rules.)

A *proof* of α is a deduction of α from ⊤. A *refutation* of α is a deduction of ⊥ from α.

An *analytic* (or for emphasis *purely analytic*) deduction is a deduction via analytic rules only; similarly for a *synthetic* (or *purely synthetic*) deduction. A *normal* deduction is a deduction divisible into two consecutive parts of which the first is purely analytic and the second purely synthetic. (Each of these parts is allowed to be empty; so purely analytic and purely synthetic deductions count as normal deductions. Note also that clearly a normal *proof* must be purely synthetic, and a normal *refutation* purely analytic.)

We say that β is *deducible* (in DC) from α, or α ⊢ β (or α ⊢$_{DC}$ β in contexts where other calculi are also being considered), if there exists some deduction of β from α. Also α ⊢$^N$ β if there exists a *normal* deduction of β from α. We use also ⊢ α for ⊤ ⊢ α, and ⊢$^N$ α for ⊤ ⊢$^N$ α.

We represent deductions by writing the relevant data in successive lines, with an indication at the end of each line (other than the last) of rule licensing passage from the given datum to the next. The names of G-rules include already an ↑ or ↓ indicating the 'direction' of the deduction-step; in writing deductions we adjoin also ↑ or ↓ (as appropriate) to the names of other rules, so that the 'directional structure' of the whole deduction becomes easily visualizable. – We sometimes skip some lines when the deduction-steps are very simple, and then indicate on the right the different rules, and/or multiple applications of a single rule, used through the omitted lines. (See e.g. the second example below.)

Here are some examples of deductions.

$\underline{p \land q \vdash^N q \lor r}$   (For we have the [normal] deduction:)

$p \land q$        $\downarrow \land$

$p, q$        $\downarrow$ Del

$q$        $\uparrow$ Exp

$q \,/\, r$        $\uparrow \lor$

$q \lor r.$

$\underline{(p \land q) \lor (p \land r) \vdash^N p, q \lor r}$

$(p \land q) \lor (p \land r)$        $\downarrow \lor$

$p \land q \,/\, p \land r$        $\downarrow \land \ (\times 2)$

$p, q \,/\, p, r$        $\uparrow \lor$

$p, q \lor r.$

$\underline{\vdash^N p \lor \neg p}$

$\top \ (= \{\{\ \}\})$                    $\uparrow$ EM

$p \,/\, \neg p \ (= \{\ \} \cup \{p\} \,/\, \{\ \} \cup \{\neg p\})$       $\uparrow \lor$

$p \lor \neg p.$

$\underline{p \land \neg p \vdash^N q}$

$p \land \neg p$            $\downarrow \land$

$p, \neg p \ (= \{\{p, \neg p\}\})$      $\downarrow$ NC

$\bot \ \ \ \ (= \{\ \})$          $\uparrow$ Exp

q        (= {{q}}).

(We also use formulas with → and ↔ in deductions, with the understanding that
A → B abbreviates ¬A ∨ B, and A ↔ B abbreviates (A ∧ B) ∨ (¬A ∧ ¬B).)

<u>p ↔ q, q ↔ r ⊢$^N$ p ↔ r</u>

p ↔ q, q ↔ r                                                                ↓ ∨

p ↔ q, q ∧ r / p ↔ q, ¬q ∧ ¬r                                               ↓ ∧ (×2)

p ↔ q, q, r / p ↔ q, ¬q, ¬r                                                 ↓ ∨ (×2)

p ∧ q, q, r / ¬p ∧ ¬q, q, r, / p ∧ q, ¬q, ¬r / ¬p ∧ ¬q, ¬q, ¬r             ↓ ∧ (×4)

p, q, q, r / ¬p, ¬q, q, r / p, q, ¬q, ¬r / ¬p, ¬q, ¬q, ¬r                   ↓ NC (×2)

p, q, r / ¬p, ¬q, ¬r                                                        ↓ Del (×2)

p, r / ¬p, ¬r                                                              ↑ ∧ (×2)

p ∧ r / ¬p ∧ ¬r                                                            ↑ ∨

p ↔ r.

<u>(p ∧ q) ∧ (r ∧ s) ⊢$^N$ (p ∨ q) ∧ (r ∨ s)</u>

(p ∧ q) ∧ (r ∧ s)              ↓ ∧ (×3)

p, q, r, s                     ↓ Del

p, r                          ↑ Exp

p, r / q, r                   ↑ ∨

p ∨ q, r                      ↑ Exp

p ∨ q, r / p ∨ q, s           ↑ ∨

p ∨ q, r ∨ s                  ↑ ∧

(p ∨ q) ∧ (r ∨ s).

$\underline{p \rightarrow q, (p \land q) \rightarrow r \vdash^N p \rightarrow r}$

$\neg p \lor q, \neg(p \land q) \lor r$                                   $\downarrow \lor$

$\neg p \lor q, \neg(p \land q) / \neg p \lor q, r$                       $\downarrow \lor (\times 2)$

$\neg p, \neg(p \land q) / q, \neg(p \land q) / \neg p, r / q, r$         $\downarrow$ Del

$\neg p / q, \neg(p \land q) / \neg p / r$                                $\downarrow \neg\land$

$\neg p / q, \neg p / q, \neg q / \neg p / r$                             $\downarrow$ NC

$\neg p / q, \neg p / \neg p / r$                                         $\downarrow$ Del

$\neg p / \neg p / \neg p / r$          $(= \neg p / r)$                  $\uparrow \lor$

$\neg p \lor r.$

$\underline{p, q \vdash p, r / q, \neg r}$

$p, q$                          $\uparrow$ EM

$p, q, r / p, q, \neg r$        $\downarrow$ Del $(\times 2)$

$p, r / q, \neg r.$

*Remarks*.  We list here summarily ten notable features of DC – actually also features of all or almost all the other calculi to be presented later in this work –, with a few brief comments.

(1) *Ground-tracking* rules; (relatedly) *directionality* of rules. And that both ways, *analytic* and *synthetic*. (In contrast to e.g. tableaux and ['cut-free'] sequent calculus, where the rules are respectively only analytic and only synthetic. Our 'Consecution Calculi' to be defined later also have only synthetic rules.) – The G-rules are analytic or synthetic in a strict, as it were 'vertical' or 'hierarchical' sense; whereas the N-rules, and certainly the S-rules, are analytic or synthetic in a less strict, as it were 'horizontal' or 'flat' sense. (I mean this for *all* our calculi. In predicate logic and modal logic it should be observed in this connection that the 'universal prefix' is meant as a kind of 'generalized *comma*', not 'generalized *conjunction*', and similarly for the existential prefix – meant as generalized *slash* not generalized *disjunction*. – In modal logic however some doubts might reasonably be raised as to the status [horizontal or vertical] of the ↕T-rules.) The difference is analogous to that between constituency as in my 'Stoichiology' (Batchelor 2013) and flat part-whole in the sense of Mereology.

(2) *Rules only* (no 'logical axioms'). (Nevertheless we can still define [as we have done] a 'proof' of a datum α as a deduction of α from the 'canonical tautology' ⊤.) – On this basis 'formalized theories' of the relevant kind can be constructed by adjunction of a set of 'extra-logical assumptions' (axioms). As usual here all basic *rules* are logical rules; but also now all basic *statements* (axioms) are extra-logical. Logic gives the 'form' (inferential connections) and the form only, in especially neat conformity with the familiar 'sentiments' about 'the nature of logic'.

(3) The rules are intrinsically *natural* and even arguably *fundamental*. – We may call *Inferential Atomism* (or maybe 'Inferential Fundamentalism') the view that logical rules of inference, in themselves and objectively (independently of our activities of system-construction), can be more or less fundamental, and form an objective structure with the basic rules at the bottom and more and more complex derived rules springing from that. And if we think that that is the case then the task arises of course of trying to determine what are those basic rules – those 'atoms of inferential connection' –, and the Data Calculus is maybe not a bad first stab at that.

(4) *Easy to find deductions* (including proofs and refutations) (and in propositional case mechanical). This is related to point (3): for surely the closer we get to the bottom of things the easier it becomes to find our way around the relevant tasks.

(5) All rules are *normal* in the sense that instances of the premisses always *imply* the corresponding instance of the conclusion. (This is especially relevant in the cases of predicate logic and modal logic, where most extant formalizations violate the requirement and have e.g. rule of 'Universal Generalization' or 'Existential Instantiation', or in modal logic the ubiquitous Rule of Necessitation.) – Non-normal rules may no doubt be convenient for the proof of purely logical truths (or refutation of purely logical falsehoods). But of course they raise problems in application (in 'formalized theories'). Also, recalling again point (3), it is hard for me to believe that a non-normal rule could be a *fundamental* inference-rule out there in the world (or up there in the Platonic World).

And anyway, one is inclined to say: surely it is not too much to ask of a rule of inference, that the premisses should *imply* the damn conclusion!

(6) *No 'temporary assumptions'* (later 'discharged') within deductions. – Again recall point (3); it is hard for me to believe that this procedure could be involved in the *ultimate* structure of deductions. (This is an important difference between our methods and 'Natural Deduction'.)

(7) *Linearity*: our rules always have a *single* datum as premiss (and a single datum as conclusion). So our deductions have a straightforward linear structure rather than the more common and more complex 'tree structure'.

(8) Extensive use of *equivalential rules*. As we have seen, in DC all G-rules and N-rules are equivalential, and weakening of strength may occur only with the S-rules (and a similar situation obtains in our other Directional Deduction calculi, although there often even some of the S-rules are also equivalential). – This is certainly an interesting aspect of our calculi; it is particularly striking that, with the introduction of our concept of datum, all Grounding-rules turn out to be equivalential. We obtain an ingenious division of inferential labour, with the G-rules and N-rules (and equivalential S-rules when there are such) responsible for analyses and syntheses of equivalent matter, and the (non-equivalential) S-rules for weakenings of strength.

(9) *Special status of negation*. And that in two ways: (i) in the fact that the G-rules roughly speaking come in pairs, one for a given 'logical constant' and the other for the corresponding negation; and (ii) in the presence of the N-rules deriving from the general character of negation. – Philosophical logicians seem to have often felt that negation has somehow a special position and is not just 'a logical constant like any other'. If one does feel this, one may be pleased by this feature of our calculi.

(10) *Complex 'formal objects'* (not mere formulas). – This has to be counted as 'on the minus side': formulas are of course simpler objects than data and so would be *ceteris paribus* preferable. Actually one could reformulate DC in terms of formulas only. (Replace commas and slashes by $\wedge$'s and $\vee$'s; drop $\updownarrow\wedge$ [since now the two sides are just identical]; and introduce formulas $\top$, $\bot$ [so we can go e.g. to $\bot$ from $p \wedge \neg p$ by the new version of NC].) But this would destroy much of its distinctive character and the direct connection with grounding including of course grounding for conjunctions. – We may note here two interesting differences between DC etc. and the sequent calculus (the most famous predecessor logical calculus with 'complex formal objects'): (i) notationally, DC etc. are more straightforward: conjunction-like and disjunction-like connections are directly notationally distinguished, whereas in sequent calculus the same symbol (the comma) is associated with different ideas depending on whether it is on the left or right; (ii) just as predicate logic or modal logic have more 'logical constants' than classical propositional logic, so in our calculi the 'structural' resources (which as in sequent calculus have relation to some of the logical constants) will also increase, whereas in sequent calculus they remain the same.

### 3. The Dual Data Calculus (DDC) and the Hybrid Calculus (HC)

*Data*, as we have seen, have the 'effect' of a disjunction of conjunctions. We consider now the dual notion of *dual data*, with the 'effect' of a conjunction of disjunctions. We will formulate both a Dual Data Calculus (DDC) where the 'formal objects' are just dual data, and also a Hybrid Calculus (HC) whose 'formal objects' are both data and dual data and where deductions start (typically analytically) with data and then transition to dual data proceeding with them (typically synthetically) till the end.

For the purposes just of DDC we can define:

Dual datum $=_{df}$ finite set of finite sets of formulas.

But in contexts, such as HC, where we are dealing both with data and with dual data, it may be desirable to change the official definitions somehow so as to differentiate data from dual data. E.g. we may define data as objects of the form $\langle 0, S \rangle$ and dual data as objects of the form $\langle 1, S \rangle$ (for S finite set of finite sets of formulas). – To avoid fuss I will normally speak of a datum or a dual datum as *being* the set of sets of formulas; where relevant I will speak of a set *qua datum* or *qua dual datum*, and then this can be understood in terms of some such differentiating definitions as the ones above.

We use $\delta$, $\varepsilon$, … as variables for dual data.

A dual datum is *elementary* if all its formulas are literals. An elementary dual datum is *purified* if no component (element) includes a variable together with its negation. A purified elementary dual datum is a *canonical elementary dual datum* if all its components have the same vocabulary.

We will usually represent dual data using semi-colons to separate the components, and a backslash (\) to separate the formulas inside each component. – This notation is designed to be both *similar* to our notation for data (given the similarity of the ideas involved) and yet *not identical* so that data and dual data are almost always easily visually distinguishable even in simple cases – e.g. p, q is a datum but p; q a dual datum.

We say that an interpretation σ *verifies* a dual datum δ (*qua* dual datum) if every component Γ of δ is s.t. σ verifies at least one formula in Γ. We then define in obvious way *valid* dual datum, *satisfiable* dual datum, etc. etc.

We use ⊤ for { } (*qua* dual datum), and ⊥ for {{ }} (*qua* dual datum). Clearly ⊤ is valid, and ⊥ unsatisfiable.

– We list now the rules of our Dual Data Calculus (DDC). These will be 'duals' to the respective rules of DC, where this 'dualization' here means changing: (1) the direction of the arrow; (2) ∧ to ∨ and vice versa; and (3) commas to backslashes and slashes to semi-colons.

<div align="center">

G-rules:

|  | Γ \ A ∨ B ; ρ |  | Γ \ A ∧ B ; ρ |  | Γ \ ¬¬A ; ρ |
|---|---|---|---|---|---|
| (↕∨) | ↕ | (↕∧) | ↕ | (↕¬¬) | ↕ |
|  | Γ \ A \ B ; ρ |  | Γ \ A ; Γ \ B ; ρ |  | Γ \ A ; ρ |

|  | Γ \ ¬(A ∧ B) ; ρ |  | Γ \ ¬(A ∨ B) ; ρ |
|---|---|---|---|
| (↕¬∧) | ↕ | (↕¬∨) | ↕ |
|  | Γ \ ¬A \ ¬B ; ρ |  | Γ \ ¬A ; Γ \ ¬B ; ρ |

N-rules:

|  | Γ \ A \ ¬A ; ρ |  | Γ \ A ; Γ \ ¬A ; ρ |
|---|---|---|---|
| (EM) | ↑ | (NC) | ↓ |
|  | ρ |  | Γ ; ρ |

S-rules:

|  | Γ \ Δ ; ρ |  | Γ ; ρ |
|---|---|---|---|
| (Exp) | ↑ | (Del) | ↓ |
|  | Γ ; ρ |  | ρ |

</div>

We can then define notions of deduction, analytic deduction, synthetic deduction, etc. etc. – all as in DC.

We define the *dual*, d(α), of a datum α as the dual datum resulting from α by dualization of the formulas involved (which here means changing ∧ to ∨ and vice versa) and replacement of commas and slashes by respectively backslashes and semi-colons. (More precisely, this defines the dual of a datum-*expression*. The dual of a datum itself is defined in corresponding way, the exact formulation depending on which exact definition of dual datum we are using.) And similarly for the notion of the *dual* d(δ) of a dual datum δ. We have then the following basic connection between DC and DDC:

*Proposition.* α ⊢_DC β iff d(β) ⊢_DDC d(α); and δ ⊢_DDC ε iff d(ε) ⊢_DC d(δ).

*Proof.* This is immediate from the connection indicated above between the primitive rules of DC and DDC. Thus e.g. given a deduction of β from α in DC, turn it upside down and dualize each datum and the result is a DDC deduction of d(α) from d(β). (Incidentally, note that the first deduction is normal iff the second is.) □

– We come now finally to the *Hybrid Calculus* (HC).

A *deduction* (in HC) of dual datum ε from datum α is a sequence ⟨α … β, δ … ε⟩ where: ⟨α … β⟩ is a DC-deduction; ⟨δ … ε⟩ is a DDC-deduction; and δ is obtainable from β by the

*Criss-Crossing Rule* (CrCr): Infer a dual datum from a datum if, for any component Γ of the datum and any component Δ of the dual datum, there exists a formula belonging to both Γ and Δ.

– This rule is clearly sound (but not equivalential). We classify it as an S-rule, since it does not involve any logical constants. But we count it also at the same time as the single member of a new category of *Transition*-rules or T-rules. Also we call it neither analytic nor synthetic but 'neutral'.

We say that an HC-deduction ⟨α … β, δ … ε⟩ is *normal* if ⟨α … β⟩ is a purely analytic DC-deduction and ⟨δ … ε⟩ a purely synthetic DDC-deduction.

# 4. The Consecution Calculus (CC)

A canonical (so analytico-synthetic) deduction of dual datum δ from datum α in HC can naturally be 'reconceptualized' as a purely synthetic proof of a 'consecution' α ⇒ δ. The beginning of the proof corresponds essentially to the transition point of the deduction, and the subsequent steps in the proof correspond to going either down in the deduction, progressively synthesizing towards δ, or up in the deduction, progressively 'dis-analyzing' towards α.

The resulting 'Consecution Calculus' has obvious similarities with the familiar Sequent Calculus. (For the benefit of readers to whom the Sequent Calculus is *not* very familiar we have included a brief note at the end of the present section describing a standard version of that calculus.) (The term 'consecution' was proposed in Anderson & Belnap 1975 § 7.2 pp. 51–52 as a replacement for 'sequent'; but the proposal did not find many followers. So I am now 'recycling' the term, using it for a related but different notion.) Indeed it would have been possible (though less convenient, especially for later [quantificational and modal] extensions of the calculus) to formulate essentially the same calculus as a (linear) calculus of (finite) sequent-classes of a special kind which we might call 'uniform' sequent-classes – i.e. satisfying the condition that, whenever both ⟨Γ, Δ⟩ and ⟨Γ′, Δ′⟩ are in the class, so is ⟨Γ, Δ′⟩ (or equivalently, so are ⟨Γ, Δ′⟩ and ⟨Γ′, Δ⟩).

Consecution $=_{df}$ pair ⟨α, δ⟩ (usually written in the present context as α ⇒ δ) where α is a datum and δ a dual datum.

So we have now the further 'structural device' ⇒, this time corresponding to material implication. Basic semantical notions are then defined in accordance with this correspondence. (So note that an interpretation σ verifies a consecution α ⇒ δ iff σ verifies the uniform sequent-class {⟨Γ, Δ⟩ : Γ ∈ α and Δ ∈ δ}.)

⊤ $=_{df}$ { } ⇒ { }. (I.e. the verum as consecution is the consecution from the falsum as datum to the verum as dual datum!) We can also define (though this will be less useful for our purposes here): ⊥ $=_{df}$ {{ }} ⇒ {{ }}.

We proceed to list the basic rules of the Consecution Calculus (CC). Again there are no 'axioms' and a *proof* will be a deduction from ⊤. However, here our purpose is *not* to be able to deduce a consecution from another whenever there is semantic implication ('deductive completeness') but only to be able to (synthetically) *prove* all valid consecutions.

*Rules of CC*:

<u>G-rules:</u>

$$(\wedge\text{-L}) \quad \frac{\Gamma, A, B \,/\, \rho \;\Rightarrow\; \delta}{\Gamma, A \wedge B \,/\, \rho \;\Rightarrow\; \delta} \qquad\qquad (\vee\text{-R}) \quad \frac{\alpha \;\Rightarrow\; \Gamma \setminus A \setminus B \,;\, \rho}{\alpha \;\Rightarrow\; \Gamma \setminus A \vee B \,;\, \rho}$$

– Similarly for ¬∨-L and ¬∧-R.

$$(\neg\neg\text{-L}) \quad \frac{\Gamma, A \,/\, \rho \;\Rightarrow\; \delta}{\Gamma, \neg\neg A \,/\, \rho \;\Rightarrow\; \delta} \qquad\qquad (\neg\neg\text{-R}) \quad \frac{\alpha \;\Rightarrow\; \Gamma \setminus A \,;\, \rho}{\alpha \;\Rightarrow\; \Gamma \setminus \neg\neg A \,;\, \rho}$$

$$(\vee\text{-L}) \quad \frac{\Gamma, A \,/\, \Gamma, B \,/\, \rho \;\Rightarrow\; \delta}{\Gamma, A \vee B \,/\, \rho \;\Rightarrow\; \delta} \qquad\qquad (\wedge\text{-R}) \quad \frac{\alpha \;\Rightarrow\; \Gamma \setminus A \,;\, \Gamma \setminus B \,;\, \rho}{\alpha \;\Rightarrow\; \Gamma \setminus A \wedge B \,;\, \rho}$$

– Similarly for ¬∧-L and ¬∨-R.

<u>N-rules:</u>

$$(\text{NC}) \quad \frac{\alpha \;\Rightarrow\; \delta}{\alpha \,/\, \Gamma, A, \neg A \;\Rightarrow\; \delta} \qquad\qquad (\text{EM}) \quad \frac{\alpha \;\Rightarrow\; \delta}{\alpha \;\Rightarrow\; \delta \,;\, \Gamma \setminus A \setminus \neg A}$$

<u>S-rules:</u>

$$\text{(CrCr-L)} \quad \frac{\alpha \Rightarrow \delta}{\alpha \,/\, \beta \Rightarrow \delta} \qquad \text{provided } \beta \text{ criss-crosses } \delta$$

$$\text{(CrCr-R)} \quad \frac{\alpha \Rightarrow \delta}{\alpha \Rightarrow \delta \,;\, \varepsilon} \qquad \text{provided } \alpha \text{ criss-crosses } \varepsilon$$

*Remarks.*  (1) All these rules have a synthetic character. So we have used the common line rather than our (generally more 'informative') arrow in the formulation of the rules; and the proofs will require no arrow annotation.

(2) Note that all rules are single-premiss. I.e. CC is still (like our earlier calculi – and unlike sequent calculus) a calculus of 'linear reasoning'.

(3) Now *all* rules are equivalential! As we have already said, our purpose with CC will be to *prove* valid consecutions, i.e. deduce them from $\top$. So of course no rule could *actually* effect weakening in such a context (as nothing is weaker than $\top$!); and so it is perhaps not very surprising that we can make do with rules which could *never* effect weakening, in whatever context.

(4) The following rules are of course correct (and have synthetic character):

$$\text{(Del)} \quad \frac{\Gamma \,/\, \rho \Rightarrow \delta}{\Gamma, \Delta \,/\, \rho \Rightarrow \delta} \qquad\qquad \text{(Exp)} \quad \frac{\alpha \Rightarrow \Gamma \,;\, \rho}{\alpha \Rightarrow \Gamma \setminus \Delta \,;\, \rho}$$

However, they are not derived rules of CC. This follows immediately from the facts that all CC rules are equivalential but neither Del nor Exp is equivalential. – Also the *inverses* of the respective CC rules are correct rules but not derivable (as can easily be shown from their *analytic* character). – None of this is any 'defect' of CC: as we have said, the purpose of CC is *proofs* of (valid) consecutions, not arbitrary *deductions* of consecutions from consecutions (which I consider beyond the call of natural duty, as of course the consecutions themselves already correspond to implications).

(5) Note that any consecution $\alpha \Rightarrow \delta$ with $\alpha$ criss-crossing $\delta$ is deducible from ⊤ by one application of CrCr-L (any datum criss-crosses { }!) and then one application of CrCr-R (or alternatively one of CrCr-R then one of CrCr-L). ⊣

*Note on the Sequent Calculus*.  As promised we describe here (for the benefit of non-expert readers) a standard version of sequent calculus. There is a single axiom-scheme $A \Rightarrow A$, and the rules:

S-rules:

(Weakening-Left)  From $\Gamma \Rightarrow \Delta$ to $\Gamma, \Gamma' \Rightarrow \Delta$.

(Weakening-Right)  From $\Gamma \Rightarrow \Delta$ to $\Gamma \Rightarrow \Delta, \Delta'$.

G-rules:

($\neg$-L)  From $\Gamma \Rightarrow A, \Delta$ to $\Gamma, \neg A \Rightarrow \Delta$.

($\neg$-R)  From $\Gamma, A \Rightarrow \Delta$ to $\Gamma \Rightarrow \neg A, \Delta$.

($\wedge$-L)  From $\Gamma, A, B \Rightarrow \Delta$ to $\Gamma, A \wedge B \Rightarrow \Delta$.

($\vee$-R)  From $\Gamma \Rightarrow A, B, \Delta$ to $\Gamma \Rightarrow A \vee B, \Delta$.

($\wedge$-R)  From $\Gamma \Rightarrow A, \Delta$ and $\Gamma \Rightarrow B, \Delta$ to $\Gamma \Rightarrow A \wedge B, \Delta$.

($\vee$-L)  From $\Gamma, A \Rightarrow \Delta$ and $\Gamma, B \Rightarrow \Delta$ to $\Gamma, A \vee B \Rightarrow \Delta$.

*Proposition* (Completeness).  Every valid sequent has a proof in this calculus.

*Proof*.  Note first that all G-rules are equivalential in the sense that, for every interpretation $\sigma$, $\sigma$ verifies the conclusion iff $\sigma$ verifies all the premisses. So in particular conclusion always implies premiss, and thus if conclusion is valid premiss is valid. – Now take arbitrary valid sequent. Start attempted construction of its proof bottom-up applying G-rules (in any order) until all 'top nodes' are *elementary* sequents (i.e. where all formulas are propositional variables). By the above observation all sequents in this tree must be valid – in particular the top nodes. But of course an *elementary* sequent can only be valid if it is of the form $\Gamma, p \Rightarrow p, \Delta$. So using the Weakening rules and the axiom-scheme A $\Rightarrow$ A we conclude the construction of the desired proof. $\square$

# 5. Completeness etc.

*Proposition* (Normal-refutation completeness of DC). If α is unsatisfiable then there exists a normal refutation of α in DC (i.e. normal DC-deduction of ⊥ from α) (even without use of the analytic S-rule Del).

*Proof.* Take unsatisfiable α and apply analytic G-rules until we reach elementary datum α′. Recall that the G-rules are equivalential – so α′ is equivalent to α, and thus α′ is unsatisfiable. So each component of α′ is unsatisfiable, which (since α′ is elementary) means that each component contains an atom together with its negation. So successive applications of NC (one for each component) will in the end yield { } = ⊥. □

*Remark.* Note that here, as well as in all the other completeness results in this section, we show not only that a deduction *exists* but also how it can be *mechanically constructed*.

*Proposition* (Normal-proof completeness of DDC). If δ is valid then there exists a normal proof of δ in DDC (even without use of the synthetic S-rule Exp).

*Proof.* Take valid δ and begin constructing its proof from the bottom up, using the relevant synthetic G-rules. We then reach valid elementary dual datum δ′; so each of its components contains an atom together with its negation; so from { } = ⊤ we get δ′ by successive applications of EM. □

(Alternatively, this Proposition can be derived from the preceding one [or vice versa] via the 'basic connection' between DC and DDC pointed out in earlier section.)

*Proposition* (Normal-deduction completeness of HC). If α ⊨ δ then there exists a normal HC-deduction of δ from α (even without use of S-rules other than CrCr).

*Proof.* Begin by reducing (top-down) α to elementary α′ via ↓ G-rules, and (bottom-up) δ to elementary δ′ via ↑ G-rules. Then (if needed) purify α′ by NC into a purified elementary datum α*, and purify δ′ by EM into purified elementary dual datum δ*. All these rules are equivalential; so we have α* ⊨ δ*. But then α* must 'criss-cross' δ*. (For suppose not, i.e. there are components Γ of α* and Δ of δ* s.t. Γ and Δ have no formula in common. But then any interpretation verifying all the literals in Γ and falsifying all the literals in Δ will verify α* but falsify δ*. [And such an interpretation exists, given the purification and the assumption of no overlap.]) So the CrCr rule yields δ* from α*. □

*Proposition* (Proof-completeness of CC). If ⊨ α ⟹ δ then ⊢$_{CC}$ α ⟹ δ.

*Proof.* Suppose ⊨ α ⟹ δ, i.e. α ⊨ δ. Then by the proof of the preceding Proposition there is a 'canonical' HC-deduction of δ from α. From this a CC proof of α ⟹ δ can be constructed in obvious way: from ⊤ (= { } ⟹ { }) by CrCr-L and CrCr-R we deduce the consecution α* ⟹ δ* corresponding to the transition step in the HC-deduction; and then we go 'upwards' from α* towards α, and 'downwards' from δ* towards δ, using the obvious correspondence between: the ↓ G-rules of DC and the left G-rules of CC; NC of DC and NC of CC; the ↑ G-rules of DDC and the right G-rules of CC; and EM of DDC and EM of CC. (Or alternatively we can of course give a direct argument, without explicit reference to HC, but still corresponding to the construction of the canonical HC-deduction.) □

Unlike HC, DC is *not* normal-deduction complete. (Nor is DDC. The whole discussion below has dual version for DDC, which we leave tacit to avoid tediousness.) For example: p, q / ¬p / ¬q is not normally deducible from ⊤ (i.e. synthetically deducible from ⊤). For a supposed synthetic deduction could not use G-rules (since p, q / ¬p / ¬q is elementary) and so would have to use only EM and/or Exp. But the last rule used in a supposed deduction could not be Exp since no 'proper subdatum' of p, q / ¬p / ¬q is valid, nor (as is 'physically' clear) EM. – Another example is: p, r / q, ¬r is not normally deducible from p, q.

Later we will formulate variants of the primitive rules of DC which may give normal-deduction completeness. But now we define a slightly relaxed notion of *quasi-normality* and show a quasi-normal--deduction completeness theorem for DC. (Since the variants are more complex than our original rules, and quasi-normality seems to be normality enough, I am inclined to prefer this basic form of DC over the possible variants.)

Let us say that a deduction (in DC) is *quasi-analytic* (resp., *quasi-synthetic*) if all its steps are licensed either by analytic (resp., synthetic) rules or by S-rules – thus we allow also use of Exp (resp., Del). Then we say that a deduction (in DC) is *quasi-normal* if it consists of a quasi-analytic part followed by a quasi-synthetic part. Thus a quasi-normal deduction is like a normal deduction except that applications of Exp are allowed in the otherwise analytic part and applications of Del are allowed in the otherwise synthetic part. We say that β is quasi-normally deducible from α, or α ⊢$^{QN}$ β, if there exists a quasi-normal deduction of β from α. And similarly for quasi-analytic and quasi-synthetic deducibility.

*Proposition* (Quasi-normal--deduction completeness of DC). If α ⊨ β then α ⊢$^{QN}$ β (indeed even without use of Exp in the analytic part).

*Proof*. We start construction of the deduction from the 'extremities', reducing at the top α to an elementary datum α′ by ↓ G-rules, then a purified elementary datum α* by NC; and at the bottom β to elementary datum β′ by ↑ G-rules, then a purified elementary datum β* by Exp, then the canonical elementary datum β$^c$ equivalent to β / β′ / β* and with L(β$^c$) = L(β*), by Del. (Note that so far this is all equivalential – the applications of Exp introduce only unsatisfiable components, and the applications of Del lead to β* from the *ex hypothesi* equivalent canonical elementary datum β$^c$. Note also that the applications of Del in question will constitute the only 'non-normal part' in the whole deduction.) Now we move from α*, by Del, to the datum γ obtained from α* by deleting all literals whose variable ∉ L(β$^c$). (Thus L(γ) = (L(α*) ∩ L(β*)) ⊆ (L(α) ∩ L(β)).) We claim that γ ⊨ β, and will use this claim in what follows; but let us postpone its justification until after we have finished constructing our deduction. But note here that of course the passage from α* to γ need not be equivalential. – Now, transitioning now to the 'quasi-synthetic' part of the deduction, we move from γ, by EM, to the equivalent *canonical* elementary datum γ$^c$ with L(γ$^c$) = (L(γ) ∪ L(β$^c$)) = L(β$^c$). Now since γ ⊨ β ≃ β$^c$ and γ$^c$ ≃ γ, we have γ$^c$ ⊨ β$^c$. But these two are *canonical* elementary data, and with same vocabulary; so in fact γ$^c$ is a 'subdatum' (subset) of β$^c$; so we can move from γ$^c$ to β$^c$ by Exp, thus completing our deduction. (This step too need not be equivalential.)

– We give now the pending justification of our claim that γ ⊨ β. We know that α* ⊨ β, or as we may write illustratively α*(p, q) ⊨ β(q, r). Then also (by Substitution) α*(¬p, q) ⊨ β(q, r); and so

$$α*(p, q) / α*(¬p, q) \ ⊨ \ β(q, r).$$

So in the l.h.s. datum for each component $\Gamma$, $\pm p$ there will also be a component $\Gamma$, $\neg\pm p$. But then this part $\Gamma$, $\pm p$ / $\Gamma$, $\neg\pm p$ is equivalent to just $\Gamma$. (In general instead of just p we may have several variables $p_1$ … $p_n$ and then we take the substitutions with all combinations of assertion and denial $\pm p_1$ … $\pm p_n$ and then eliminate all $p_i$ as illustrated above.) So the result of these simplifications, which is precisely the datum $\gamma$, also implies $\beta(q, r)$. □

  *Proposition* (Directional Interpolation).  If $\alpha \vDash \beta$ then there exists $\gamma$ s.t.: $L(\gamma) \subseteq$ $L(\alpha) \cap L(\beta)$; $\gamma$ is analytically deducible from $\alpha$; and $\beta$ is quasi-synthetically deducible from $\gamma$.

  *Proof*.  The $\gamma$ in the preceding proof will be such a 'directional interpolant'.  □

# 6. Miscellaneous topics

(1) *Modularity Lemmas*.

*Proposition*. *Inner Modularity Lemma for DC*: If $\Gamma \vdash \Delta$ (i.e. $\{\Gamma\} \vdash \{\Delta\}$) then $\Gamma$, $\Sigma \vdash \Delta$, $\Sigma$.

*Outer Modularity Lemma for DC*: If $\alpha \vdash \beta$ then $\alpha / \rho \vdash \beta / \rho$.

*Proof*. Immediate from the completeness (and soundness) of DC. (E.g. in Outer case: $\alpha \vdash \beta$ implies $\alpha \vDash \beta$, which implies $\alpha / \rho \vDash \beta / \rho$, which implies $\alpha / \rho \vdash \beta / \rho$.) In the case of Outer Modularity it is also clear from the basic rules that, given a deduction of $\beta$ from $\alpha$, simply affixing '/ $\rho$' to each line will yield a deduction of $\beta / \rho$ from $\alpha / \rho$. (Presumably for Inner Modularity too there is some simple syntactic argument, though I don't myself immediately see it.) □

(2) *Proofs, refutations and deductions*.

Using the Modularity Lemmas we can easily pass from deductions to proofs or refutations, and vice versa:

*Proposition*. (i) $\top \vdash A \rightarrow B$ iff $A \vdash B$.

(ii) $A, \neg B \vdash \bot$ iff $A \vdash B$.

*Proof*. <u>(i) $\Rightarrow$</u>. Suppose hypothesis. Note that 'the datum A' is officially $\{\{A\}\} = \{\{A\} \cup \{\ \}\} =$ 'the datum A, $\top$'. We have then a deduction:

A

…                              (by hypothesis and Inner Modularity)

A, $\neg$A $\vee$ B            ↓ $\vee$

A, $\neg$A / A, B            ↓ NC

A, B                              ↓ Del

B.

(i) ⇐.  Suppose hypothesis. Then:

⊤  (= {{ }})   ↑ EM

¬A / A

…               (by hypothesis and Outer Modularity)

¬A / B        ↑ ∨

¬A ∨ B.

(ii) ⇒.  Suppose hypothesis. Then:

A                       ↑ EM

A, B / A, ¬ B

…                       (by hypothesis and Outer Modularity)

A, B                   ↓ Del

B.

(ii) ⇐.  Suppose hypothesis. Then:

A, ¬B

…               (by hypothesis and Inner Modularity)

B, ¬B         ↓ NC

⊥  (= { }).    □

A more general form of (i) is the following form of 'Deduction Theorem':

*Proposition.*  Γ ⊢ A → B iff Γ, A ⊢ B.

*Proof.*  (⇒).  Suppose hypothesis. Then:

Γ, A

…                          (by hypothesis and Inner Modularity)

¬A ∨ B, A              ↓∨

A, ¬A / A, B           ↓NC

A, B                   ↓Del

B.

(⇐). Suppose hypothesis. Then:

Γ                      ↑EM

Γ, A / Γ, ¬A

…                      (by hypothesis and Outer Modularity)

B / Γ, ¬A              ↓Del

B / ¬A                 ↑∨

¬A ∨ B.                □

The connections above between proofs, refutations and deductions permit us to easily pass from refutation-completeness or proof-completeness to deductive completeness (or vice versa, though that is less useful since it is refutation-completeness or proof-completeness that are easier to show than deductive completeness). Suppose e.g. we know that DC is refutation-complete, and know the Outer Modularity Lemma and hence (ii)⇒ above. To show deductive completeness of DC we can then argue as follows. Suppose α ⊨ β. Then for A, B the characteristic formulas of respectively α, β, we have A ⊨ B, i.e. {A, ¬B} is unsatisfiable. (For cases involving data which lack characteristic formulas, viz. ⊤ and ⊥, a simple supplementary argument can be given.) So by refutation-completeness A, ¬B ⊢ ⊥. So by (ii)⇒ we have A ⊢ B. But of course α ⊢ A (by ↑∧, ↑∨) and B ⊢ β (by ↓∨, ↓∧); and so α ⊢ β.

(3) *Derived rules*.

There is of course a natural notion of 'derived rule' of DC (or DDC etc.). E.g. for MP (Modus Ponens) we have the derivation:

A, ¬A ∨ B          ↓ ∨

A, ¬A / A, B          ↓ NC

A, B          ↓ Del

B.

It is interesting that here rules of all three kinds (a G-rule, an N-rule, an S-rule) are used. Also, note that they are all analytic rules; so MP is in the obvious sense of the term an *analytic derived rule* of DC.

∨-Introductions provide examples of *synthetic* derived rules of DC:

A          ↑ Exp

A / B          ↑ ∨

A ∨ B.

This is ∨-Introduction from the Left Disjunct; and ∨-Introduction from the Right Disjunct is derived similarly.

However, ∨-Introduction from Both Disjuncts (i.e.: infer A ∨ B from the datum A, B), although of course a derived rule, is *not* purely synthetic (nor of course purely analytic). The natural derivation is this (or the same but with B in the second line):

A, B          ↓ Del

A          ↑ Exp

A / B          ↑ ∨

A ∨ B.

And it is easily seen that there can be no synthetic derivation: for A ∨ B can only come synthetically (immediately) from A / B, and this can only come synthetically from A or from B, neither of which can come synthetically from A, B.

– All such rules have a more general form as in (for MP): From Γ, A, ¬A ∨ B / ρ, infer Γ, B / ρ. That these are also derived rules follows from the Modularity Lemmas.

*Question*.  Are all correct structural rules (in natural definition of the term) derivable from Del and Exp?

(4)  *Independence of rules*.

It seems clear that the basic rules of DC are independent in the following 'direction-sensitive' sense: no basic analytic rule is a derived *analytic* rule of the calculus given by the other basic DC rules; and same for synthetic.

Independence in the stronger, 'direction-insensitive' sense, does *not* hold, as is seen from the following cases. (I don't know whether there are other cases besides these six.)

(i) <u>Derivation of ↑¬¬ from EM, NC and Del:</u>

Γ, A / ρ                         ↑ EM

Γ, A, ¬A / Γ, A, ¬¬A / ρ    ↓ NC,  ↓ Del

Γ, ¬¬A / ρ.

(ii) <u>Derivation of ↓¬¬ from EM, NC, Del:</u>

Γ, ¬¬A / ρ                       ↑ EM

Γ, ¬¬A, A / Γ, ¬¬A, ¬A / ρ     ↓ NC,  ↓ Del

Γ, A / ρ.

(iii) <u>Derivation of ↑¬∨ from EM, NC, Del, ↓∨:</u>

Γ, ¬A, ¬B / ρ         ↑ EM

Γ, ¬A, ¬B, A ∨ B / Γ, ¬A, ¬B, ¬(A ∨ B) / ρ       ↓ ∨

Γ, ¬A, ¬B, A / Γ, ¬A, ¬B, B / Γ, ¬A, ¬B, ¬(A ∨ B) / ρ          ↓ NC

Γ, ¬A, ¬B, ¬(A ∨ B) / ρ             ↓ Del

Γ, ¬(A ∨ B) / ρ.

(iv) <u>Derivation of ↓¬∨ from EM, NC, Del, ↑∨:</u>

Γ, ¬(A ∨ B) / ρ             ↑ EM (×3)

Γ, ¬(A ∨ B), A, B / Γ, ¬(A ∨ B), A, ¬B / Γ, ¬(A ∨ B), ¬A, B / Γ, ¬(A ∨ B), ¬A, ¬B / ρ         ↓ Del (×3)

Γ, ¬(A ∨ B), A / Γ, ¬(A ∨ B), B / Γ, ¬(A ∨ B), ¬A, ¬B / ρ       ↑ ∨

Γ, ¬(A ∨ B), A ∨ B / Γ, ¬(A ∨ B), ¬A, ¬B / ρ                    ↓ NC,  ↓ Del

Γ, ¬A, ¬B / ρ.

(v) <u>Derivation of ↑¬∧ from EM, NC, Del, ↓∧:</u>

Γ, ¬A / Γ, ¬B / ρ             ↑ EM (×2)

Γ, ¬A, A ∧ B / Γ, ¬A, ¬(A ∧ B) / Γ, ¬B, A ∧ B / Γ, ¬B, ¬(A ∧ B) / ρ

                                          ↓ ∧ (×2),  ↓ NC (×2)

Γ, ¬A, ¬(A ∧ B) / Γ, ¬B, ¬(A ∧ B) / ρ       ↓ Del (×2)

Γ, ¬(A ∧ B) / ρ.

(vi) <u>Derivation of ↓¬∧ from EM, NC, Del, ↑∧:</u>

Γ, ¬(A ∧ B) / ρ             ↑ EM

Γ, ¬(A ∧ B), A / Γ, ¬(A ∧ B), ¬A / ρ       ↑ EM

Γ, ¬(A ∧ B), A, B / Γ, ¬(A ∧ B), A, ¬B / Γ, ¬(A ∧ B), ¬A / ρ

                                ↑ ∧,  ↓ NC,  ↓ Del (×2)

Γ, ¬A / Γ, ¬B / ρ.

Given these derivations (and the completeness of DC), we see that the following calculus is already complete: ↕∧, ↕∨, NC, EM, Del, Exp. Note that the ↓∧ and ↑∧ rules are 'determinative' of conjunction in the sense that conjunction is the only truth-function which (when assigned as 'meaning' to the symbol '∧') makes both these rules correct; and that likewise ↓∨ and ↑∨ are 'determinative' of disjunction, and NC and EM determinative of negation. So it is not so surprising that the calculus just described should be complete. – Note also that in this calculus (unlike in DC) every rule is 'pure' in the sense of containing at most one occurrence-of-connective.

(It is worth noting here also that any *application* of ↑Exp to a *non-empty* datum can be replaced by applications of ↑EM and ↓Del, on the pattern of the following example:

Γ                                                                    ↑ EM (×3)

Γ, A, B / Γ, A, ¬B / Γ, ¬A, B / Γ, ¬A, ¬B          ↓ Del (×4)

Γ / A, B.

But this does require at least one component [for the first application of EM]; and it is not here a case of a single derived rule but rather [generalizing the above example in the obvious way] a scheme, which does however cover all cases of expansion from non-empty datum.)

(5) *Invertible deductions*.

Consider e.g. the following deduction:

p → q   i.e. ¬p ∨ q           ↓ ∨

¬p / q                         ↑ ¬¬

¬¬q / ¬p                       ↑ ∨

¬¬q ∨ ¬p     i.e. ¬q → ¬p.

Like any deduction using only G-rules, this is what we may call an *invertible deduction* – i.e. a deduction s.t. the inverted ('upside-down') sequence of data is also a deduction. In the case of a 'G-deduction' (i.e. deduction using only G-rules) as in the example above, the steps in the inverted deduction are of course justified by the same G-rules as before only now in the opposite direction.

Let us temporarily modify the primitive rules of DC by replacing ↓Del with a strengthened ↓Del$^+$ where simultaneous deletions in any number of different components are allowed. (I.e.: From Γ, Γ′ / Δ, Δ′ / …, infer Γ / Δ / …) Then also all 'GN-deductions' (i.e. deductions using only G-rules and N-rules [not necessarily both]) are invertible deductions: use of a G-rule becomes use of the same G-rule in the opposite direction; use of NC becomes use of Exp; and use of EM becomes use of Del$^+$. (Thus note that of course the inverse deduction need not be itself a GN-deduction.)

Obviously an invertible deduction, and in particular a G-deduction or a GN-deduction, can only connect *equivalent* data. However, equivalent data may easily not be connectable by a G-deduction, and may even not be connectable by a GN-deduction, and *even* not connectable by any invertible deduction at all (e.g. p / ¬p and p / ¬p / q do not seem to be connectable by invertible deduction).

We will see in Chapter 4 that existence of G-deduction of β from α (α ⊢$_{STDC}$ β) has a natural semantic equivalent (α ≃$_{ST}$ β).

*Questions.* Try to find some interesting necessary and/or sufficient conditions for (for equivalent data α and β):

(i) existence of GN-deduction of β from α;

(ii) existence of invertible deduction of β from α.

(6) *Equivalential DC.*

Let EqDC be the calculus where we omit the non-equivalential rules Del and Exp of DC, and add the inverses ↑NC and ↓EM of ↓NC and ↑EM. I.e. the primitive rules of EqDC are the equivalential G-rules and the equivalential (now 'primitively two-directional') N-rules ↕NC and ↕EM. – This calculus EqDC is (correct and) complete for 'equivalential classical propositional logic': i.e. for any data α, β: If α ≃ β then α ⊢$_{EqDC}$ β (and, clearly, if α ⊢$_{EqDC}$ β then α ≃ β). For a suitable deduction can be constructed as follows. We first reduce α to an equivalent purified elementary datum α′ using the ↓ G-rules and ↓NC; then turn α′ into the equivalent canonical elementary datum α* in L(α) ∪ L(β) by ↑EM. Then, bottom up, we reduce β to an equivalent purified elementary datum β′ using the ↑ G-rules and ↑NC; and then turn β′ into the equivalent canonical elementary datum β* in L(α) ∪ L(β) by ↓EM. Since α ≃ β and α ≃ α*, β ≃ β*, also α* ≃ β*; which, since α*, β* are *canonical* elementary data in the same variables, means that in fact α* = β*, and so the deduction is complete.

Note that since in this calculus all rules are in the '↕ style', now *all* deductions are invertible!

We see also that, in (ordinary 'implicational') DC itself, if we 'primitivize' the inverses of ↓NC and ↑EM, we may omit Del. For the same construction above, now under the hypothesis merely that α ⊨ β, will give canonical elementary data α*, β* with α* ⊆ β*, and so Exp is enough to conclude the deduction. Indeed it is enough to 'primitivize' ↓EM, since ↑NC is of course derivable from Exp. I.e.: in DC we might (without changing the extension of the deducibility relation) replace ↓Del by ↓EM.

(7) *Alternative sets of (truth-functional) connectives.*

Given any set (even infinite set) of truth-functional connectives including negation, it is easy to adapt DC to give a calculus suitable to those connectives. (The set need not be 'functionally complete'. The calculus will be deductively complete for the data with the connectives in question.) E.g. take {¬, ↔}. The calculus then has the N- and S-rules as before, and ↕¬¬ as before, and:

$$(↕↔) \quad \frac{Γ, A ↔ B / ρ}{↕} \quad\quad\quad (↕¬↔) \quad \frac{Γ, ¬(A ↔ B) / ρ}{↕}$$
$$Γ, A, B / Γ, ¬A, ¬B / ρ \quad\quad\quad Γ, A, ¬B / Γ, ¬A, B / ρ$$

In general, for a truth-functional connective f, ↕f is in terms of the 'truth-conditions' of f(A, B, …), and ↕¬f in terms of its 'falsity-conditions'.

(DC itself can be reformulated in these terms, with the full three 'truth-conditions' in ↕∨ and ↕¬∧. This seems to me less natural than DC; but it is, as is readily seen, deductively equivalent.)

This procedure assumes the *presence of negation*. But we can reproduce it even in the absence of negation if we use 'signed formulas' instead of ordinary, 'unsigned' formulas. (See the reformulation of DC in terms of signed formulas described in the next section, which adapts immediately to arbitrary set of basic truth-functional connectives.)

# 7. Variants of the main calculi

(1) *EM\* in DC, NC\* in DDC*.

Let DC\* be DC with EM replaced by:

$$\Gamma, A \, / \, \Delta, \neg A \, / \, \rho$$
$$(EM^*) \qquad \uparrow$$
$$\Gamma, \Delta \, / \, \rho$$

Note that EM\*, unlike EM, is non-equivalential. EM\* is of course a derived rule of DC (use EM, then Del), but not (qua derived rule of DC) purely synthetic. Conversely, EM is a derived purely synthetic rule of DC\*: it is simply the special case of EM\* where $\Gamma = \Delta$.

The obvious counterexamples to full normal-form deduction completeness of DC are readily dealt with by EM\*. *Question*: Does full normal-form deduction completeness hold for DC\*?

– Similarly, DDC\* is DDC with NC replaced by:

$$\Gamma \setminus A \, ; \, \Delta \setminus \neg A \, ; \, \rho$$
$$(NC^*) \qquad \downarrow$$
$$\Gamma \setminus \Delta \, ; \, \rho$$

Again, NC\* is non-equivalential, etc.

It is also worth noting here the close similarity with the more familiar case of the following two forms of cut-rule in sequent calculus:

$$(Cut) \quad \frac{\Gamma, A \Rightarrow \Delta \quad \Gamma \Rightarrow A, \Delta}{\Gamma \Rightarrow \Delta} \qquad\qquad (Cut^*) \quad \frac{\Gamma, A \Rightarrow \Delta \quad \Gamma' \Rightarrow A, \Delta'}{\Gamma, \Gamma' \Rightarrow \Delta, \Delta'}$$

Indeed given the obvious translation of sequents or finite sets thereof (relevantly here set of two premisses of rule) into dual data (with e.g. p, q ⇒ r, s translated as ¬p \ ¬q \ r \ s), Cut and Cut\* became practically the same as NC and NC\*.

NC\* is also very close to the 'Resolution Rule' in the eponymous refutation method (viz.: From disjunctions of literals A ∨ p and B ∨ ¬p, infer A ∨ B [including inference of ⊥ from p and ¬p as limiting case]).

Indeed this connection with Resolution suggests a way of trying to answer the question above concerning normal-form deduction completeness: (i) proof of completeness of the Resolution Rule for refutation of unsatisfiable sets of 'clauses' (disjunctions of literals) should translate into proof of normal-form i.e. analytic refutation completeness for elementary dual data, hence for dual data in general since of course the ↓ G-rules reduce arbitrary dual datum to elementary dual datum; (ii) dually, we would have proof of normal-form i.e. synthetic proof completeness for elementary data; and finally (iii) an adaptation of the relevant construction may be possible to show, given α ⊨ β, how to synthetically deduce β (now not from ⊤ but) from the appropriate elementary datum obtainable after analyses from α.

(2) *Distribution rather than Criss-Crossing in HC*.

Our formulation of HC with the single Transition-Rule CrCr certainly allows very strikingly simple normal forms for deductions. However, in certain respects the CrCr rule is not very satisfactory, at least in the context of investigations aiming at unveiling the Ultimate Structure of Logical Deduction. Note e.g. that CrCr permits the inference from datum p, q to dual datum p – which would certainly seem like an analytic step –, but also from datum p to dual datum p \ q – which would certainly seem like a synthetic step. Thus this single atomic (basic) inferential connection would have both some analytic and some synthetic instances, contrary to the whole tendency of the rest of our schemes.

From this point of view the following alternative version of HC is preferable. We replace CrCr by the following single new Transition-Rule (also an S-rule):

(Distr) From datum α, infer dual datum whose components are the sets which choose one formula from each component of α.

(This is the correlate for data and dual data of the usual transformation by 'distribution' of a disjunction of conjunctions into an equivalent conjunction of disjunctions – e.g. of $(p \wedge q) \vee (r \wedge s \wedge u)$ into

$$(p \vee r) \wedge (p \vee s) \wedge (p \vee u) \wedge (q \vee r) \wedge (q \vee s) \wedge (q \vee u).)$$

Thus we simply 'rewrite' the datum as a dual datum. As before with CrCr, here too we call Distr neither analytic nor synthetic but 'neutral' – now with better justification for the term (for CrCr in view of the considerations above 'hybrid' might be more appropriate). It is also nice that this Distribution rule (unlike CrCr) is equivalential, so that the ordinary S-rules of DC and DDC remain the sole loci of weakenings.

The normal-form deduction completeness theorem for this new formulation of HC (in the sense: If $\alpha \vDash \delta$ then there exists a deduction of $\delta$ from $\alpha$ where all applications of analytic rules precede all applications of synthetic rules) can be shown as follows. Suppose $\alpha \vDash \delta$. As before, reduce $\alpha$ by $\downarrow$ G-rules and NC to elementary purified datum $\alpha'$, and $\delta$ (bottom-up) by $\uparrow$ G-rules and EM to elementary purified dual datum $\delta'$. Then $\alpha'$ must criss-cross $\delta'$. Now apply Distr to $\alpha'$; let us call the result of that $\varepsilon$. Since $\alpha'$ criss-crosses $\delta'$, each component of $\delta'$ must be an expansion (i.e. superset) of some component of $\varepsilon$ (i.e. of some choice-set of $\alpha'$). So to get from $\varepsilon$ to $\delta'$ it is enough first to omit by Del such components as may occur in $\varepsilon$ without having an expansion in $\delta'$, and then to prolong by Exp the remaining components to yield exactly the components of $\delta'$. (Note that since we are dealing with a set here, we can apply Exp several times to the same component to obtain different prolongations if necessary.)

Note however a slight impurity here. The deduction described is indeed a normal form deduction in the sense that all applications of analytic rules precede all applications of synthetic rules. But the synthetic rules do not always start being applied immediately after the application of the Transition-rule as one might expect (and as is the case in the normal forms for HC with CrCr), but we may need to use $\downarrow$Del immediately after the transition.

But I think the main consideration in favour of the formulation of HC with Criss-Crossing rather than Distribution is the following. The whole motivation for 'going hybrid', i.e. constructing a calculus involving both data and dual data, is that, just as there is a simple syntactic criterion for the unsatisfiability of an elementary datum, and another for the validity of an elementary dual datum, so also there is a simple syntactic criterion for implication from an elementary datum to an elementary dual datum – and what the CrCr rule does is precisely (with help from the presence of NC and EM) to enshrine that criterion, and thus permit strikingly simple normal form deductions in HC, just as we have strikingly simple normal form refutations in DC, and strikingly simple normal form proofs in DDC. – All this is lost in the formulation with Distr. One might just as well have a hybrid calculus starting with dual data and ending with data, with a Distribution rule now going the other way – from dual datum to equivalent distributed datum. Or even better not go hybrid at all and have a single type of formal objects (only data, or only dual data) and do the intermediate portion of deductions with Del, Exp etc. rather than with Del, Exp etc. *and* Distr!

(3) *'Atomic' versions of N-rules and S-rules*.

There would be no loss in the 'deductive power' of our calculi if we reformulated the N-rules in terms of an atomic formula p rather than an arbitrary formula A, or even the ordinary S-rules in terms of deletions of literals rather than arbitrary formulas and expansion by components consisting entirely of literals rather than arbitrary components, or the T-rules in terms of criss-crossing or distribution from an elementary datum to an elementary dual datum. This is clear from our proofs of the normal and quasi-normal deduction completeness theorems, where the deductions begin and end with G-rules and have a middle portion consisting entirely of manipulation of elementary data and/or elementary dual data by N-rules and S-rules alone, including a T-rule in case of HC.

This is an interesting fact, worth bearing in mind. But as for actually adopting the restricted 'atomic' versions of the rules instead of the 'generic' versions, I don't see any good reason for doing that ('because we can' is obviously a lousy reason); and on the other hand our 'generic' versions are in themselves simpler and more natural (why the artificial restrictions if they are not *needed*?), and often allow simpler and more natural deductions than what would be possible with the atomic versions.

(4) *No negation in primitive notation except for atoms*.

Thus ¬¬A is understood as an *abbreviation* (in the 'logical' sense!) of A, ¬(A ∧ B) of ¬A ∨ ¬B, and ¬(A ∨ B) of ¬A ∧ ¬B. The rules ↕¬¬, ↕¬∧, ↕¬∨ are then dropped, and NC and EM assume the 'atomic' version (now not 'because we can' but 'because we must'!). This is I guess 'elegant', but really rather artificial and obscures the whole conceptual situation. (Also, in practice the work saved in one place reappears in another, since we will have to rewrite natural formulas to conform to this artificial primitive notation.) At least we should say this if we think, as I do, that 'metaphysically speaking', any proposition can be negated (in a direct sense). One might alternatively take the recondite (but not entirely absurd) view that, 'metaphysically speaking', only atomic propositions can be negated (or have their 'polarity' 'inverted'), and 'negation of molecular propositions' should be understood as shorthand for the 'negation normal form' transform and so on – then of course the logical devices here *would* correspond to what one thinks is 'the whole conceptual situation'.

(5) *Minor variants of N- and S-rules depending on 'how much we can do at once'.*

We already saw before (in connection with 'invertible deductions') the variant Del$^+$ where we can make simultaneous deletions in multiple components. Likewise there is Exp$^+$ where we can add multiple components; and NC$^+$ where we can simultaneously drop multiple contradictory components; and EM$^+$ where we can 'bifurcate' multiple components.

In the opposite direction, there is Del$^-$ where we can delete just *one formula*; and NC$^-$ where we can drop just a component A, ¬A (rather than any component A, ¬A, Γ). To be thorough, there is even Del$^{-\,+}$ where we can make simultaneous deletions in multiple components but of one formula in each; and similarly for NC$^{-\,+}$.

In all such cases there is of course easy interderivability of the relevant rules. I have adopted Del, Exp, NC, EM as the 'official' rules rather than their 'plus' counterparts in the spirit of 'doing one thing (inferential step) at a time'. But one might say of course that Del$^-$ and maybe NC$^-$ are even more in *that* spirit … (It does seem to me quite unsatisfactory though to have first to delete say q from a component p, ¬p, q so as to only then be able to apply Non-Contradiction to drop the component.)

Note that 'doing one thing at a time' is particularly inefficient if we have an *infinite* number of things to do. So in *infinitary* logic we *need* the 'plus' versions. Maybe that is an indication that even in the context of finitary logic (as here) they should be adopted as basic rules.

(6) *Signed formulas*.

Interesting variants of our calculi are obtained by use of 'signed formulas' TA, FA, as in Smullyan's book. Then instead of $\updownarrow\wedge$, $\updownarrow\neg\wedge$, $\updownarrow\vee$, $\updownarrow\neg\vee$ we have $\updownarrow$T$\wedge$, $\updownarrow$F$\wedge$, $\updownarrow$T$\vee$, $\updownarrow$F$\vee$ (with the obvious formulations), and instead of $\updownarrow\neg\neg$ we have $\updownarrow$T$\neg$ and $\updownarrow$F$\neg$. Also we use the 'signs' T and F instead of negation in both EM and NC; so that the G-rules are now the only non-structural rules. – Here we use two further 'structural' devices (the 'signs'); but on the other hand the resulting calculus is slightly more symmetric, and all the 'introduction' and 'elimination' rules are 'pure' in the sense of involving just one occurrence-of-connective, and (as indicated in the preceding section) the calculus is now automatically generalizable to a *completely arbitrary* set of basic truth-functional connectives (thus including the cases where negation is absent). – I will not discuss here further motivations for using, or for not using, this device; on this see e.g. Rumfitt 2000, Smiley 1996, etc.

(7) *Multiple-Conclusion Calculi.*

DC may be regarded as resulting from the familiar method of analytic tableaux by two modifications:

(i) Instead of having just formulas as the 'formal objects' of our calculus, and deductions with a tree structure with some of the rules involving 'branching', we now have more complex formal objects (data) but simpler, purely linear deductions. (The complexity is as it were transferred from the structure of deductions to the structure of the 'formal objects'.)

(ii) Instead of a purely analytic method designed for refutation, we now have a richer analytico-synthetic method for deduction in general.

– It is also possible however to make just one (either one) of these two modifications. Making just modification (i) gives the 'purely analytic fragment' of DC, perfectly suitable as a general method for refutations. (As we have already mentioned, Smullyan himself considers this in the last chapter of his book.) – Making just modification (ii) yields a kind of analytico-synthetic tableau method, which we proceed to describe here.

Here we deal with just formulas (rather than data, or sequents, etc.), but allow deductions with 'branching' through rules with 'multiple conclusions'. Thus a deduction has tree form exactly as in tableaux. The rules are as follows (note that here there are no S-rules):

↓ G-rules:  From ¬¬A infer A; From A ∧ B infer A; From A ∧ B infer B; From A ∨ B branch with A on left and B on right (just as in tableaux); and similarly rules for ¬∨ and ¬∧.

↑ G-rules:  From A infer ¬¬A; From A infer A ∨ B; From B infer A ∨ B; From A and B infer A ∧ B; and similarly rules for ¬∧ and ¬∨.

N-rules:  (NC) From A and ¬A infer B; (EM) From B branch with A on left and ¬A on right.

A *deduction* of $\{B_1, …, B_k\}$ from $\{A_1, …, A_n\}$ ($n ≥ 1$, $k ≥ 1$) is a tableau starting with some of the formulas $A_1 … A_n$ and with every branch containing at least one of the formulas $B_1 … B_k$. Also a *proof* of A can be defined (*faute de mieux*) as a deduction of A from ¬A; and a *refutation* of A as a deduction of ¬A from A. – As is easily seen, there is a proof of A in this sense iff for every formula B there is a deduction of A from B; and there is a refutation of A iff for every formula B there is a deduction of B from A.

The deduction-completeness theorem for this calculus can be shown as follows. Suppose $\{A_1, …, A_n\} ⊨ \{B_1, …, B_k\}$, i.e. $A_1 ∧ … ∧ A_n ⊨ B_1 ∨ … ∨ B_k$. We first fully analyze the premisses $A_1 … A_n$, i.e. construct the 'exhaustive tableau' for them. Then to the branches with contradictory formulas we adjoin say $B_1$ by NC. Then in the remaining branches we apply EM so as to have in each such branch literals for each variable occurring in $\{B_1, …, B_k\}$. Next we synthesize $B_1 ∨ … ∨ B_k$ in each such branch using the ↑ G-rules. And finally we apply ↓∨ to this disjunction in each branch.

(The above calculus is given, under the name 'deduction trees', in Jeffrey 1981 § 2.5. [It does not seem to be present either in the earlier or the later editions of that book.] The calculus is also very close in spirit to Kneale's original multiple-conclusion calculus. [See Kneale 1956, and Kneale and Kneale 1984 Ch. 9 § 3 and note (3) on p. vi.])

The 'general format' of the above rules is: *From* a certain (finite) non-empty set of meta-formulas *to* a certain (finite) non-empty set of meta-formulas – or briefly,

(A) ⟨Non-empty Set, Non-empty Set⟩.

A more liberal format drops the restriction of non-emptiness from both sides:

(B) ⟨Set, Set⟩.

And a third format, more liberal still, permits 'multiply multiple conclusions':

(C) ⟨Set, Set-of-Sets⟩.

– In each case the natural concept of *deduction* for the calculus will have then a format matching that of the basic rules.

Format (B) permits reformulation of EM as a zero-premiss rule, and of NC as a zero-conclusion ('zero-branching') rule (i.e. we 'destroy' a branch). Thus a *proof* of formula A can now be more satisfactorily defined as a zero-premiss tableau where every branch contains A; and similarly a *refutation* of A as an empty tableau (not that there is more than one) constructible from the set of premisses {A}.

Format (C) is particularly convenient for treatment of connectives beyond just ¬, ∧, ∨. For example with ↔, the natural rule for analysis of A ↔ B is: From A ↔ B, branch into A, B and ¬A, ¬B. And similarly for analysis of ¬(A ↔ B). In general, for truth-functional connective f, and v, v′, … representing the 'truth-conditions' for the corresponding truth-function and u, u′, … its 'falsity-conditions', we have the natural rules:

From f(A̲), branch into v(A̲), v′(A̲), …

From ¬f(A̲), branch into u(A̲), u′(A̲), …

– The natural *synthetic* rules on the other hand are already in format (B) (indeed in format (A)):

From v(A̲) to f(A̲);

From v′(A̲) to f(A̲);  etc.

From u(<u>A</u>) to ¬f(<u>A</u>);

From u′(<u>A</u>) to ¬f(<u>A</u>);  etc.

Thus we have a natural formulation of a complete calculus for any set of truth-functional connectives including negation – or for *any* set of truth-functional connectives, if we use 'signed formulas'.

Note however that it is in fact possible to replace the above analytic format (C) rules by analytic format (B) rules:

From f(<u>A</u>), u(<u>A</u>), zero-branch;

From f(<u>A</u>), u′(<u>A</u>), zero-branch;  etc.

From ¬f(<u>A</u>), v(<u>A</u>), zero-branch;

From ¬f(<u>A</u>), v′(<u>A</u>), zero-branch;  etc.

(E.g.: From A ↔ B, A, ¬B, zero-branch.) Thus to obtain the effect of the original natural format (C) rules, we first apply EM to branch w.r.t. all truth-value combinations for the formulas <u>A</u>, and then destroy the appropriate branches using the new, format (B) rules. (Indeed essentially the same method can be used to eliminate arbitrary proper format (C) rules in favour of a set of zero-conclusion format (B) rules [assuming presence of EM].) – Still this procedure is of course less natural, and less purely analytic, than the original procedure with format (C) rules.

– Incidentally, note that in DC itself, and in its variants for other truth-functional connectives, there is a similar phenomenon. E.g. instead of the natural 'component-fissioning' ↓↔ rule, one might have two 'component-destroying' rules, and the effect of the natural rule would then be obtainable (though synthetico-analytically) through EM and the new rules.

CHAPTER 2

DIRECTIONAL DEDUCTION FOR CLASSICAL PREDICATE LOGIC

## 1. Basic notions

We extend here our methods of Directional Deduction to classical predicate logic.

We use a standard formulation of classical first-order predicate logic (with functional variables but without identity): the basic symbols (for construction of ordinary formulas) are the 'logical constants' $\neg$, $\wedge$, $\vee$, $\forall$, $\exists$, individual variables x, y, z, …, and for each $n \geq 0$ a denumerable stock of n-ary predicate variables and a denumerable stock of n-ary functional variables. *Formulas* are then defined in the usual way.

A *datum* now is an expression of the form

$$(E\underline{f}) \ (\underline{x}): \Gamma \ / \ \Delta \ / \ …$$

where $\Gamma$, $\Delta$, … are finite sets of formulas. Formally, this can be taken as an ordered triple consisting of: (i) a finite set of functional variables; (ii) a finite set of individual variables; and (iii) a finite set of finite sets of formulas.

In such a datum $(E\underline{f})$ may be called the *existential (functional) prefix*; $(\underline{x})$ the *universal (individual) prefix*; $(E\underline{f}) \ (\underline{x})$ as a whole the (full) *prefix*; and $\Gamma \ / \ \Delta \ / \ …$ the *matrix*. (Note that, differently from the familiar terminology for prenex formulas, here what we call the matrix need not be quantifier-free.) We say that a datum is *elementary* if all formulas in the matrix are either atomic formulas or negations of atomic formulas.

We consider data as literally *identical* if they so to speak 'differ only in the renaming of bound variables'. I.e. we define data as really equivalence-classes of expressions w.r.t. the equivalence-relation of obtainability via renamings of bound variables. Thus e.g. the datum $\forall xPx$ (i.e. the equivalence-class generated from $\langle \{ \ \}, \{ \ \}, \{\{\forall xPx\}\}\rangle$) is literally identical with the datum $\forall yPy$; and the datum (Ef) (x): Pf(x) is literally identical with the datum (Eg) (y): Pg(y).

*Remark.* A more precise definition of this equivalence-relation can be given as follows. We may speak of 'proto-data' in the sense of data 'prior to' this 'abstraction'. (I.e. a proto-datum is a triple consisting of (i) a finite set of functional variables, etc. [as in the initial formulation above].) Then we say that proto-datum X is *bound-variable similar* to proto-datum Y if either Y can be obtained from X by some renaming of bound variables (with the usual proviso that free variables should not be 'captured' by the renaming) or vice versa. – Note that this is *not* an equivalence-relation, as it is not transitive (it is clearly reflexive and symmetric): e.g. $\forall x_1 P x_1$, $\forall x_2 P x_2$, $\forall y Q y$ (i.e. $\langle \{ \ \}, \{ \ \}, \{\{\forall x_1 P x_1, \forall x_2 P x_2, \forall y Q y\}\}\rangle$) is bound-variable similar to $\forall x P x$, $\forall y Q y$, which in turn is bound-variable similar to $\forall x P x$, $\forall y_1 Q y_1$, $\forall y_2 Q y_2$, but the first of these proto-data is clearly *not* bound-variable similar to the third. – We say then finally that X is *bound-variable equivalent* to Y if X stands in the transitive closure of the bound-variable similarity relation to Y. This gives our desired equivalence-relation. – Intuitively, we may think of proto-data as having an associated 'archetype' in a more 'pure' notation without bound variables at all (e.g. we might use Quine's 'arcs'); and then proto-data X and Y are bound-variable equivalent iff their associated archetype is the same. E.g. in the previous example above the archetype is $\forall(P)$, $\forall(Q)$ (= $\forall(P)$, $\forall(P)$, $\forall(Q)$ = $\forall(P)$, $\forall(Q)$, $\forall(Q)$). ⊣

Semantical notions are defined in obvious way, with the existential prefix understood as corresponding to functional existential quantification, the universal prefix to individual universal quantification, and the matrix to a disjunction of conjunctions.

We say that a datum is *standard* if the following conditions are satisfied:

(1) all functional variables in the matrix appear also in the prefix;

(2) functional variables in the matrix appear only in application to individual *variables* (as opposed to functional terms); and

(3) there exists an ordering (total ordering) of the universal prefix which is s.t.: in every functional term in the matrix, the variables (as ordered in the term) constitute an initial segment (without repetitions) of the ordering.

– Our main interest here is to develop methods for deduction of formulas of Q (first-order predicate logic without identity and without functional variables) from formulas of Q. So we will be mainly concerned with *standard* data, which have the following basic connections with Q-formulas:

(i) For each Q-formula A, there exists an equivalent elementary standard datum. (Prenex A, then Skolemize that prenex formula, then put the matrix in DNF: the result is in effect an elementary standard datum.)

(ii) For each standard datum $\alpha$, there exists an equivalent Q-formula. (Just 'de-Skolemize' $\alpha$, or the corresponding second-order logic formula $\exists \underline{f} \forall \underline{x}(M)$: drive in universal quantifiers until the universal prefix coincides with the longest argument-sequence $x_1 \ldots x_n$ appearing in functional terms $f(x_1 \ldots x_n)$, $g(x_1 \ldots x_n)$, $\ldots$ in M; then drive in $\exists f \exists g \ldots$ across the (current) universal prefix as existentially quantified variables $\exists y \exists z \ldots$, replacing $f(x_1 \ldots x_n)$, $g(x_1 \ldots x_n)$, $\ldots$ in the matrix by y, z, $\ldots$; and continue like this until all functional variables have been eliminated.)

However, it proves convenient to permit non-standard data to appear in the course of deductions. (This is the case particularly with instantiations of universal variables.)

Note that *non*-standard data need not be 'first-orderizable'. A typical example is

$$(\text{E}f, \underline{g}) \ (x, z): R(x, f(x), z, g(z)),$$

which corresponds to the Henkin or branching quantification $\forall x \exists y$ 'and' $\forall z \exists u$ s.t. R(x, y, z, u) (with u 'depending' only on z), which as is well known is not first-orderizable.

– Now similarly a *dual datum* is an expression of the form

$$(\underline{f}) \ (\text{E}\underline{x}): \Gamma; \Delta; \ldots$$

Similar definitions apply as with data. In particular there is the notion of *standard* dual datum, defined as above for data (with (3) now referring of course to the existential individual prefix).

A *consecution* is now an expression

$$(\underline{f}) \ (\underline{g}) \ (\text{E}\underline{x}) \ (\text{E}\underline{y}): \alpha \Rightarrow \delta$$

where $\alpha$ is a datum-matrix and $\delta$ a dual-datum--matrix. We say that such a consecution is *standard* if both $(\text{E}\underline{f})(\underline{x})$: $\alpha$ is a standard datum and $(\underline{g})(\text{E}\underline{y})$: $\delta$ is a standard dual datum. We say also that the consecution is *regular* if it is standard in this sense and moreover the variables $\underline{f}$ do not occur in δ, nor the $\underline{g}$ in α, nor the $\underline{x}$ free in δ, nor the $\underline{y}$ free in α. In Quantificational CC we will be interested in the *proof* of regular consecutions (which corresponds to: deductions of standard dual data from standard data); but in the construction of such proofs it is convenient to permit occurrence of non-regular consecutions.

— Numerous definitions from propositional case transfer straightforwardly to quantificational case and so will not be given here explicitly but are to be assumed as implicitly present.

## 2. The calculi QDC, QDDC, QHC, QCC

As in propositional case, we give four calculi: one for data, one for dual data, one hybrid, and one for consecutions: we call these respectively Quantified Data Calculus (QDC), Quantified Dual Data Calculus (QDDC), Quantified Hybrid Calculus (QHC), and Quantified Consecution Calculus (QCC). Again the rules come in three groups, viz. G-rules, N-rules and S-rules; and in QHC one of the S-rules is a T-rule (transition rule). The G-rules include the obvious counterparts of the earlier propositional rules, plus rules in very similar style for analysis and synthesis w.r.t. $\forall$, $\neg\exists$, $\exists$, $\neg\forall$. The N-rules are the obvious counterparts of NC and EM. The S-rules include the obvious counterparts of the earlier weakening rules of ↓Deletion (for [in data case] commas) and ↑Expansion (for slashes), plus: two new weakening rules corresponding to the new 'structural' resources, viz. ↓Instantiation (for universal prefix) and ↑Generalization (for existential prefix); and two equivalential S-rules, viz. ↓EVE (Elimination of Vacuous Existential variables) and ↑IVU (Introduction of Vacuous Universal variables). And the T-rule for QHC will be an adaptation of the earlier Criss-Crossing rule.

*Rules of QDC*:

<u>G-rules:</u>

↕∧, ↕∨, ↕¬¬, ↕¬∧, ↕¬∨ as before. E.g. ↕∧ is now:

$$
\begin{array}{c}
(E\underline{f}) \ (\underline{x}): A \wedge B, \Gamma \ / \ \rho \\
\updownarrow \\
(E\underline{f}) \ (\underline{x}): A, B, \Gamma \ / \ \rho
\end{array}
$$

$$
(\updownarrow\forall) \qquad
\begin{array}{c}
(E\underline{f}) \ (\underline{x}): \forall yA, \Gamma \ / \ \rho \\
\updownarrow \\
(E\underline{f}) \ (\underline{x}, y): A, \Gamma \ / \ \rho
\end{array}
\qquad \text{provided y does not occur free in } \Gamma, \rho.
$$

$$
(\updownarrow\neg\exists) \qquad
\begin{array}{c}
(E\underline{f}) \ (\underline{x}): \neg\exists yA, \Gamma \ / \ \rho \\
\updownarrow \\
(E\underline{f}) \ (\underline{x}, y): \neg A, \Gamma \ / \ \rho
\end{array}
\qquad \text{provided y not free in } \Gamma, \rho.
$$

$(\updownarrow\exists)$     (E$\underline{f}$) ($\underline{x}$): $\exists yA(y)$, $\Gamma$ / $\rho$
              $\updownarrow$                    provided g not in A(y), $\Gamma$, $\rho$.
          (E$\underline{f}$, g) ($\underline{x}$): $A(g\underline{x})$, $\Gamma$ / $\rho$

$(\updownarrow\neg\forall)$     (E$\underline{f}$) ($\underline{x}$): $\neg\forall yA(y)$, $\Gamma$ / $\rho$
              $\updownarrow$                    provided g not in A(y), $\Gamma$, $\rho$.
          (E$\underline{f}$, g) ($\underline{x}$): $\neg A(g\underline{x})$, $\Gamma$ / $\rho$

N-rules:

(NC)     (E$\underline{f}$) ($\underline{x}$): $\Gamma$, A, $\neg A$ / $\rho$          (EM)     (E$\underline{f}$) ($\underline{x}$): $\Gamma$, A / $\Gamma$, $\neg A$ / $\rho$
              $\downarrow$                                              $\uparrow$
          (E$\underline{f}$) ($\underline{x}$): $\rho$                              (E$\underline{f}$) ($\underline{x}$): $\Gamma$ / $\rho$

S-rules:

(Del)     (E$\underline{f}$) ($\underline{x}$): $\Gamma$, $\Delta$ / $\rho$          (Exp)     (E$\underline{f}$) ($\underline{x}$): $\Gamma$ / $\rho$
              $\downarrow$                                              $\uparrow$
          (E$\underline{f}$) ($\underline{x}$): $\Gamma$ / $\rho$                          (E$\underline{f}$) ($\underline{x}$): $\rho$

(Inst)     (E$\underline{f}$) ($\underline{x}$, y): $\alpha(y)$          (Gen)     (E$\underline{f}$, g) ($\underline{x}$): $\alpha(g\underline{x})$
              $\downarrow$                                              $\uparrow$
          (E$\underline{f}$) ($\underline{x}$): $\alpha(t)$                            (E$\underline{f}$) ($\underline{x}$): $\alpha(t)$

(EVE)     (E$\underline{f}$, g) ($\underline{x}$): $\alpha$
              $\downarrow$                    provided g not in $\alpha$
          (E$\underline{f}$) ($\underline{x}$): $\alpha$

(IVU)     (E$\underline{f}$) ($\underline{x}$, $\underline{y}$): $\alpha$
              $\uparrow$                    provided $\underline{y}$ not free in $\alpha$
          (E$\underline{f}$) ($\underline{x}$): $\alpha$

(In Inst and Gen, t is an arbitrary term.)

Note that the inverses of EVE and IVU are readily derivable by respectively Gen and Inst.

– As in propositional case, the rules of QDDC are obtained by dualizing the rules of QDC. (Here we won't bother to write this all out explicitly, since there are a lot of rules now and the dualization is entirely straightforward.)

Now, in QHC, as before in propositional case the deductions start with data and end with dual data and there is a single T-rule:

*Criss-Crossing*: Infer dual datum (g̲) (Ey̲): δ from datum (Ef̲) (x̲): α, if there are instances $α_0$, $δ_0$ of respectively (x̲): α, (Ey̲): δ s.t. $α_0$ 'weakly criss-crosses' $δ_0$ – i.e. every pure component of $α_0$ has a formula in common with every pure component of $δ_0$.

(An 'instance' here means of course any expression resulting from the given one by dropping the prefix and uniformly substituting arbitrary terms for the indicated variables in the matrix. And a 'pure' component is one without pair of contradictory formulas A, ¬A.) This rule is clearly correct: – For suppose there are such instances $α_0$, $δ_0$ with $α_0$ weakly criss-crossing $δ_0$. Then of course (with some harmless *abus de langage*):

$$⊨ \ α_0 → δ_0.$$

And so $⊨$ (x̲): α .→. (Ey̲): δ,

whence $⊨$ (f̲) (g̲): (x̲): α .→. (Ey̲): δ,

whence $⊨$ (Ef̲) (x̲): α .→. (g̲) (Ey̲): δ.

We come finally to QCC. We use as basic rules consecution-versions of the QDC and QDDC derived rules of Generalized Instantiation and Generalized Generalization (see section after next). We use below $α_0$ to indicate a 'conjunction' (by 'distribution') of instances of (z̲): α, and $δ_0$ for a 'disjunction' of instances of (Ez̲): δ.

*Rules of QCC*:

G̲-rules: The left and right rules for ∧, ∨, ¬∧, ¬∨, ¬¬ as in CC; and the following plus the obvious analogues ¬∃-L, ¬∀-R, ¬∀-L, ¬∃-R.

$(\forall\text{-L})$ $\dfrac{(\underline{f})\ (\underline{g})\ (\text{E}\underline{x},\ z)\ (\text{E}\underline{y})\colon A(z),\ \Gamma\ /\ \rho\ \Rightarrow\ \delta}{(\underline{f})\ (\underline{g})\ (\text{E}\underline{x})\ (\text{E}\underline{y})\colon \forall z A(z),\ \Gamma\ /\ \rho\ \Rightarrow\ \delta}$ provided $z$ not free in $\Gamma$, $\rho,\ \delta$

$(\exists\text{-R})$ $\dfrac{(\underline{f})\ (\underline{g})\ (\text{E}\underline{x})\ (\text{E}\underline{y},\ z)\colon \alpha\ \Rightarrow\ A(z)\setminus\Gamma\ ;\ \rho}{(\underline{f})\ (\underline{g})\ (\text{E}\underline{x})\ (\text{E}\underline{y})\colon \alpha\ \Rightarrow\ \exists z A(z)\setminus\Gamma\ ;\ \rho}$ provided $z$ not free in $\alpha$, $\Gamma,\ \rho$

$(\exists\text{-L})$ $\dfrac{(\underline{f},\ h)\ (\underline{g})\ (\text{E}\underline{x})\ (\text{E}\underline{y})\colon A(h\underline{x}),\ \Gamma\ /\ \rho\ \Rightarrow\ \delta}{(\underline{f})\ (\underline{g})\ (\text{E}\underline{x})\ (\text{E}\underline{y})\colon \exists z A(z),\ \Gamma\ /\ \rho\ \Rightarrow\ \delta}$ provided $h$ not in $\Gamma, \rho, \delta$, and $\underline{y}$ not free in $A(h\underline{x})$, $\Gamma,\ \rho$

$(\forall\text{-R})$ $\dfrac{(\underline{f})\ (\underline{g},\ h)\ (\text{E}\underline{x})\ (\text{E}\underline{y})\colon \alpha\ \Rightarrow\ A(h\underline{y})\setminus\Gamma\ ;\ \rho}{(\underline{f})\ (\underline{g})\ (\text{E}\underline{x})\ (\text{E}\underline{y})\colon \alpha\ \Rightarrow\ \forall z A(z)\setminus\Gamma\ ;\ \rho}$ provided $h$ not in $\alpha, \Gamma, \rho$, and $\underline{x}$ not free in $A(h\underline{y})$, $\Gamma,\ \rho$

<u>N-rules:</u> NC and EM as before in CC.

<u>S-rules:</u> CrCr-L and CrCr-R as before, and:

$(\text{EVE})$ $\dfrac{(\underline{f})\ (\underline{g})\ (\text{E}\underline{x})\ (\text{E}\underline{y})\colon \alpha\ \Rightarrow\ \delta}{(\underline{f},\ \underline{h})\ (\underline{g})\ (\text{E}\underline{x})\ (\text{E}\underline{y})\colon \alpha\ \Rightarrow\ \delta}$ provided $\underline{h}$ not in $\alpha \Rightarrow \delta$

$(\text{IVU})$ $\dfrac{(\underline{f})\ (\underline{g})\ (\text{E}\underline{x})\ (\text{E}\underline{y})\colon \alpha\ \Rightarrow\ \delta}{(\underline{f})\ (\underline{g},\ \underline{h})\ (\text{E}\underline{x})\ (\text{E}\underline{y})\colon \alpha\ \Rightarrow\ \delta}$ provided $\underline{h}$ not in $\alpha \Rightarrow \delta$

$(\text{GenInst})$ $\dfrac{(\underline{f})\ (\underline{g})\ (\text{E}\underline{x})\ (\text{E}\underline{y})\colon \alpha_0\ \Rightarrow\ \delta}{(\underline{f})\ (\underline{g})\ (\text{E}\underline{x},\ \underline{z})\ (\text{E}\underline{y})\colon \alpha\ \Rightarrow\ \delta}$ provided $\underline{z}$ not free in $\alpha_0 \Rightarrow \delta$

$(\text{GenGen})$ $\dfrac{(\underline{f})\ (\underline{g})\ (\text{E}\underline{x})\ (\text{E}\underline{y})\colon \alpha\ \Rightarrow\ \delta_0}{(\underline{f})\ (\underline{g})\ (\text{E}\underline{x})\ (\text{E}\underline{y},\ \underline{z})\colon \alpha\ \Rightarrow\ \delta}$ provided $\underline{z}$ not free in $\alpha \Rightarrow \delta_0$

Note that all these rules are equivalential except for GenInst and GenGen. The equivalence in the cases of the rules $\exists$-L and $\forall$-R is perhaps not immediately evident, but can easily be checked using the 'rules of passage' for 'quantifiers' and Herbrandization. – It is funny that through two 'inversions' what we still call EVE effects now the *introduction* of a vacuous *universal* variable!

### 3. Examples of deductions

(All examples here are of QDC deductions unless otherwise stated. Note that we use a, b, … as 0-ary functional variables, and abbreviate a(⟨ ⟩), b(⟨ ⟩), … to just a, b, …)

∀xPx ⊢$^N$ ∃xPx

| ∀xPx | ↓ ∀ |
| (x): Px | ↓ Inst |
| Px | ↑ Gen |
| (Ea): Pa | ↑ ∃ |
| ∃xPx. | |

∀x(Px ∧ Qx) ⊢$^N$ ∃x(Qx ∨ Rx)

| ∀x(Px ∧ Qx) | ↓ ∀ |
| (x): Px ∧ Qx | ↓ ∧ |
| (x): Px, Qx | ↓ Del |
| (x): Qx | ↓ Inst |
| Qx | ↑ Exp |
| Qx / Rx | ↑ ∨ |
| Qx ∨ Rx | ↑ Gen |
| (Ea): Qa ∨ Ra | ↑ ∃ |
| ∃x(Qx ∨ Rx). | |

⊢$^N$ ∀x(Px ∨ ¬Px)

| {{ }} (= ⊤) | ↑ IVU |
| (x): {{ }} | ↑ EM |

(x): Px / ¬Px         ↑ ∨

(x): Px ∨ ¬Px      ↑ ∀

∀x(Px ∨ ¬Px).

∃x(Px ∧ ¬Px) ⊢$^N$ ⊥

∃x(Px ∧ ¬Px)      ↓ ∃

(Ea): Pa ∧ ¬Pa      ↓ ∧

(Ea): Pa, ¬Pa      ↓ NC

(Ea): { }           ↓ EVE

{ }    (= ⊥).

⊢$^N$ ¬∃x(Px ∧ ¬Px)

{{ }}           ↑ IVU

(x): {{ }}        ↑ EM

(x): Px / ¬Px        ↑ ¬¬

(x): ¬¬Px / ¬Px     ↑ ¬∧

(x): ¬(Px ∧ ¬Px)    ↑ ¬∃

¬∃x(Px ∧ ¬Px).

⊢$^N$ p ↔ ∀x p

{{ }}  (= ⊤)             ↑ EM

p / ¬p  (= p, p / ¬p, ¬p)    ↑ IVU

(x): p, p / ¬p, ¬p       ↑ ∀

p, ∀x p / ¬p, ¬p        ↑ Gen

(Ea): p, ∀x p / ¬p, ¬p    ↑ ¬∀

p, ∀x p / ¬p, ¬∀x p            ↑∧ (×2), ↑∨

p ↔ ∀x p.

‖⊢ ∀y(∀xPx → Py)‖

{{ }}                   ↑EM

¬∀xPx / ∀yPy          ↓∀

(y): ¬∀xPx / Py         ↑∨

(y): ∀xPx → Py         ↑∀

∀y(∀xPx → Py).

‖⊢ ∃y(∃xPx → Py)‖

{{ }}                   ↑EM

¬∃xPx / ∃yPy          ↓∃

(Ea): ¬∃xPx / Pa        ↑∨

(Ea): ∃xPx → Pa        ↑∃

∃y(∃xPx → Py).

‖⊢ ∃y(Py → ∀xPx)‖

{{ }}                   ↑EM

¬∀yPy / ∀xPx          ↓¬∀

(Ea): ¬Pa / ∀xPx        ↑∨

(Ea): Pa → ∀xPx        ↑∃

∃y(Py → ∀xPx).

‖¬∃xPx ⊢ ∀y(∃xPx → Py)‖

¬∃xPx                 ↑Exp

¬∃xPx / ∀yPy           ↓ ∀

(y): ¬∃xPx / Py        ↑ ∨

(y): ∃xPx → Py       ↑ ∀

∀y(∃xPx → Py).

∀xPx ∨ ∀xQx ⊢ᴺ ∀x(Px ∨ Qx)

∀xPx ∨ ∀xQx          ↓ ∨

∀xPx / ∀xQx          ↓ ∀ (×2)

(x, y): Px / Qy        ↓ Inst

(x): Px / Qx          ↑ ∨

(x): Px ∨ Qx         ↑ ∀

∀x(Px ∨ Qx).

∃x(Px ∧ Qx) ⊢ᴺ ∃xPx ∧ ∃xQx

∃x(Px ∧ Qx)         ↓ ∃, ↓ ∧

(Ea): Pa, Qa         ↑ Gen

(Ea, b): Pa, Qb      ↑ ∃ (×2)

∃xPx, ∃xQx         ↑ ∧

∃xPx ∧ ∃xQx.

∀x(Px ∧ Qx) ⊢ᴺ ∀x∀y(Px ∧ Qy)

∀x(Px ∧ Qx), ∀y(Py ∧ Qy)      ↓ ∀ (×2), ↓ ∧ (×2)

(x, y): Px, Qx, Py, Qy         ↓ Del

(x, y): Px, Qy            ↑ ∧, ↑ ∀ (×2)

∀x∀y(Px ∧ Qy).

$\forall xPx \land \forall xQx \vdash^N \forall x(Px \land Qx)$

| | |
|---|---|
| $\forall xPx \land \forall xQx$ | $\downarrow \land$ |
| $\forall xPx, \forall xQx$ | $\downarrow \forall$ |
| (x): Px, $\forall xQx$ | $\downarrow \forall$ |
| (x, y): Px, Qy | $\downarrow$ Inst |
| (x): Px, Qx | $\uparrow \land$ |
| (x): Px $\land$ Qx | $\uparrow \forall$ |
| $\forall x(Px \land Qx)$. | |

$\forall x(Px \land Qx) \vdash^N \forall xPx \land \forall xQx$

| | |
|---|---|
| $\forall x(Px \land Qx), \forall y(Py \land Qy)$ | $\downarrow \forall (\times 2)$ |
| (x, y): Px $\land$ Qx, Py $\land$ Qy | $\downarrow \land (\times 2)$ |
| (x, y): Px, Qx, Py, Qy | $\downarrow$ Del |
| (x, y): Px, Qy | $\uparrow \forall (\times 2), \uparrow \land$ |
| $\forall xPx \land \forall yQy$. | |

$\exists xPx \lor \exists xQx \vdash^N \exists x(Px \lor Qx)$

| | |
|---|---|
| $\exists xPx \lor \exists xQx$ | $\downarrow \lor$ |
| $\exists xPx \ / \ \exists xQx$ | $\downarrow \exists$ |
| (Ea): Pa $/ \exists xQx$ | $\downarrow \exists$ |
| (Ea, b): Pa $/$ Qb | $\uparrow$ Exp |
| (Ea, b): Pa $/$ Qa $/$ Pb $/$ Qb | $\uparrow \lor (\times 2)$ |
| (Ea, b): Pa $\lor$ Qa $/$ Pb $\lor$ Qb | $\uparrow \exists (\times 2)$ |
| $\exists x(Px \lor Qx) \ / \ \exists x(Px \lor Qx)$ | |

[= ∃x(Px ∨ Qx)].

∃x(Px ∨ Qx)  ⊢ᴺ  ∃xPx ∨ ∃xQx

∃x(Px ∨ Qx)                          ↓ ∃, ↓ ∨

(Ea): Pa / Qa                        ↑ Gen

(Ea, b): Pa / Qb                     ↑ ∃ (×2)

∃xPx / ∃xQx                          ↑ ∨

∃xPx ∨ ∃xQx.

∀x(Px → Qx), ∀x(Qx → Rx)  ⊢ᴺ  ∀x(Px → Rx)

∀x(Px → Qx), ∀x(Qx → Rx)                              ↓ ∀

(x): ¬Px ∨ Qx, ∀x(¬Qx ∨ Rx)                           ↓ ∀

(x, y): ¬Px ∨ Qx, ¬Qy ∨ Ry                            ↓ ∨

(x, y): ¬Px, ¬Qy ∨ Ry / Qx, ¬Qy ∨ Ry                  ↓ ∨

(x, y): ¬Px, ¬Qy / ¬Px, Ry / Qx, ¬Qy / Qx, Ry         ↓ Inst

(x): ¬Px, ¬Qx / ¬Px, Rx / Qx, ¬Qx / Qx, Rx            ↓ NC

(x): ¬Px, ¬Qx / ¬Px, Rx / Qx, Rx                      ↓ Del

(x): ¬Px / Rx                                         ↑ ∨

(x): ¬Px ∨ Rx                                         ↑ ∀

∀x(¬Px ∨ Rx).

∀x(Px → Qx), ∃xPx  ⊢ᴺ  ∃xQx

∀x(Px → Qx), ∃xPx                    ↓ ∃

(Ea): ∀x(Px → Qx), Pa                ↓ ∀

(Ea) (x): ¬Px ∨ Qx, Pa               ↓ ∨

(Ea) (x): Pa, ¬Px / Pa, Qx          ↓ Inst

(Ea): Pa, ¬Pa / Pa, Qa              ↓ NC, ↓ Del

(Ea): Qa                            ↑ ∃

∃xQx.

∀x Rxx  ⊢ᴺ  ∀x∃y Rxy

∀x Rxx                ↓ ∀

(x): Rxx              ↑ Gen

(Ef) (x): Rxf(x)      ↑ ∃

(x): ∃y Rxy           ↑ ∀

∀x∃y Rxy.

∃y∀x Rxy  ⊢  ∀x∃y Rxy

∃y∀x Rxy                ↓ ∃, ↓ ∀

(Ea) (x): Rxa           ↑ Gen

(Ea, f) (x): Rxf(x)     ↓ EVE

(Ef) (x): Rxf(x)        ↑ ∃

(x): ∃y Rxy            ↑ ∀

∀x∃y Rxy.

– And in QHC (now normal deduction):

∃y∀x Rxy                ↓ ∃, ↓ ∀

(Ea) (x): Rxa           CrCr [for: Rba is a common instance]

(b) (Ey): Rby           ↑ ∃, ↑ ∀

∀x∃y Rxy.

$\vdash^N \exists y \forall x\ Rxy \rightarrow \forall x \exists y\ Rxy$

| | |
|---|---|
| ⊤ | ↑ IVU, ↑ EM |
| (x, y): ¬Rxy / Rxy | ↑ Gen |
| (Ef) (x, y): ¬Rxy / Rxf(x, y) | ↑ Gen |
| (Ef, g) (x, y): ¬Rg(x, y)y / Rxf(x, y) | ↑ ∃ |
| (Eg) (x, y): ¬Rg(x, y)y / ∃z Rxz | ↑ ¬∀ |
| (x, y): ¬∀u Ruy / ∃z Rxz | ↑ ∀ |
| (y): ¬∀u Ruy / ∀x∃z Rxz | ↑ ¬∃ |
| ¬∃y∀u Ruy / ∀x∃z Rxz | ↑ ∨ |

$\exists y \forall u\ Ruy \rightarrow \forall x \exists z\ Rxz.$

$\vdash_{QCC} \exists y \forall x\ Rxy \Rightarrow \forall x \exists y\ Rxy$

| | |
|---|---|
| { } ⇒ { } | IVU, EVE |
| (a) (b): { } ⇒ { } | CrCr-L |
| (a) (b): Rba ⇒ { } | CrCr-R |
| (a) (b): Rba ⇒ Rba | GenInst |
| (a) (b) (Ex): Rxa ⇒ Rba | GenGen |
| (a) (b) (Ex) (Ey): Rxa ⇒ Rby | ∃-R |
| (a) (b) (Ex): Rxa ⇒ ∃y(Rby) | ∀-L |
| (a) (b): ∀x(Rxa) ⇒ ∃y(Rby) | ∀-R |
| (a): ∀x(Rxa) ⇒ ∀x∃y Rxy | ∃-L |

$\exists y \forall x\ Rxy \Rightarrow \forall x \exists y\ Rxy.$

$\forall x \forall y(Rxy \rightarrow Ryx),\ \forall x \forall y \forall z(Rxy \land Ryz \rightarrow Rxz)\ \vdash^N\ \forall x \forall y(Rxy \rightarrow Rxx)$

∀x∀y(Rxy → Ryx), ∀x∀y∀z(Rxy ∧ Ryz → Rxz)                    ↓ ∀ (×5)

(x, y, u, v, z): Rxy → Ryx, Ruv ∧ Rvz → Ruz                    ↓ Inst

(x, y): Rxy → Ryx, Rxy ∧ Ryx → Rxx          [then here a ↓↑ deduction following the pattern of the ↓↑ deduction which we saw in propositional DC from p → q, (p ∧ q) → r to p → r]

(x, y): Rxy → Rxx                                              ↑ ∀ (×2)

∀x∀y(Rxy → Rxx).

<u>∃x∃y(Rxy ∧ Ryx), ∀x∀y∀z(Rxy ∧ Ryz → Rxz)  ⊢<sup>N</sup> ∃x Rxx</u>

∃x∃y(Rxy ∧ Ryx), ∀x∀y∀z(Rxy ∧ Ryz → Rxz)                    ↓ ∃, ↓ ∀, ↓ ∧

(Ea, b) (x, y, z): Rab, Rba, Rxy ∧ Ryz → Rxz                  ↓ ∨

(Ea, b) (x, y, z): Rab, Rba, ¬(Rxy ∧ Ryz) / Rab, Rba, Rxz    ↓ ¬∧

(Ea, b) (x, y, z): Rab, Rba, ¬Rxy / Rab, Rba, ¬Ryz / Rab, Rba, Rxz    ↓ Inst

(Ea, b): Rab, Rba, ¬Rab / Rab, Rba, ¬Rba / Rab, Rba, Raa      ↓ Del

(Ea, b): Rab, Rba, ¬Rab / Rab, Rba, ¬Rba / Raa              ↓ NC, ↓ EVE

(Ea): Raa                                                      ↑ ∃

∃x Rxx.

<u>∀x∃y Rxy, ∀x∀y(Rxy → Ryx)  ⊢<sup>N</sup>  ∀x∃y(Rxy ∧ Ryx)</u>

∀x∃y Rxy, ∀x∀y(Rxy → Ryx)                    ↓ ∀

(x): ∃y Rxy, ∀x∀y(Rxy → Ryx)                 ↓ ∃

(Ef) (x): Rxf(x), ∀x∀y(Rxy → Ryx)            ↓ ∀ (×2)

(Ef) (x, y, z): Rxf(x), Ryz → Rzy            ↓ ∨

(Ef) (x, y, z): Rxf(x), ¬Ryz / Rxf(x), Rzy   ↓ Inst (×2)

(Ef) (x): Rxf(x), ¬Rxf(x) / Rxf(x), Rf(x)x   ↓ NC

(Ef) (x): Rxf(x), Rf(x)x                    ↑ ∧

(Ef) (x): Rxf(x) ∧ Rf(x)x                   ↑ ∃

(x): ∃y(Rxy ∧ Ryx)                          ↑ ∀

∀x∃y(Rxy ∧ Ryx).

### 4.  Some derived rules

The formula say ∀xPx implies not only an instance, like Pu, but also any conjunction of instances, like Pu ∧ Pv. This (Pu ∧ Pv) is not literally an instance of ∀xPx; it is however literally an instance of the strange formula ∀x∀y(Px ∧ Py), which as is easily seen is in fact equivalent to the simple ∀xPx. We may think of ∀x∀y(Px ∧ Py) as a kind of 'unfolding' of ∀xPx, and of ∀xPx as a kind of 'folding' of ∀x∀y(Px ∧ Py). The 'import' is the same, but in the one case it is presented 'compactly' and in the other in 'spread out' form. – One way of arguing from ∀xPx to Pu ∧ Pv is then first 'unfolding' and then applying ordinary instantiation.

Dually, ∃xPx is implied not just by an instance like Pu but also by any disjunction of instances like Pu ∨ Pv. This (Pu ∨ Pv) is literally an instance of the strange formula ∃x∃y(Px ∨ Py), equivalent to ∃xPx. Again we may think of ∃x∃y(Px ∨ Py) as 'unfolding' ∃xPx, and ∃xPx as 'folding' ∃x∃y(Px ∨ Py). We may argue from Pu ∨ Pv to ∃xPx by first generalizing and then folding.

– We give in the present section some derived rules of our calculi of directional deduction related to these curious phenomena. (In addition to being interesting in its own right, this material will also be useful later in other developments.)

– Let us consider first the following QDC deduction of the datum Pu, Pv from (x): Px, including as proper part a deduction of (x, y): Px, Py from (x): Px.

(x): Px                    ↑ ∀

∀xPx, ∀yPy            ↓ ∀ (×2)

(x, y): Px, Py          ↓ Inst (×2)

Pu, Pv.

– A few things are worth noting about this deduction:

(1) Like with ordinary Instantiation, one might expect this 'Multiple Instantiation' to be *purely analytic.* But the above deduction is ↑↓, and nor *is* there a ↓-deduction for these data. For (x): Px is already an elementary datum: so no ↓ G-rule can be applied to it, nor of course ↓NC (which is the only ↓ N-rule), nor ↓EVE, nor ↓Del if we hope to deduce Pu, Pv; so there remains ↓Inst, which however once applied will take us at best to Pu or to Pv singly, from which the other does not even *follow* (semantically) and so certainly cannot be deduced.

(2) As for the 'Unfolding' deduction represented by the first three lines, one might have expected a *purely synthetic* deduction. But the given deduction is ↑↓, and nor is there a ↑-deduction. (For (x, y): Px, Py is elementary, etc. etc.)

(3) *This* deduction (first three lines) uses only G-rules and so is invertible – giving us a ↑↓-deduction for the corresponding 'Folding'. Again, one would expect for Folding a *purely analytic* deduction. And now actually there *is* one: just instantiate y as x.

– All this applies not just in the above simple example but also more generally. Let us give just one more example before stating the fully general rules.

| | |
|---|---|
| (x): Px / Qx | ↑ ∨, ↑ ∀ |
| ∀x(Px ∨ Qx), ∀y(Py ∨ Qy) | ↓ ∀ (×2) |
| (x, y): Px ∨ Qx, Py ∨ Qy | ↓ ∨ |
| (x, y): Px, Py ∨ Qy / Qx, Py ∨ Qy | ↓ ∨ (×2) |
| (x, y): Px, Py / Px, Qy / Qx, Py / Qx, Qy | ↓ Inst (×2) |
| Pu, Pv / Pu, Qv / Qu, Pv / Qu, Qv. | |

The points (1)–(3) above apply equally here; only note that in the present case (which is closer to 'the general case') we use also ↓Del, after instantiating y as x, in the purely analytic deduction for 'Folding'. Note that the last line is in effect a kind of 'conjunction' (through 'distribution') of the instances Pu / Qu and Pv / Qv of the initial datum (x): Px / Qx.

– The dual case of Existential Folding and Unfolding and 'Generalized Generalization' in QDDC is exactly analogous. (There are of course correlates of these existential principles in QDC itself, as well as correlates of the universal principles in QDDC. Here we are considering only the 'fullest' dualization.) Thus e.g. corresponding to the first example above we have:

Pu \ Pv                    ↑ Gen (×2)

(Ex, y): Px \ Py           ↑ ∃ (×2)

∃xPx \ ∃yPy   (= ∃xPx)     ↓ ∃

(Ex): Px.

– Similar remarks as before apply here too.

– We can now state the fully general derived rules. Their derivation should be obvious from the preceding examples. – We write the variables w̲ separately to emphasize that (like the f̲) they are 'inactive' in the rule.

*Universal Unfolding Rule* (derived ↑↓ QDC-rule):

$$\text{(E}\underline{f}\text{) }(\underline{w})\,(\underline{x})\text{: }\Gamma(\underline{x})\,/\,\Delta(\underline{x})\,/\,\dots$$

(↑↓UnivUnfold)   $\overline{\text{(E}\underline{f}\text{) }(\underline{w})\,(\underline{x},\,\underline{y},\,\underline{z},\,\dots)\text{:}\,\big/\{(\text{C}(\underline{x}),\,\text{D}(\underline{y}),\,\text{E}(\underline{z}),\,\dots)\text{:}}$

$$\{\text{C, D, E, }\dots\}\subseteq\{\Gamma,\,\Delta,\,\dots\}\}$$

*Universal Folding Rule* (↓UnivFold) (derived ↓ QDC-rule): this is the inverse of the preceding rule. (In the purely analytic derivation we instantiate all of y̲, z̲, … as x̲, then perform appropriate deletions.)

*Existential Unfolding Rule* (derived ↑ QDDC-rule):

$$\text{(}\underline{f}\text{) }(\text{E}\underline{w})\,(\text{E}\underline{x})\text{: }\Gamma(\underline{x})\,;\,\Delta(\underline{x})\,;\,\dots$$

(↑ExistUnfold)   $\overline{\text{(}\underline{f}\text{) }(\text{E}\underline{w})\,(\text{E}\underline{x},\,\underline{y},\,\underline{z},\,\dots)\text{:}\,⁏\{(\text{C}(\underline{x})\setminus\text{D}(\underline{y})\setminus\text{E}(\underline{z})\setminus\dots)\text{:}}$

$$\{\text{C, D, E, }\dots\}\subseteq\{\Gamma,\,\Delta,\,\dots\}\}$$

(In the purely synthetic derivation we make appropriate expansions of the components, then generalize the appropriate variables.)

*Existential Folding Rule* (↑↓ExistFold) (derived ↑↓ QDDC-rule): the inverse of the preceding rule.

*Generalized Instantiation Rule* (derived ↑↓ QDC-rule): (where t̲, u̲, v̲, … are appropriately-sized sequences of arbitrary terms):

$(\uparrow\downarrow\text{GenInst})$ 
$$\frac{(\text{E}\underline{f})\ (\underline{w})\ (\underline{x})\colon \Gamma(\underline{x})\ /\ \Delta(\underline{x})\ /\ \ldots}{(\text{E}\underline{f})\ (\underline{w})\colon\big/\{(\text{C}(\underline{t}),\ \text{D}(\underline{u}),\ \text{E}(\underline{v}),\ \ldots)\ :\ \{\text{C},\ \text{D},\ \text{E},\ \ldots\}\subseteq\{\Gamma,\ \Delta,\ \ldots\}\}}$$

*Generalized Generalization Rule* (derived $\uparrow\downarrow$ QDDC-rule):

$(\uparrow\downarrow\text{GenGen})$ 
$$\frac{(\underline{f})\ (\text{E}\underline{w})\colon\,\raisebox{0.3ex}{\bf ;}\{(\text{C}(\underline{t})\setminus\text{D}(\underline{u})\setminus\text{E}(\underline{v})\setminus\ldots)\ :\ \{\text{C},\ \text{D},\ \text{E},\ \ldots\}\subseteq\{\Gamma,\ \Delta,\ \ldots\}\}}{(\underline{f})\ (\text{E}\underline{w})\ (\text{E}\underline{x})\colon \Gamma(\underline{x});\ \Delta(\underline{x});\ \ldots}$$

# 5. Herbrand theorems

Given a formula $\forall \underline{x} A \underline{x}$ or $\exists \underline{x} A \underline{x}$ (typically though not necessarily with $A\underline{x}$ quantifier-free), by a *Herbrand instance* of such a formula we mean a formula obtainable from $A\underline{x}$ by uniform substitution, of terms built with arbitrary individual variables and the functional variables (if any) in $A\underline{x}$, for the variables $\underline{x}$ (in their free occurrences in $A\underline{x}$). More generally, we may also speak in the similar sense of a Herbrand instance *w.r.t. a set of functional variables* (the other notion being the special case where the set consists of the functional variables in $A\underline{x}$).

By the *Herbrand universe* for set S of functional variables we mean the set of all terms built from arbitrary individual variables and the functional variables in S.

*Proposition* (Herbrand's theorem, satisfiability form).  For $A\underline{x}$ quantifier-free: $\forall \underline{x} A \underline{x}$ is satisfiable iff the set $\Gamma$ of its Herbrand instances is *propositionally* satisfiable.

*Proof.*  ($\Rightarrow$) is immediate, since $\forall \underline{x} A \underline{x} \vDash \Gamma$ and so of course if $\forall \underline{x} A \underline{x}$ is satisfiable then $\Gamma$ is satisfiable and hence a fortiori propositionally satisfiable.

($\Leftarrow$). Suppose $\Gamma$ is propositionally satisfiable, and let $\sigma$ be a propositional interpretation (i.e. assignment of truth-values to the atomic formulas) verifying $\Gamma$. We define an interpretation (now in the sense of predicate logic) $\sigma'$ as follows. Let S be the set of predicate variables and functional variables in A. The domain of the interpretation is the Herbrand universe H for the set of functional variables in A. For individual variable v, we put: $\sigma'(v) = v$. For functional variable f $\in$ S, we take $\sigma'(f)$ as the function giving the term $f(t_1 \ldots t_n)$ to each sequence of terms $t_1 \ldots t_n \in H$. For g $\notin$ S we may take $\sigma'(g)$ as e.g. the constant function with identical value x (the alphabetically first variable). And for predicate variable P, we put $\sigma'(P) = \{\langle t_1 \ldots t_n \rangle \in H^n : \sigma(Pt_1 \ldots t_n) = T\}$. – Thus clearly $\sigma'$ will agree with $\sigma$ on the truth-values of atomic formulas occurring in formulas of $\Gamma$, whence also on the truth-values of truth-functional compounds therefrom, and thus on all formulas in $\Gamma$. So $\sigma'$ verifies $\Gamma$. – Now note that, w.r.t. the interpretation $\sigma'$, the set $\Gamma$ is not an 'ordinary' set of instances of the universal quantification $\forall \underline{x} A \underline{x}$ but in fact includes one instance for each sequence of individuals (possible values for the variables $\underline{x}$) in the domain. So w.r.t. $\sigma'$ the set $\Gamma$ has in fact the exact same import as the formula $\forall \underline{x} A \underline{x}$, i.e.

σ′ verifies Γ iff σ′ verifies $\forall \underline{x} A \underline{x}$. So we may conclude that σ′ verifies $\forall \underline{x} A \underline{x}$ and so $\forall \underline{x} A \underline{x}$ is satisfiable. □

*Remark*. By the Compactness Theorem for Propositional Logic, the r.h.s. of the Proposition above is equivalent to: All finite subsets of Γ are propositionally satisfiable. – The similar point applies also below to the two other forms of the Herbrand theorem.

*Proposition* (Herbrand's theorem, validity form). For $A\underline{x}$ quantifier-free: $\exists \underline{x} A \underline{x}$ is valid iff the set Γ of its Herbrand instances is *propositionally* disjunctively valid (i.e. every propositional interpretation verifies at least one formula in Γ).

*Proof*. $\exists \underline{x} A \underline{x}$ is valid ⇔ $\forall \underline{x} \neg A \underline{x}$ is unsatisfiable ⇔ (by the satisfiability form of Herbrand's theorem) the corresponding set of Herbrand instances $\{\neg A\underline{t}, \neg A\underline{u}, \ldots\}$ is propositionally unsatisfiable ⇔ Γ = $\{A\underline{t}, A\underline{u}, \ldots\}$ is propositionally disjunctively valid. □

*Proposition* (Herbrand's theorem, implicational form). For $A\underline{x}$ and $B\underline{y}$ quantifier-free and S = the set of functional variables occurring in either $A\underline{x}$ or $B\underline{x}$: $\forall \underline{x} A \underline{x}$ implies $\exists \underline{y} B \underline{y}$ iff the sets Γ, Δ of their respective Herbrand instances w.r.t. S are s.t. every propositional interpretation which verifies all the formulas in Γ verifies at least one formula in Δ.

*Proof*. $\forall \underline{x} A \underline{x}$ implying $\exists \underline{y} B \underline{y}$ is of course equivalent to the unsatisfiability of $\forall \underline{x} A \underline{x} \land \neg \exists \underline{y} B \underline{y}$, or equivalently of (renaming variables if needed) $\forall \underline{x} \forall \underline{y} (A\underline{x} \land \neg B\underline{y})$. This in turn, by the satisfiability form of Herbrand's theorem, is equivalent to the propositional unsatisfiability of the set of instances $\{A\underline{t} \land \neg B\underline{u}, A\underline{t}' \land \neg B\underline{u}', \ldots\}$, i.e. the propositional unsatisfiability of $\{A\underline{t}, A\underline{t}', \ldots, \neg B\underline{u}, \neg B\underline{u}', \ldots\}$; i.e. the statement: every propositional interpretation which verifies all the formulas in $\{A\underline{t}, A\underline{t}', \ldots\}$ = Γ verifies at least one formula in $\{B\underline{u}, B\underline{u}', \ldots\}$ = Δ. □

# 6. Completeness etc.

*Proposition* (Refutation-completeness of QDC).  Every unsatisfiable standard datum is QDC-refutable.

*Proof.*  We begin the refutation by applying the $\downarrow$ G-rules until we reach an elementary datum (E$\underline{f}$) ($\underline{x}$): $\alpha$. Since this is equivalent to the original datum, it is also unsatisfiable, whence so is ($\underline{x}$): $\alpha$. Then by Herbrand's theorem (satisfiability form) and the Compactness Theorem for propositional logic, there exist Herbrand instances $\alpha_1$ … $\alpha_n$ of ($\underline{x}$): $\alpha$ s.t. $\{\alpha_1 \ldots \alpha_n\}$ is propositionally unsatisfiable. Now letting $\beta$ be the datum which 'conjoins' (by 'distribution') $\alpha_1 \ldots \alpha_n$, we can then deduce (E$\underline{f}$): $\beta$ from (E$\underline{f}$) ($\underline{x}$): $\alpha$ by $\uparrow\downarrow$GenInst. Since $\beta$ is propositionally unsatisfiable, so is every component of $\beta$ – which since $\beta$ is elementary means that every such component contains some explicit contradiction A, ¬A. So from (E$\underline{f}$): $\beta$ successive applications of $\downarrow$NC and $\downarrow$EVE will eventually yield $\perp$.  □

*Remarks.*  (1) Because of the use of $\uparrow\downarrow$GenInst this refutation will not in general be purely analytic. It would be good if this could somehow be 'rectified'.

(2) Instead of appealing to the previously proved Herbrand theorem, we could reformulate the above proof so as to include the main ideas in the proof of the Herbrand theorem. (And the same holds for the other completeness proofs below.) One advantage of this reformulation is that we can then strengthen the Proposition to: If standard datum $\gamma$ is irrefutable (in QDC), then $\gamma$ is satisfiable *in denumerable domain*. (This, together with the [obvious] soundness of QDC, immediately yields Löwenheim's theorem, i.e. the fact that if a formula [or standard datum] is satisfiable then it is satisfiable in denumerable domain.) The reformulated proof is then in sketch as follows. Suppose $\gamma$ is irrefutable. Then in particular it cannot be refuted by reduction via $\downarrow$ G-rules to elementary datum (E$\underline{f}$) ($\underline{x}$): $\alpha$ followed by application of Generalized Instantiation to pass to a 'conjunction' of Herbrand instances. So all such 'conjunctions' are propositionally satisfiable; and so by propositional Compactness the whole set of all Herbrand instances is propositionally satisfiable. But then ($\underline{x}$): $\alpha$ is verified by interpretation (the Herbrand interpretation) with denumerable domain, whence so is (E$\underline{f}$) ($\underline{x}$): $\alpha$, whence so is $\gamma$.

*Questions.* (1) To what extent do 'Modularity Lemmas' hold for our Directional Deduction calculi for predicate logic?

(2) Presumably QDC is *deductively* complete (i.e. for any standard data $\alpha$, $\beta$: if $\alpha \vDash \beta$ then $\alpha \vdash_{QDC} \beta$)? If at least a mild form of Modularity holds then this may be derived from the above refutation-completeness theorem as follows. Suppose $\alpha \vDash \beta$. Put A $:=$ cf($\alpha$) and B $:=$ cf($\beta$). (The degenerate cases of data without characteristic formulas are readily dealt with by a separate argument.) Obviously $\alpha \vdash$ A and B $\vdash \beta$ (by G-rules); so it is enough to show A $\vdash$ B. This we can then do (modulo Modularity) thus:

A $\qquad\qquad$ ↑ EM

A, B / A, ¬B

… $\qquad\qquad$ [refutation of A, ¬B]

A, B $\qquad$ ↓ Del

B.

(3) What about *normal* or *quasi-normal* deduction completeness?

– We proceed now to the analogous case of *dual data*.

*Proposition* (Proof-completeness of QDDC). Every valid standard dual datum is QDDC-provable.

*Proof.* We construct the proof from the bottom up starting by reducing via ↑ G-rules the final dual datum to equivalent elementary dual datum (f) (E<u>x</u>): $\delta$. Since this is valid, so is (E<u>x</u>): $\delta$. Then by Herbrand's theorem (validity form) and Compactness for propositional logic, there exist Herbrand instances $\delta_1 \ldots \delta_n$ of (E<u>x</u>): $\delta$ s.t. $\{\delta_1 \ldots \delta_n\}$ is propositionally disjunctively valid. Now letting $\varepsilon$ be the dual datum which 'disjoins' (by 'distribution') $\delta_1 \ldots \delta_n$, we can deduce (f) (E<u>x</u>): $\delta$ from (f): $\varepsilon$ by ↑↓GenGen. Since $\varepsilon$ is propositionally valid, so is every component of $\varepsilon$ – which since $\varepsilon$ is elementary means that every such component contains pair of contradictory formulas A, ¬A. So (f): $\varepsilon$ can be gotten from ⊤ by successive applications of ↑IVU and ↑EM. □

(The Remarks and Questions above apply also here to dual data mutatis mutandis.)

*Proposition* (Deductive completeness of QHC). For every standard datum α and standard dual datum δ: If α ⊨ δ then α ⊢<sub>QHC</sub> δ.

*Proposition* (Deductive completeness of QHC). For every standard datum α and standard dual datum δ: If $\alpha \vDash \delta$ then $\alpha \vdash_{QHC} \delta$.

*Proof.* We begin by reducing α via ↓ G-rules to equivalent elementary datum (E$\underline{f}$) ($\underline{x}$): β, and δ from the bottom up via ↑ G-rules to equivalent elementary dual datum ($\underline{g}$) (E$\underline{y}$): ε. Since (E$\underline{f}$) ($\underline{x}$): β implies ($\underline{g}$) (E$\underline{y}$): ε, also ($\underline{x}$): β implies (E$\underline{y}$): ε. So by Herbrand's theorem (implicational form) and propositional Compactness, there are Herbrand instances (w.r.t. the functional variables in β, ε) $\beta_1 \ldots \beta_n$ and $\varepsilon_1 \ldots \varepsilon_k$ s.t. every propositional interpretation which verifies all of $\beta_1 \ldots \beta_n$ verifies at least one of $\varepsilon_1 \ldots \varepsilon_k$. Now let β′ be the datum which 'conjoins' (by 'distribution') $\beta_1 \ldots \beta_n$, and similarly let ε′ be the dual datum which 'disjoins' $\varepsilon_1 \ldots \varepsilon_k$. Now by ↑↓UnivUnfold we pass from (E$\underline{f}$) ($\underline{x}$): β to (E$\underline{f}$) ($\underline{x}$, …): β\* with the unfolding which would be appropriate for the subsequent direct instantiation to (E$\underline{f}$): β′. And similarly we pass bottom-up from ($\underline{g}$) (E$\underline{y}$): ε via ↑↓ExistFold to ($\underline{g}$) (E$\underline{y}$, …): ε\* where this is the unfolded dual datum which would be appropriate to be obtained by direct generalization from ($\underline{g}$): ε′. – Now note that, since β′ propositionally implies ε′, and they are both elementary, β′ must then weakly criss-cross ε′. So the Criss-Crossing rule licenses the passage from (E$\underline{f}$) ($\underline{x}$, …): β\* to ($\underline{g}$) (E$\underline{y}$, …): ε\*, with β′ and ε′ being the relevant instances. ◻

*Remarks.* (1) Again, because of ↑↓UnivUnfold and ↑↓ExistFold the deduction is not quite normal and it would be good if this could be 'rectified'.

(2) The present formulation of QHC, and in particular the treatment of transition (i.e. the Criss-Crossing rule), do not seem to me to be very satisfactory; although I do not at present see how to improve this. The trouble is that the Criss-Crossing rule bundles together too many inferences which, in a more satisfactory treatment, should be separable. Essentially, from the unfolded elementary datum one would like to *instantiate*, then remove impure components by *NC*; and similarly go *to* the unfolded elementary dual datum by *generalization*, with application of *EM* before to introduce impure components; and then for the transition one would use some relatively straightforward form of *criss-crossing*. – There *is* in fact a simple way of proceeding along these lines, but at the cost of introducing further resources in the prefixes: thus writing (E$\underline{f}$) ($\underline{x}$): α and ($\underline{g}$) (E$\underline{y}$): δ for the already unfolded datum and dual datum, and $\alpha_0$, $\delta_0$ for the relevant Herbrand instances (w.r.t. $\underline{f}$, $\underline{g}$) and $\alpha_1$, $\delta_1$ for their corresponding purified forms, we would proceed thus (where CrCr now is formulated in terms of propositional criss-crossing and quantifier-interchange [so a fastidious person might still say that it is bundling together more than one principle – but of course bundling together two is much better than bundling together some half-dozen!]):

| | |
|---|---|
| (E$\underline{f}$) ($\underline{x}$): $\alpha$ | ↑ IVU |
| (E$\underline{f}$) ($\underline{g}$) ($\underline{x}$): $\alpha$ | ↓ Inst |
| (E$\underline{f}$) ($\underline{g}$): $\alpha_0$ | ↓ NC |
| (E$\underline{f}$) ($\underline{g}$): $\alpha_1$ | CrCr |
| ($\underline{g}$) (E$\underline{f}$): $\delta_1$ | ↑ EM |
| ($\underline{g}$) (E$\underline{f}$): $\delta_0$ | ↑ Gen |
| ($\underline{g}$) (E$\underline{f}$) (E$\underline{y}$): $\delta$ | ↓ EVE |
| ($\underline{g}$) (E$\underline{y}$): $\delta$. | |

– One might consider liberalizing the notions of data and dual data so as to allow for something like the above deduction. But note that this requires further changes to the calculus: for e.g. how are we going to 'analyze' existential individual quantification in data? If there is a universal *functional* prefix, the existential individual quantifier cannot 'go out' as existential quantifier over functions *of individuals*! The obvious change here is to restrict the G-rules for quantifiers to context where the second part of the prefix is purely 'individual'. This still permits canonical deductions of a suitably 'standard' dual datum from a suitably 'standard' datum.

Another possible approach may be along the lines of the formulation of propositional HC with Distr instead of CrCr. Thus the transition rule would turn a (Skolemite) elementary datum into an equivalent (Herbrandite) elementary dual datum, by (roughly speaking) turning the (initial) existential functional prefix into the (final) existential individual prefix, and similarly for the (final) universal individual prefix and the (initial) universal functional prefix, and 'inverting' the dependency structure in the (distributed) matrix. It can then be hoped that (as in propositional case), given the hypothesis of semantic implication, some natural rules may lead from this 'reformulation' of the given elementary datum to the desired elementary dual datum.

*Proposition* (Proof-completeness of QCC for regular consecutions). All valid regular consecutions are provable in QCC.

*Proof*. We construct the proof from the bottom up, beginning by reducing the given regular consecution via the G-rules to an equivalent elementary regular consecution ($\underline{f}$) ($\underline{g}$) (E$\underline{x}$) (E$\underline{y}$): $\alpha \Rightarrow \delta$. Since this is valid, we have that the datum ($\underline{x}$): $\alpha$ implies the dual datum (E$\underline{y}$): $\delta$; and so by Herbrand's theorem there are 'conjunctions' and 'disjunctions' of Herbrand instances, $\alpha_0(\underline{f}, \underline{g})$, $\delta_0(\underline{f}, \underline{g})$, with $\alpha_0(\underline{f}, \underline{g}) \vDash \delta_0(\underline{f}, \underline{g})$. So from { } $\Rightarrow$ { } by IVU and EVE we deduce ($\underline{f}$) ($\underline{g}$): { } $\Rightarrow$ { }, and then (note that $\alpha_0$ and $\delta_0$ are both elementary) by CrCr-L, CrCr-R, NC, EM we deduce ($\underline{f}$) ($\underline{g}$): $\alpha_0(\underline{f}, \underline{g}) \Rightarrow \delta_0(\underline{f}, \underline{g})$. Then finally we have:

($\underline{f}$) ($\underline{g}$): $\alpha_0(\underline{f}, \underline{g}) \Rightarrow \delta_0(\underline{f}, \underline{g})$     GenInst

($\underline{f}$) ($\underline{g}$) (E$\underline{x}$): $\alpha \Rightarrow \delta_0(\underline{f}, \underline{g})$     GenGen

($\underline{f}$) ($\underline{g}$) (E$\underline{x}$) (E$\underline{y}$): $\alpha \Rightarrow \delta$.     □

CHAPTER 3

DIRECTIONAL DEDUCTION FOR MODAL LOGIC

## 1. Directional Deduction for modal propositional logic (S5): Basic notions

We proceed now to Directional Deduction for modal propositional logic (S5). The basic idea here is to treat modal formulas $\Box A$ and $\Diamond A$ as if they were the world-quantifications $\forall w\, T^w A$ and $\exists w\, T^w A$, so that we can use essentially the same ideas as in the case of classical predicate logic, with a few adaptations and supplementations – in particular supplementations to deal with the world-relative truth-operator $T^w$.

– *Formulas* are the usual formulas of propositional S5 with primitive connectives $\neg, \wedge, \vee, \Box, \Diamond$. Now, for constructing data (and dual data) we will use also: a denumerable stock of *world-variables* w, v, u, …; and for each $n \geq 0$, a denumerable stock of *n-ary functional variables* (variables for functions from n-tuples of worlds to worlds). We usually write f, g, … for such variables, letting their arity be inferred from the context; but we use also m, n, … specifically for 0-ary functional variables. Then *world-terms* are constructed from world-variables and functional variables in the natural way. A (world-relative) *truth-operator* is an expression of the form $T^t$ where t is a world-term. A *qualified formula* is a string of zero or more truth-operators followed by a formula. (Thus a formula on its own counts also as a 'qualified formula'.) Finally a *datum* is an expression of the form

$$(E\underline{f})\ (\underline{w})\colon \Gamma\ /\ \Delta\ /\ \ldots$$

where $\underline{f}$ indicates a finite set of functional variables, $\underline{w}$ a finite set of world-variables, and $\Gamma\ /\ \Delta\ /\ \ldots$ a finite set of finite sets of qualified formulas. We say that a datum is *elementary* if all its qualified formulas are literals preceded by at most one truth-operator.

Basic semantical notions are defined in obvious way within the scheme of possible world semantics. Thus a *model* consists of a non-empty set ('worlds') with designated element, together with an assignment of worlds to the world-variables, sets of worlds to the propositional variables, and functions (on worlds) of suitable arity to the functional variables; then *verification* (of datum by model) is defined in obvious way, and so on.

As in predicate logic, we use ⊥ for the (unsatisfiable) prefix-free datum with { } as matrix, and ⊤ for the (valid) prefix-free datum with {{ }} as matrix.

We define $T^t[\alpha]$ as the datum resulting from the datum $\alpha$ by insertion of $T^t$ at the beginning of all qualified formulas in all components. This is so to speak semantically justified, since a conjunction (disjunction) being true in a world w amounts to every conjunct (some disjunct) being true in w, and similarly for constant-domain universal or existential quantifications (which indeed themselves amount to [possibly infinite] conjunctions and disjunctions). – Note here the limiting cases $T^t[\bot] = \bot$ and $T^t[\top] = \top$.

– Finally as one would expect a *dual datum* is an expression

$$(\underline{f}) \ (E\underline{w}): \Gamma; \Delta; \dots$$

and a *consecution* an expression

$$(\underline{f}) \ (\underline{g}) \ (E\underline{w}) \ (E\underline{v}): \alpha \Rightarrow \delta$$

where $\alpha$ is datum-matrix and $\delta$ dual-datum--matrix.

– Again here various definitions, conventions etc. will be taken over from the previous cases (classical propositional logic, classical predicate logic) without explicit note.

## 2. The calculi S5DC, S5DDC, S5HC, S5CC

As in earlier cases, we give here four calculi: Modal Data Calculus (S5DC), Modal Dual Data Calculus (S5DDC), Modal Hybrid Calculus (S5HC), and Modal Consecution Calculus (S5CC). We state explicitly only the basic rules of the Modal Data Calculus, and leave to the reader their straightforward dualization yielding the basic rules of the Modal Dual Data Calculus. Also, S5HC is defined from S5DC and S5DDC in the obvious way, with a Criss-Crossing rule exactly as in QHC (in the verbal expression of the rule, just change x̲, y̲ to w̲, v̲); and so also S5CC is defined in the obvious way.

The rules of S5DC themselves are very similar to the rules of QDC. The main difference is the structural ↕T rules for introduction and elimination of non-rightmost truth-operators.

– We use τ to indicate an arbitrary (possibly empty) string of truth-operators.

*Rules of S5DC*:

<u>G-rules:</u>

↕∧, ↕∨, ↕¬¬, ↕¬∧, ↕¬∨ as before: E.g. ↕∧ is now:

$$(E\underline{f}) \ (\underline{w}){:}\ \tau\ A \wedge B, \Gamma\ /\ \rho$$
$$\updownarrow$$
$$(E\underline{f}) \ (\underline{w}){:}\ \tau A, \tau B, \Gamma\ /\ \rho$$

(↕□)
$$(E\underline{f}) \ (\underline{w}){:}\ \tau\square A, \Gamma\ /\ \rho$$
$$\updownarrow$$
$$(E\underline{f}) \ (\underline{w}, v){:}\ \tau T^v A, \Gamma\ /\ \rho$$
provided v does not occur in τ, Γ, ρ.

(↕¬◇)
$$(E\underline{f}) \ (\underline{w}){:}\ \tau\neg\Diamond A, \Gamma\ /\ \rho$$
$$\updownarrow$$
$$(E\underline{f}) \ (\underline{w}, v){:}\ \tau T^v \neg A, \Gamma\ /\ \rho$$
provided v not in τ, Γ, ρ.

(↕◇)
$$(E\underline{f}) \ (\underline{w}){:}\ \tau\Diamond A, \Gamma\ /\ \rho$$
$$\updownarrow$$
$$(E\underline{f}, g) \ (\underline{w}){:}\ \tau T^{g(\underline{w})} A, \Gamma\ /\ \rho$$
provided g not in τ, Γ, ρ.

$$(\updownarrow\neg\square) \quad \begin{array}{c} (\text{E}\underline{f})\ (\underline{w})\text{: } \tau\neg\square A,\ \Gamma\ /\ \rho \\ \updownarrow \\ (\text{E}\underline{f},\ g)\ (\underline{w})\text{: } \tau T^{g(w)}\neg A,\ \Gamma\ /\ \rho \end{array} \quad \text{provided g not in } \tau,\ \Gamma,\ \rho$$

<u>N-rules:</u>

$$(\text{NC}) \quad \begin{array}{c} (\text{E}\underline{f})\ (\underline{w})\text{: } \Gamma,\ \tau A,\ \tau\neg A\ /\ \rho \\ \downarrow \\ (\text{E}\underline{f})\ (\underline{w})\text{: } \rho \end{array} \qquad (\text{EM}) \quad \begin{array}{c} (\text{E}\underline{f})\ (\underline{w})\text{: } \Gamma,\ \tau A\ /\ \Gamma,\ \tau\neg A\ /\ \rho \\ \uparrow \\ (\text{E}\underline{f})\ (\underline{w})\text{: } \Gamma\ /\ \rho \end{array}$$

<u>S-rules:</u>

↓Del, ↑Exp, ↓EVE, ↑IVU as before, and:

$$(\text{Inst}) \quad \begin{array}{c} (\text{E}\underline{f})\ (\underline{w},\ v)\text{: } \alpha(v) \\ \downarrow \\ (\text{E}\underline{f})\ (\underline{w})\text{: } \alpha(t) \end{array} \qquad (\text{Gen}) \quad \begin{array}{c} (\text{E}\underline{f},\ g)\ (\underline{w})\text{: } \alpha(g\underline{w}) \\ \uparrow \\ (\text{E}\underline{f})\ (\underline{w})\text{: } \alpha(t) \end{array}$$

$$(\updownarrow\text{T}) \quad \begin{array}{c} (\text{E}\underline{f})\ (\underline{w})\text{: } \tau'\tau T^v A,\ \Gamma\ /\ \rho \\ \updownarrow \\ (\text{E}\underline{f})\ (\underline{w})\text{: } \tau T^v A,\ \Gamma\ /\ \rho \end{array}$$

– An *absolute deduction* (in S5DC) of β from α is a sequence of data, starting with α and ending with β, where each datum after the first is immediately inferable from the immediately preceding one by one of the above rules. A *world-relative deduction* of β from α is an absolute deduction of $T^w[\beta]$ from $T^w[\alpha]$, where w is any world-variable not occurring in either α or β. A *deduction* of β from α is something which is either an absolute or a world-relative deduction of β from α.

– The Herbrand theorems from classical predicate logic case can easily be adapted to the present modal propositional case; as can the results and questions concerning completeness etc. (Indeed the Herbrand theorems can no doubt be strengthened here so that only finitely many instances need to be considered, so that e.g. an attempted analytic refutation will, after a finite mechanical construction, either succeed as refutation or yield countermodel. This should be similar to monadic classical predicate logic. But I have not worked out the details.)

## 3. Examples of deductions

□p ⊢$^N$ p

T$^w$□p                          ↓ □

(v): T$^w$T$^v$p                 ↓ T

(v): T$^v$p                      ↓ Inst

T$^w$p.

p ⊢$^N$ ◇p

T$^w$p                           ↑ Gen

(Em): T$^m$p                     ↑ T

(Em): T$^w$T$^m$p                ↑ ◇

T$^w$◇p.

□p ⊢$^N$ ◇p

□p                              ↓ □

(w): T$^w$p                      ↓ Inst

T$^w$p                           ↑ Gen

(Em): T$^m$p                     ↑ ◇

◇p.

□(p ∧ q) ⊢$^N$ ◇(q ∨ r)

□(p ∧ q)                        ↓ □, ↓ ∧

(w): T$^w$p, T$^w$q              ↓ Del

(w): T$^w$q                      ↓ Inst

T$^w$q                           ↑ Gen

(Em): T$^m$q                     ↑ Exp

(Em): T$^m$q / T$^m$r $\qquad$ ↑ ∨

(Em): T$^m$(q ∨ r) $\qquad$ ↑ ◇

◇(q ∨ r).

<u>□(p → q) ⊢$^N$ □p → □q</u>

□(¬p ∨ q) $\qquad\qquad$ ↓ □, ↓ ∨

(w): T$^w$¬p / T$^w$q $\qquad\qquad$ ↑ Gen

(Ef) (w): T$^{f(w)}$¬p / T$^w$q $\qquad$ ↑ ¬□

(w): ¬□p / T$^w$q $\qquad\qquad$ ↑ □

¬□p / □q $\qquad\qquad\qquad$ ↑ ∨

¬□p ∨ □q.

<u>◇p ⊢ □◇p</u>

◇p $\qquad\qquad\qquad$ ↓ ◇

(Em): T$^m$p $\qquad\qquad$ ↑ IVU

(Em) (w): T$^m$p $\qquad\qquad$ ↑ T

(Em) (w): T$^w$T$^m$p $\qquad\qquad$ ↑ Gen

(Em, f) (w): T$^w$T$^{f(w)}$p $\qquad$ ↓ EVE

(Ef) (w): T$^w$T$^{f(w)}$p $\qquad$ ↑ ◇

(w): T$^w$◇p $\qquad\qquad$ ↑ □

□◇p.

<u>◇□p ⊢$^N$ □p</u>

◇□p $\qquad\qquad\qquad$ ↓ ◇

(Em): T$^m$□p $\qquad\qquad$ ↓ □

(Em) (w): $T^mT^wp$             $\downarrow$ T

(Em) (w): $T^wp$                $\downarrow$ EVE

(w): $T^wp$                  $\uparrow$ □

□p.

$\vdash$ □p $\rightarrow$ ◇p

⊤                      $\uparrow$ EM

¬□p / □p             $\downarrow$ □

(w): ¬□p / $T^wp$        $\downarrow$ Inst, $\uparrow$ Gen

(Em): ¬□p / $T^mp$      $\uparrow$ ◇, $\uparrow$ $\vee$

¬□p $\vee$ ◇p.

$\vdash$ □p $\rightarrow$ p

$T^w[⊤]$ (= ⊤)           $\uparrow$ EM

$T^w$¬□p / $T^w$□p       $\downarrow$ □

(v): $T^w$¬□p / $T^wT^vp$    $\downarrow$ T

(v): $T^w$¬□p / $T^vp$      $\downarrow$ Inst

$T^w$¬□p / $T^wp$        $\uparrow$ $\vee$

$T^w$¬□p $\vee$ p.

□p, ¬p $\vdash^N$ ⊥

$T^w$□p, $T^w$¬p         $\downarrow$ □

(v): $T^wT^vp$, $T^w$¬p     $\downarrow$ T, $\downarrow$ Inst

$T^wp$, $T^w$¬p          $\downarrow$ NC

⊥ (= $T^w[⊥]$).

<u>□(p ∧ ◇¬p) ⊢ ⊥</u>

□(p ∧ ◇¬p)                                     ↓□, ↓∧

(w): T$^w$p, T$^w$◇¬p                           ↓◇, ↓T

(Ef) (w): T$^w$p, T$^{f(w)}$¬p                  ↑↓ GenInst

(Ef): T$^w$p, T$^{f(w)}$¬p, T$^{f(w)}$p, T$^{ff(w)}$¬p    ↓NC, ↓EVE

⊥.

<u>□(p ∧ ◇q) ⊢$^{N, Invertible}$ □p ∧ ◇q</u>

□(p ∧ ◇q)                          ↓□

(w): T$^w$(p ∧ ◇q)                  ↓∧

(w): T$^w$p, T$^w$◇q                ↓◇

(Ef) (w): T$^w$p, T$^w$T$^{f(w)}$q  ↓T

(Ef) (w): T$^w$p, T$^{f(w)}$q       ↑◇

(w): T$^w$p, ◇q                     ↑□

□p, ◇q                              ↑∧

□p ∧ ◇q.

<u>p ∧ ◇¬p ⊢$^N$ ◇(¬p ∧ ◇p)</u>

T$^w$(p ∧ ◇¬p)                      ↓∧

T$^w$p, T$^w$◇¬p                    ↓◇

(Em): T$^w$p, T$^w$T$^m$¬p          ↓T

(Em): T$^w$p, T$^m$¬p              ↑Gen

(Em, n): T$^n$p, T$^m$¬p           ↑T

(Em, n): T$^m$T$^n$p, T$^m$¬p      ↑◇

(Em): $T^m \Diamond p$, $T^m \neg p$        ↑ ∧

(Em): $T^m(\neg p \wedge \Diamond p)$        ↑ T

(Em): $T^w T^m(\neg p \wedge \Diamond p)$        ↑ ◇

$T^w \Diamond(\neg p \wedge \Diamond p)$.

### 4. Directional Deduction for modal predicate logic: Basic notions

We extend now Directional Deduction to 'modal predicate logic', i.e. quantified S5 with constant domains (and no identity in the language). As it happens it is enough to just put together in obvious way our earlier calculi for classical predicate logic and modal propositional logic. We proceed then very briskly here to avoid boredom.

The *formulas* are built with primitive connectives and quantifiers ¬, ∧, ∨, □, ◇, ∀, ∃ in usual way. A *datum* is now an expression

$$\text{(E}\underline{f}\text{) (}\underline{w}\text{, }\underline{x}\text{): } \Gamma \text{ / } \Delta \text{ / } \ldots$$

where Γ, Δ, … are (finitely many) finite sets of qualified formulas, and where now the universal prefix can contain both world-variables and individual variables, and the functional variables in the existential prefix correspond to functions from worlds and/or individuals to either worlds or individuals. (The *value* of such a function is either always a world or always an individual; also what the function takes at a *given* argument-place is always world or always individual; but it may be world at one argument-place and individual at another.)

Similarly a *dual datum* is an expression

$$\text{(}\underline{f}\text{) (E}\underline{w}\text{, }\underline{x}\text{): } \Gamma \text{; } \Delta \text{; } \ldots$$

and a *consecution* an expression

$$\text{(}\underline{f}\text{) (}\underline{g}\text{) (E}\underline{w}\text{, }\underline{x}\text{) (E}\underline{v}\text{, }\underline{y}\text{): } \alpha \Rightarrow \delta$$

where α is datum-matrix and δ dual-datum--matrix.

## 5. The calculi S5QDC, S5QDDC, S5QHC, S5QCC

Again, it will be enough to list the rules of S5QDC (Quantified Modal Data Calculus); from this the dual and hybrid and consecution calculi are defined as before.

*Rules of S5QDC*:

<u>G-rules:</u> $\updownarrow\wedge$, $\updownarrow\vee$, $\updownarrow\neg\neg$, $\updownarrow\neg\wedge$, $\updownarrow\neg\vee$, $\updownarrow\forall$, $\updownarrow\exists$, $\updownarrow\neg\forall$, $\updownarrow\neg\exists$, $\updownarrow\square$, $\updownarrow\Diamond$, $\updownarrow\neg\square$, $\updownarrow\neg\Diamond$ (with the obvious formulations).

<u>N-rules:</u> $\downarrow$NC, $\uparrow$EM (obvious formulations).

<u>S-rules:</u> $\downarrow$Del, $\uparrow$Exp, $\downarrow$Inst, $\uparrow$Gen, $\downarrow$EVE, $\uparrow$IVU, $\updownarrow$T.

– Again, deductions can be 'absolute' or 'world-relative'. – Herbrand theorems and completeness results (and questions) are no doubt straightforward extensions of classical predicate logic case.

*Question* (Datum Interpolation). Is it the case that: $\forall$ formulas A, B: If A $\vDash$ B, or equivalently $T^wA \vDash T^wB$, then $\exists$ datum $\gamma$ with $L(\gamma) \subseteq L(T^wA) \cap L(T^wB)$ and $T^wA \vDash \gamma \vDash T^wB$. (Perhaps we can just reduce $T^wA$ to elementary datum by $\downarrow$ G-rules and then delete all qualified literals with predicate-variables $\notin L(B)$?) (See below the last of our examples of deductions.)

## 6. Examples of deductions

$\forall x\Box Px \vdash^{N,\ Inv.} \Box\forall xPx$

| | |
|---|---|
| $\forall x\Box Px$ | $\downarrow \forall$ |
| $(x): \Box Px$ | $\downarrow \Box$ |
| $(x, w): T^wPx$ | $\uparrow \forall$ |
| $(w): T^w\forall xPx$ | $\uparrow \Box$ |
| $\Box\forall xPx.$ | |

$\exists x\Diamond Px \vdash^{N,\ Inv.} \Diamond\exists xPx$

| | |
|---|---|
| $\exists x\Diamond Px$ | $\downarrow \exists$ |
| $(Ea): \Diamond Pa$ | $\downarrow \Diamond$ |
| $(Ea, m): T^mPa$ | $\uparrow \exists$ |
| $(Em): T^m\exists xPx$ | $\uparrow \Diamond$ |
| $\Diamond\exists xPx.$ | |

$\forall x\Box Px \vdash^N \exists x\Diamond Px$

| | |
|---|---|
| $\forall x\Box Px$ | $\downarrow \forall, \downarrow \Box$ |
| $(x, w): T^wPx$ | $\downarrow$ Inst |
| $T^wPx$ | $\uparrow$ Gen |
| $(Ea, m): T^mPa$ | $\uparrow \Diamond$ |
| $(Ea): \Diamond Pa$ | $\uparrow \exists$ |
| $\exists x\Diamond Px.$ | |

$\Box\forall x(Px \rightarrow Qx), \Diamond\exists xPx \vdash^N \Diamond\exists xQx$

| | |
|---|---|
| $\Box\forall x(\neg Px \vee Qx), \Diamond\exists xPx$ | $\downarrow \Diamond, \downarrow \exists$ |

(Em, a): □∀x(¬Px ∨ Qx), T$^m$Pa                           ↓ □, ↓ ∀

(Em, a) (w, x): T$^w$(¬Px ∨ Qx), T$^m$Pa                   ↓ ∨

(Em, a) (w, x): T$^w$¬Px, T$^m$Pa / T$^w$Qx, T$^m$Pa       ↓ Inst

(Em, a): T$^m$¬Pa, T$^m$Pa / T$^m$Qa, T$^m$Pa              ↓ NC, ↓ Del

(Em, a): T$^m$Qa                                           ↑ ∃

(Em): T$^m$∃xQx                                            ↑ ◇

◇∃xQx.

The following example corresponds to the counterexample to Interpolation for quantified S5 with constant domains in Fine 1979. – We use Δ(Q) as short for ∀x(□Qx ∨ ¬◇Qx), and P ↔ R as short for ∀x(Px ↔ Rx).

Note that, abbreviating the l.h.s. formula below as φ(P, Q) and the r.h.s. one as ψ(P, R), we have (in second-order quantified S5) ∃Q φ(P, Q) ≃ ∀R ψ(P, R). So a putative interpolant θ(P) would have to follow from ∃Q φ(P, Q), and imply ∀R ψ(P, R), and hence be equivalent to them. I.e. θ(P) would 'say' that 'The extension of P might have been disjoint from its true extension'. – Now, although there can be no such *formula* θ(P), it is easy to give a *datum*: (Em) (x): T$^m$¬Px / ¬Px. It is this datum, or rather in effect T$^w$[it], that appears as the critical midpoint in the deduction below.

Δ(Q) ∧ ∀x(Px ↔ Qx) ∧ ◇¬∃x(Px ∧ Qx)  ⊢  Δ(R) ∧ ∀x(Px ↔ Rx) **.→.** ◇¬∃x(Px ∧ Rx)

T$^w$ Δ(Q) ∧ ∀x(Px ↔ Qx) ∧ ◇¬∃x(Px ∧ Qx)                  ↓ ∧

T$^w$Δ(Q), T$^w$∀x(Px ↔ Qx), T$^w$◇¬∃x(Px ∧ Qx)            ↓ ◇, ↓ T

(Em): T$^m$¬∃x(Px ∧ Qx), T$^w$∀x(Px ↔ Qx), T$^w$ΔQ        ↓ ¬∃

(Em) (x): T$^m$¬(Px ∧ Qx), T$^w$∀x(Px ↔ Qx), T$^w$ΔQ      ↓ ∀ (×2)

(Em) (x, y, z): T$^m$¬(Px ∧ Qx), T$^w$ Py ↔ Qy, T$^w$ □Qz ∨ ¬◇Qz   ↓ ¬∧

(Em) (x, y, z): T$^m$¬Px, T$^w$ Py ↔ Qy, T$^w$ □Qz ∨ ¬◇Qz / T$^m$¬Qx, T$^w$ Py ↔ Qy, T$^w$ □Qz ∨ ¬◇Qz                                          ↓ Inst, ↓ Del

(Em) (x): $T^m \neg Px$ / $T^m \neg Qx$, $T^w Px \leftrightarrow Qx$, $T^w \Box Qx \lor \neg \Diamond Qx$         $\downarrow \lor$

(Em) (x): $T^m \neg Px$ / $T^w \Box Qx$, $T^m \neg Qx$, $T^w Px \leftrightarrow Qx$ / $T^w \neg \Diamond Qx$, $T^m \neg Qx$, $T^w Px \leftrightarrow Qx$

                                                 $\downarrow \Box$, $\downarrow T$, $\downarrow$ Inst, $\downarrow$ NC

(Em) (x): $T^m \neg Px$ / $T^w \neg \Diamond Qx$, $T^m \neg Qx$, $T^w Px \leftrightarrow Qx$      $\downarrow \lor$, $\downarrow \land$, $\downarrow \neg \Diamond$, $\downarrow$ Inst, $\downarrow$ NC

(Em) (x): $T^m \neg Px$ / $T^w \neg \Diamond Qx$, $T^m \neg Qx$, $T^w \neg Px$, $T^w \neg Qx$      $\downarrow$ Del

(Em) (x): $T^m \neg Px$ / $T^w \neg Px$                        $\uparrow$ EM, $\downarrow$ Del

(Em) (x): $T^m \neg Px$ / $T^w \neg Px$, $\Delta R$ / $\neg \Delta R$            $\uparrow$ EM, $\downarrow$ Del

(Em) (x): $T^m \neg Px$ / $T^w \neg Px$, $\Delta R$, $T^w P \leftrightarrow R$ / $T^w P \nleftrightarrow R$ / $\neg \Delta R$

                                            $\downarrow \forall$, $\downarrow$ Inst etc., $\downarrow$ NC, $\downarrow$ Del

(Em) (x): $T^m \neg Px$ / $T^w \neg Rx$, $\Delta R$ / $T^w P \nleftrightarrow R$ / $\neg \Delta R$      $\downarrow \forall$, $\downarrow$ Inst

(Em) (x): $T^m \neg Px$ / $T^w \neg Rx$, $\Box Rx \lor \neg \Diamond Rx$ / $T^w P \nleftrightarrow R$ / $\neg \Delta R$

                                            $\downarrow \lor$, $\downarrow \Box$ etc., $\downarrow$ NC, $\downarrow$ Del

(Em) (x): $T^m \neg Px$ / $\neg \Diamond Rx$ / $T^w P \nleftrightarrow R$ / $\neg \Delta R$      $\downarrow \neg \Diamond$ etc.

(Em) (x): $T^m \neg Px$ / $T^m \neg Rx$ / $T^w P \nleftrightarrow R$ / $\neg \Delta R$      $\uparrow \neg \land$, $\uparrow \neg \exists$, $\uparrow \Diamond$

$\Diamond \neg \exists x (Px \land Rx)$ / $T^w P \nleftrightarrow R$ / $\neg \Delta R$              etc.

$T^w \Delta(R) \land \forall x (Px \leftrightarrow Rx)$ **.→.** $\Diamond \neg \exists x (Px \land Rx)$.

DIRECTIONAL DEDUCTION FOR RELEVANCE LOGICS: BASIC SYSTEMS

## 1. Introduction

We consider now directional deduction calculi corresponding to some relevantistic subsystems of Classical Propositional Logic; many of these are equivalential, so in fact subsystems of EqCPL.

The systems include some which will perhaps be already known to the reader: the Belnap–Dunn logic of First-Degree Entailment or FDE, and its equivalential correlate EqFDE; R. B. Angell's equivalential logic which is a proper subsystem of EqFDE (e.g. p and p ∧ (p ∨ q) are EqFDE equivalent but not equivalent in Angell's logic); Correia's 'logic of factual equivalence' which is a proper subsystem of Angell's logic (e.g. p ∨ (q ∧ r) and (p ∨ q) ∧ (p ∨ r) are equivalent in Angell's logic but not in Correia's); and the logic of 'exact entailment' in Fine's sense, and its equivalential correlate which is a proper subsystem of Correia's logic (e.g. p and p ∧ p are equivalent in Correia's logic but not in this Finean logic of mutual exact entailment). (See the Historical Notes at the end of this section for references.)

In the light of Fine's recent work on truthmaker semantics, it becomes clear that these systems (despite their sometimes heterogeneous historical origins) can be regarded as typical representatives in a very natural family of systems – or better very natural class of families of systems.

(I will assume the reader is at least vaguely familiar with the basic ideas of truthmaker semantics. Otherwise one may wish to read first e.g. the survey-paper Fine 2017b before proceeding here.)

There will be fourteen systems in all, divided into four families (of 2 + 4 + 4 + 4 members). Each family corresponds to a type of *concept of truthmaking* (and falsity-making). And *within* each family the different members correspond to different *semantic relations* defined (in uniform manner) from the given concepts of truthmaking and falsity-making.

The four concepts of truthmaking (and falsity-making) are:

*Exact* truthmaker: the 'bare' notion of truthmaker, without any closure conditions.

*Semiregular* truthmaker: with closure under *fusions* of truthmakers.

*Regular* truthmaker: with closure under *fusions* of truthmakers and *intermediates* between truthmakers.

*Isotone* truthmaker: with closure under *superstates* (i.e. an isotone truthmaker is a state *including* an exact truthmaker).

And the different semantic relations are (say for data):

*Truthmaker equivalence*: (W.r.t. every model:) The truthmakers of α are precisely the truthmakers of β.

*Falsity-maker equivalence*: The falsity-makers of α are precisely the falsity-makers of β.

*Entailment*: Every truthmaker of α is truthmaker of β.

*Containment*: Every falsity-maker of β is falsity-maker of α.

– Truthmaker equivalence is of course mutual entailment, and falsity-maker equivalence mutual containment. The Containment system is (in each case) dual to the Entailment system, and the Falsity-maker equivalence system dual to the Truthmaker equivalence system.

Isotone Entailment happens to coincide in extension with Isotone Containment; i.e. we have here a 'self-dual' system – which is none other than the familiar FDE. Consequently Isotone Truthmaker equivalence also coincides in extension with Isotone Falsity-maker equivalence. So in this 'isotone family' we have only two systems: Isotone

Entailment = Isotone Containment = FDE; and Isotone Truthmaker equivalence = Isotone Falsity-maker equivalence = EqFDE.

In each of the other three cases (Regular, Semiregular, Exact) we get four extensionally distinct systems. (As far as data and dual data are concerned. Restricting attention to *formulas* we get a further self-duality collapse with Regular Truthmaker equivalence and Regular Falsity-maker equivalence [*not* with Regular Entailment and Regular Containment]. See below the section 'Note on duality'.) Our fourteen systems can then be represented in the following table. – We name the systems in the obvious way from the initials of Isotone / Regular / Semiregular / Exact and then Entailment / Containment / Truthmaker equivalence / Falsity-maker equivalence. (But for FDE we stick with this traditional name, and then use EqFDE for the corresponding equivalential system.)

| Truthmaking & falsity-making | ⊨ | ≃ |
|---|---|---|
| Isotone | FDE (= IE = IC) | EqFDE (= IT = IF) |
| Regular | RE | RT (= EqRE = Angell's log.) |
| | RC | RF (= EqRC) |
| Semiregular | SE | ST (= EqSE = Correia's 'log. of factual equiv.') |
| | SC | SF (= EqSC) |
| Exact | EE (= Fine's log. of exact ent.) | ET (= EqEE) |
| | EC | EF (= EqEC) |

Except for the case of the isotone family, the Entailment systems have a *synthetic* character, and the Containment systems an *analytic* character – we have: p entails p ∨ q; p ∧ q does *not* entail p; p ∧ q contains p; p does *not* contain p ∨ q. But note that this analyticity and syntheticity are now w.r.t. a thick notion of 'truthmaker content', not ultra-fine structured hyperintensional content as before. Thus e.g. between p and ¬¬p there is always both mutual entailment and mutual containment.

*Historical notes.* (1) Belnap 1959 and Anderson & Belnap 1962 give a simple syntactic criterion for FDE consequence, which amounts to a formulaic formulation of our Hybrid Calculus for FDE given below (i.e. formulation in terms of formulas only rather than data and dual data – essentially, one reduces premiss-formula to a corresponding DNF and conclusion-formula to corresponding CNF, rather than premiss-datum to corresponding elementary datum and conclusion dual datum to corresponding elementary dual datum). The paper van Fraassen 1969 gives a striking truthmaker semantics for FDE (which we will use here, in somewhat modified form). There is also a better-known four-valued semantics for FDE, first given in print in Dunn 1976 and Belnap 1977a,b; see also Priest 2008 ch. 8 for a snappy exposition.

(2) RT originates in the work of R. B. Angell (1977, 1989). The paper Fine 2016 gives a very ingenious truthmaker semantics for the system (again, what we will use here, in somewhat modified form). Fine gives also (ibid., § 10) a many-valued semantics. See also Correia 2004 for related work.

(3) ST originates with Correia 2016, where inter alia a natural Fine-style truthmaker semantics is given (again, what we will use here, in somewhat modified form).

(4) EE is mentioned briefly by Fine in various places. Also the recent papers Fine and Jago 2019 and Krämer forthcoming are particularly relevant here.

## 2. The calculi FDEDC, RTDC, STDC, ETDC, etc.

(1) *The calculi FDEDC, FDEDDC, FDEHC, FDECC.* Each of these is obtained from our corresponding calculus for classical propositional logic (DC, DDC, HC, CC) simply by omission of the NC and EM rules.

(2) *The calculi EqFDEDC and EqFDEDDC.* EqFDEDC is the (data) calculus whose primitive rules are the G-rules and the following rules for SuperSet inclusion ($\uparrow$) and exclusion ($\downarrow$):

$$\Gamma \,/\, \Gamma, \Delta \,/\, \rho$$
$$(\updownarrow SS) \qquad \updownarrow$$
$$\Gamma \,/\, \rho$$

Note that from these $\updownarrow$SS rules we can immediately derive what we may call the Mutual Inclusion rule or MutInc: Proceed from $\alpha$ to $\beta$ provided every component of $\alpha$ includes some component of $\beta$ and vice versa. For under these conditions, and writing $\alpha$ as $\Gamma_1 \,/\, \Gamma_2 \,/\, \ldots$ and $\beta$ as $\Delta_1 \,/\, \Delta_2 \,/\, \ldots$, we have:

$\Gamma_1 \,/\, \Gamma_2 \,/\, \ldots \qquad \uparrow$ SS  (since every $\Delta$ is superset of some $\Gamma$)

$\Gamma_1 \,/\, \Gamma_2 \,/\, \ldots \,/\, \Delta_1 \,/\, \Delta_2 \,/\, \ldots \qquad \downarrow$ SS  (since every $\Gamma$ is superset of some $\Delta$)

$\Delta_1 \,/\, \Delta_2 \,/\, \ldots$

– There is then the similar (dual data) calculus EqFDEDDC, with primitive rules the G-rules and $\updownarrow$SS in the dual data form

$$\Gamma \,;\, \Gamma \setminus \Delta \,;\, \rho$$
$$(\updownarrow SS) \qquad \updownarrow$$
$$\Gamma \,;\, \rho$$

(3) *The calculi RTDC, RFDDC, RTHC, RTCC.*

(i) *The calculus RTDC.* We omit (from the basic rules of DC) NC, EM, Del, and Exp (thus all basic rules other than the G-rules); and add the following new S-rules of Fusion-Inclusion (↑Fus), Fusion-Exclusion (↓Fus), Intermediate-Inclusion (↑Int), and Intermediate-Exclusion (↓Int):

$$\Gamma / \Delta / \Gamma, \Delta / \rho \qquad\qquad \Gamma / \Gamma, \Delta / \Gamma, \Delta, \Theta / \rho$$

(↕Fus)    ↕         (↕Int)    ↕

$$\Gamma / \Delta / \rho \qquad\qquad\qquad \Gamma / \Gamma, \Delta, \Theta / \rho$$

*Remark.* In DC (indeed even FDEDC) these rules are readily derivable by ↑Exp and ↓Del. Ditto in EqFDEDC by ↕SS. – The Fusion rules are also derivable by G-rules alone:

$\Gamma / \Delta / \rho \qquad ↑ ∧$

$∧\Gamma / ∧\Delta / \rho \qquad ↑ ∨$

$∧\Gamma ∨ ∧\Delta, ∧\Gamma ∨ ∧\Delta / \rho \qquad\qquad ↓ ∨$

$∧\Gamma ∨ ∧\Delta, ∧\Gamma / ∧\Gamma ∨ ∧\Delta, ∧\Delta / \rho \qquad\qquad ↓ ∨$

$∧\Gamma / ∧\Gamma, ∧\Delta / ∧\Delta / \rho \qquad\qquad ↓ ∧$

$\Gamma / \Delta / \Gamma, \Delta / \rho.$

(This is of course invertible, giving then derivation also of ↓Fus.) These derivations are however doubly anomalous: connectives intrude in derivations of purely structural rules; and the derivations are synthetico-analytic whereas Fusion-Inclusion is ideally regarded as purely synthetic rule and Fusion-Exclusion as purely analytic. So I prefer to count ↑Fus and ↓Fus as primitive S-rules of RTDC.

(ii) *The calculus RFDDC.* This is the calculus for dual data with primitive rules the exact duals of those of RTDC. I.e. we have the G-rules for dual data (as in DDC) and the dual data versions of ↕Fus and ↕Int:

$$\Gamma\,;\Delta\,;\Gamma\setminus\Delta\,;\rho \qquad\qquad\qquad \Gamma\,;\Gamma\setminus\Delta\,;\Gamma\setminus\Delta\setminus\Theta\,;\rho$$

($\updownarrow$Fus) $\qquad\qquad\quad \updownarrow \qquad\qquad$ ($\updownarrow$Int) $\qquad\qquad\quad \updownarrow$

$$\Gamma\,;\Delta\,;\rho \qquad\qquad\qquad\qquad\qquad \Gamma\,;\Gamma\setminus\Delta\setminus\Theta\,;\rho$$

(iii) *The calculus RTHC*. We consider the formulation of HC with the Distribution rule Distr rather than Criss-Crossing as the Transition-rule; and from this we now omit NC, EM, Del, and Exp (in each case of course both the rule for data and that for dual data); and add as new S-rules $\updownarrow$Fus and $\updownarrow$Int for both data and dual data. Thus an RTHC deduction is a sequence $\langle \alpha \ldots \beta, \delta \ldots \varepsilon \rangle$ where $\langle \alpha \ldots \beta \rangle$ is an RTDC deduction, $\langle \delta \ldots \varepsilon \rangle$ is an RFDDC deduction, and $\delta$ is obtainable from $\beta$ by Distr.

(iv) *The calculus RTCC*. We omit from CC the rules NC, EM, CrCr-L, and CrCr-R; add the S-rules Fus-L, Int-L, Fus-R, Int-R with the obvious formulations (for inclusion of fusions or intermediates in the datum on the left or the dual datum on the right); and now define a *proof* of a consecution as a deduction of it, by the given rules (G-rules as before, plus the new S-rules), from a 'distribution axiom', i.e. consecution $\alpha \Rightarrow \delta$ where $\delta$ is obtainable from $\alpha$ by Distribution.

(4) *The calculi REDC and RCDDC*. REDC is RTDC plus $\uparrow$Exp. RCDDC is RFDDC plus $\downarrow$Del.

(5) *The calculi STDC and SFDDC*. STDC is DC minus the NC, EM, Del and Exp rules – i.e. we have now the G-rules and nothing else. Similarly SFDDC is DDC stripped down to just its G-rules.

(6) *The calculi SEDC and SCDDC*. SEDC is STDC plus $\uparrow$Exp. SCDDC is SFDDC plus $\downarrow$Del.

(7) *The calculi ETDC and EFDDC*. We now change somewhat the very concepts of data and dual data. Instead of finite set of finite sets of formulas as before, the definition now (for both concepts) is: finite set of finite *multiset* of formulas. To avoid confusion with the 'old' data and dual data we will call these resp. *sms data* and *sms dual data*.

Let us say that a finite multiset of formulas $\Gamma_0$ is an *inferior replica* of finite multiset of formulas $\Gamma$, or that $\Gamma$ is a *superior replica* of $\Gamma_0$, if the underlying set of formulas is the same in both cases, and for any formula A in this set, the number of occurrences of A in $\Gamma_0$ is less than or equal to the number of occurrences of A in $\Gamma$.

The rules of ETDC are then the G-rules plus the following structural rules of Inferior Replica inclusion (↑) and exclusion (↓):

(↕IR)  For $\Gamma_0$ inferior replica of $\Gamma$:

$$\Gamma \: / \: \Gamma_0 \: / \: \rho$$
$$\updownarrow$$
$$\Gamma \: / \: \rho$$

And the rules of EFDDC are the G-rules (for dual data, now reconceptualized as sms dual data) plus the structural ↕IR rules:

(↕IR)  For $\Gamma_0$ inferior replica of $\Gamma$:

$$\Gamma \: ; \: \Gamma_0 \: ; \: \rho$$
$$\updownarrow$$
$$\Gamma \: ; \: \rho$$

(8) *The calculi EEDC and ECDDC.*  EEDC is ETDC plus ↑Exp. ECDDC is EFDDC plus ↓Del.

# 3. Examples of deductions

Note that, in all the calculi ETDC, STDC, RTDC, EqFDEDC, EqDC (and also in the dual data calculi EFDDC, SFDDC, RFDDC, EqFDEDDC, EqDDC) – in all these calculi *all* the primitive rules are 'two-directional'. So here a deduction witnessing $\alpha \vdash \beta$ (or witnessing $\beta \vdash \alpha$) is always automatically invertible and so in effect witnesses the mutual deducibility $\alpha \dashv\vdash \beta$.

We have here:

$$\dashv\vdash_{\text{ETDC}} \subset \dashv\vdash_{\text{STDC}} \subset \dashv\vdash_{\text{RTDC}} \subset \dashv\vdash_{\text{EqFDEDC}} \subset \dashv\vdash_{\text{EqDC}}.$$

(And dually:

$$\dashv\vdash_{\text{EFDDC}} \subset \dashv\vdash_{\text{SFDDC}} \subset \dashv\vdash_{\text{RFDDC}} \subset \dashv\vdash_{\text{EqFDEDDC}} \subset \dashv\vdash_{\text{EqDDC}}.)$$

We have chosen the examples below so that, w.r.t. this series, the (inter)deducibility *first* holds in the calculus in question (thus does *not* hold in the previous ones if any).

$\underline{p \ \dashv\vdash_{\text{ETDC}} \ p \vee p}$

$p \quad$ (i.e. $\{[p]\} = \{[p], [p]\} = p \mathbin{/} p$) $\qquad \Uparrow \vee$

$p \vee p.$

(Dually: $p \dashv\vdash_{\text{EFDDC}} p \wedge p$. Now we have: $\{[p], [p]\} = p \mathbin{;} p.$)

$\underline{(p \wedge p) \vee p \ \dashv\vdash_{\text{ETDC}} \ p \wedge p}$

$(p \wedge p) \vee p \qquad \Downarrow \vee, \ \Downarrow \wedge$

$p, p \mathbin{/} p \qquad\qquad \Downarrow \text{IR}$

$p, p \qquad\qquad \Uparrow \wedge$

$p \wedge p.$

(Dually: $(p \lor p) \land p \dashv\vdash_{EFDDC} p \lor p$.)

$\underline{p \land (q \lor r) \dashv\vdash_{ETDC} (p \land q) \lor (p \land r)}$

$p \land (q \lor r)$       $\downarrow \land$

$p, q \lor r$         $\downarrow \lor$

$p, q \,/\, p, r$      $\uparrow \land,\ \uparrow \lor$

$(p \land q) \lor (p \land r)$.

(Dually: $p \lor (q \land r) \dashv\vdash_{EFDDC} (p \lor q) \land (p \lor r)$.)

$\underline{p \dashv\vdash_{STDC} p \land p}$

$p$    $(= p, p)$      $\uparrow \land$

$p \land p$.

(Dually: $p \dashv\vdash_{SFDDC} p \lor p$.)

$\underline{p \lor (q \land r) \dashv\vdash_{RTDC} (p \lor q) \land (p \lor r)}$

$(p \lor q) \land (p \lor r)$        $\downarrow \land$

$p \lor q, p \lor r$         $\downarrow \lor$

$p, p \lor r \,/\, q, p \lor r$     $\downarrow \lor\ (\times 2)$

$p \,/\, p, r \,/\, p, q \,/\, q, r$     $\uparrow$ Fus

$p \,/\, p, r \,/\, p, q \,/\, q, r \,/\, p, q, r$     $\downarrow$ Int $(\times 2)$

$p \,/\, q, r \,/\, p, q, r$       $\downarrow$ Fus

$p \,/\, q, r$            $\uparrow \land,\ \uparrow \lor$

$p \lor (q \land r)$.

(Dually: p ∧ (q ∨ r) ⊣⊢_RFDDC (p ∧ q) ∨ (p ∧ r), and this is the 'first time' we have this interdeducibility, w.r.t. the dual series above.)

p  ⊣⊢_EqFDEDC  p ∨ (p ∧ q)

p                    ↑ SS

p / p, q             ↑ ∧, ↑ ∨

p ∨ (p ∧ q).

p  ⊣⊢_EqFDEDC  p ∧ (p ∨ q)

p ∧ (p ∨ q)          ↓ ∧

p, p ∨ q             ↓ ∨

p / p, q             ↓ SS

p.

p ∧ ¬p  ⊣⊢_EqDC  q ∧ ¬q

p ∧ ¬p          ↓ ∧

p, ¬p           ↓ NC

{ }             ↑ NC

q, ¬q           ↑ ∧

q ∧ ¬q.

# 4. Truthmaker semantics for FDE, ET, etc.

By an *elementary situation space* (e.s.s.) we mean simply any set. We call the elements of an e.s.s. S *elementary situations*. A *state* (w.r.t. e.s.s. S) is a set of elementary situations (i.e. subset of S). (Intuitively, we think of such elementary situations as atomic situations and negations thereof, and of states as conjunctions of elementary situations. Here I am slightly departing from the standard Fine-style construction; see item (1.1) in the section 'Variations in truthmaker semantics' below for discussion.) A *bilateral proposition* is an ordered pair $\langle X, Y \rangle$ where X and Y are sets of states with $X \neq Y$. The states in X (resp., Y) are called the *truthmakers* (resp., *falsity-makers*) of the bilateral proposition $\langle X, Y \rangle$.

A *model* is a pair $(S, \sigma)$ where S is an e.s.s. and $\sigma$ a function assigning a bilateral proposition (w.r.t. S) to each propositional variable. The function $\sigma$ can then be extended to assignment of bilateral propositions to all formulas by the following clauses (we write $\sigma^+(A)$ and $\sigma^-(A)$ for resp. the first and second terms of the bilateral proposition $\sigma(A)$):

$$\sigma^+(\neg A) \ = \ \sigma^-(A).$$

$$\sigma^-(\neg A) \ = \ \sigma^+(A).$$

$$\sigma^+(A \wedge B) \ = \ \{\Sigma_1 \cup \Sigma_2 : \Sigma_1 \in \sigma^+(A), \Sigma_2 \in \sigma^+(B)\}.$$

$$\sigma^-(A \wedge B) \ = \ \sigma^-(A) \cup \sigma^-(B).$$

$$\sigma^+(A \vee B) \ = \ \sigma^+(A) \cup \sigma^+(B).$$

$$\sigma^-(A \vee B) \ = \ \{\Sigma_1 \cup \Sigma_2 : \Sigma_1 \in \sigma^-(A), \Sigma_2 \in \sigma^-(B)\}.$$

Exactly similar clauses serve also to extend $\sigma$ to data expressions and dual data expressions. We have then in particular for $\top$ and $\bot$ (whether as data expressions or dual data expressions): $\forall$ model $(S, \sigma)$: $\sigma(\top) = \langle \{\varnothing\}, \{\ \}\rangle$ and $\sigma(\bot) = \langle \{\ \}, \{\varnothing\}\rangle$. For data expressions and dual data expressions which have characteristic formulas, what $\sigma$ gives to the datum expression or dual datum expression is what $\sigma$ gives to its characteristic formula. – All this applies also of course to sms data expressions and sms dual data expressions, which indeed are the same expressions as resp. data expressions and dual data expressions.

We may also consider interpretation $\sigma$ as extended to data, dual data etc. themselves, as opposed to data expressions etc.: thus say $\sigma(\alpha)$ where $\alpha$ is a datum is what $\sigma$ gives to the 'canonical' datum expression for $\alpha$ – where component expressions and formulas inside component expressions are disposed in alphabetical order and without multiple occurrences of the same component expression or the same formula inside a single component expression.

We then define FDE consequence for data by:

$\alpha \vDash_{FDE} \beta \;=_{df}\; \forall$ model $(S, \sigma)$: $\forall \Sigma \in \sigma^+(\alpha)\; \exists \Sigma' \in \sigma^+(\beta)$: $\Sigma' \subseteq \Sigma$.

The same definition applies also for FDE consequence with dual data ($\delta \vDash_{FDE} \varepsilon$), or datum and dual datum ($\alpha \vDash_{FDE} \delta$), or for FDE validity of a consecution ($\vDash_{FDE} \alpha \Rightarrow \delta$). Also

$\alpha \simeq_{FDE} \beta$ (or: $\alpha \simeq_{EqFDE} \beta$) $\;=_{df}\; \alpha \vDash_{FDE} \beta$ and $\beta \vDash_{FDE} \alpha$

and similarly for dual data.

Now for sms data $\alpha$, $\beta$ we define:

$\alpha \simeq_{ET} \beta \;=_{df}\; \forall$ model $(S, \sigma)$: $\sigma^+(\alpha) = \sigma^+(\beta)$.

$\alpha \simeq_{EF} \beta \;=_{df}\; \forall$ model $(S, \sigma)$: $\sigma^-(\alpha) = \sigma^-(\beta)$.

$\alpha \vDash_{EE} \beta \;=_{df}\; \forall$ model $(S, \sigma)$: $\sigma^+(\alpha) \subseteq \sigma^+(\beta)$.

$\alpha \vDash_{EC} \beta \;=_{df}\; \forall$ model $(S, \sigma)$: $\sigma^-(\beta) \subseteq \sigma^-(\alpha)$.

And similarly for sms dual data. It should be clear that $\alpha \simeq_{EF} \beta$ iff $d(\alpha) \simeq_{ET} d(\beta)$, and (for sms dual data $\delta$, $\varepsilon$) $\delta \simeq_{EF} \varepsilon$ iff $d(\delta) \simeq_{ET} d(\varepsilon)$. Also $\alpha \vDash_{EE} \beta$ iff $d(\beta) \vDash_{EC} d(\alpha)$, and $\delta \vDash_{EC} \varepsilon$ iff $d(\varepsilon) \vDash_{EE} d(\delta)$.

– The *canonical model* is the model (S, σ) where S is the set of literals and, for each propositional variable p, $\sigma(p) = \langle \{\{p\}\}, \{\{\neg p\}\} \rangle$.

– In addition to the above 'basic' or 'exact' truthmaker semantics for FDE, we will consider also two variants, both still within the general truthmaker semantics framework: the *regular* semantics to be given in the next section, and the *isotone* semantics which we proceed to give now.

For set of states X, we define:

Ss Cl (X) ['the closure of X under superstates'] $=_{df}$ $\{\Sigma : \Sigma' \subseteq \Sigma$ for some $\Sigma' \in$ X$\}$.

We say that set of states X is *closed under superstates*, or *upward closed*, if X = Ss Cl (X). And we say that a bilateral proposition $\langle X, Y \rangle$ is *isotone* if both X and Y are upward closed.

An *isotone model* is a pair $\langle S, \sigma \rangle$ where S is e.s.s. and σ is function assigning an isotone bilateral proposition to each propositional variable. σ is then extended to other formulas and to data etc. exactly as before in the 'basic' semantics. – Note that still all values of this extended σ are isotone propositions (so that sprinkling Ss Cl through the definition of extended σ, in the style of the regular and semiregular semantics in the next section, would not change any values of this extended σ function). This is obvious for negations (of course if $\langle X, Y \rangle$ is isotone then so is $\langle Y, X \rangle$), and for truthmakers of disjunctions and the like. And for the remaining case of truthmakers of conjunctions (and the like) we can argue as follows: Suppose $\Sigma'$ is superstate of a state $\Sigma$ which is fusion of truthmaker $\Sigma_1$ of A with truthmaker $\Sigma_2$ of B; then $\Sigma'$ is superstate of $\Sigma_1$ and so (by hypothesis that $\sigma^+(A)$ is upward closed) $\Sigma'$ is itself truthmaker of A, and so it is a fusion of some truthmaker of A with some truthmaker of B, viz. of $\Sigma'$ with $\Sigma_2$.

We then define:

$\alpha \vDash_{\text{FDEi}} \beta$ $=_{df}$ $\forall$ isotone model (S, σ): $\sigma^+(\alpha) \subseteq \sigma^+(\beta)$.

Also of course

$\alpha \simeq_{\text{FDEi}} \beta =_{\text{df}} \alpha \vDash_{\text{FDEi}} \beta$ and $\beta \vDash_{\text{FDEi}} \alpha$.

This is equivalent to: $\forall$ isotone model $(S, \sigma)$: $\sigma(\alpha) = \sigma(\beta)$.

In the context of these isotone models, the *canonical model* is the model $(S, \sigma)$ where S is the set of literals and, for each propositional variable p, $\sigma(p) = \langle \{\Sigma \subseteq S : p \in \Sigma\}, \{\Sigma \subseteq S : \neg p \in \Sigma\} \rangle$.

## 5. Truthmaker semantics for RT, ST, etc.

For set of states X, we define:

Fus Cl (X) ['the closure of X under fusions'] $=_{df}$ {$\cup$(Y) : Y $\subseteq$ X, Y $\neq \varnothing$}.

Int Cl (X) ['the closure of X under intermediates'] $=_{df}$ {$\Sigma$ : $\exists \Sigma_1, \Sigma_2 \in$ X ($\Sigma_1 \subseteq \Sigma$ $\subseteq \Sigma_2$)}.

And we say that set of states X is *closed under fusions* (resp., *closed under intermediates*) if X = Fus Cl (X) (resp., X = Int Cl (X)).

A bilateral proposition $\langle$X, Y$\rangle$ is *regular* if both X and Y are closed under fusions and closed under intermediates and non-empty.

*Remarks*. (1) Note the non-emptiness requirement. We will have (and 'want' to have here) p $\vee$ (q $\wedge$ r) $\simeq_{RT}$ (p $\vee$ q) $\wedge$ (p $\vee$ r). But if the non-emptiness requirement were dropped (and nothing else changed) this would no longer hold: if say q is $\langle \varnothing$, Y$\rangle$ then fusion of a truthmaker of p with a truthmaker of r will be truthmaker of the r.h.s. but not in general of the l.h.s.

(2) Since, intuitively, a proposition might have no truthmakers (e.g. the empty disjunction), or no falsity-makers (empty conjunction), it may be of interest to investigate the system defined like RT here except that we drop this non-emptiness requirement.

(3) In a calculus corresponding to such a system, the G-rules (and hence $\updownarrow$Fus) are still correct, but *not* so the $\updownarrow$Int rules. E.g. if r has no truthmaker, then the truthmakers of p / p, q, r are just the truthmakers of p, but the truthmakers of p / p, q / p, q, r are the truthmakers of p *and* the truthmakers of p, q.

(4) Indeed suppose we were to add to our language the new formulas $\top$, taken to be made true by the empty state alone and made false by no state, and $\bot$, taken to be made true by no state and made false by the empty state alone. There are then the natural rules:

$$\Gamma, \top / \rho \qquad\qquad\qquad \Gamma, \bot / \rho$$

$$(\updownarrow\top) \quad \updownarrow \qquad\qquad\qquad (\updownarrow\bot) \quad \updownarrow$$

$$\Gamma / \rho \qquad\qquad\qquad\qquad \rho$$

$$\Gamma, \neg\bot / \rho \qquad\qquad\qquad \Gamma, \neg\top / \rho$$

$$(\updownarrow\neg\bot) \quad \updownarrow \qquad\qquad\qquad (\updownarrow\neg\top) \quad \updownarrow$$

$$\Gamma / \rho \qquad\qquad\qquad\qquad \rho$$

(Addition of these rules to those of a complete calculus in the original language will habitually yield a complete calculus in the extended language: in the middle portion of deduction [after analyses and 'reverse analyses' by G-rules] we may use these rules to summarily remove all $\top$'s and $\bot$'s, and then proceed as before.) Then $\updownarrow$SS becomes derivable from $\updownarrow$Int and $\updownarrow\bot$: we have (invertibly):

$$\Gamma / \rho \qquad\qquad\qquad \uparrow \bot$$

$$\Gamma / \Gamma, \Delta, \bot / \rho \qquad\qquad \uparrow \text{Int}$$

$$\Gamma / \Gamma, \Delta / \Gamma, \Delta, \bot / \rho \quad \downarrow \bot$$

$$\Gamma / \Gamma, \Delta / \rho. \ \dashv$$

A *regular model* is a pair $(S, \sigma)$ where S is e.s.s. and $\sigma$ is function assigning a regular bilateral proposition to each propositional variable. $\sigma$ is then extended to negations as before and to conjunctions and disjunctions by:

$$\sigma^+(A \wedge B) \ = \ \text{Int Cl (Fus Cl } (\{\Sigma_1 \cup \Sigma_2 : \Sigma_1 \in \sigma^+(A), \Sigma_2 \in \sigma^+(B)\})).$$

$$\sigma^-(A \wedge B) \ = \ \text{Int Cl (Fus Cl } (\sigma^-(A) \cup \sigma^-(B))).$$

$$\sigma^+(A \vee B) \ = \ \text{Int Cl (Fus Cl } (\sigma^+(A) \cup \sigma^+(B))).$$

$$\sigma^-(A \vee B) \ = \ \text{Int Cl (Fus Cl } (\{\Sigma_1 \cup \Sigma_2 : \Sigma_1 \in \sigma^-(A), \Sigma_2 \in \sigma^-(B)\})).$$

And similarly for data expressions etc.

We then define RT equivalence for data by:

$\alpha \simeq_{RT} \beta \ =_{df} \ \forall$ regular model $(S, \sigma)$: $\sigma^+(\alpha) = \sigma^+(\beta)$.

And the same definition applies for dual data ($\delta \simeq_{RT} \varepsilon$), datum and dual datum ($\alpha \simeq_{RT} \delta$), and for RT validity of consecution ($\vDash_{RT} \alpha \Rightarrow \delta$).

We define RF equivalence for data by:

$\alpha \simeq_{RF} \beta \ =_{df} \ \forall$ regular model $(S, \sigma)$: $\sigma^-(\alpha) = \sigma^-(\beta)$.

And similarly for dual data etc.

*Remarks.* (1) For *formulas* we have: A $\simeq_{RT}$ B iff A $\simeq_{RF}$ B. – This is also almost always, though not quite always, true for *data* (and dual data etc.). Thus if each of $\alpha$ and $\beta$ has 'characteristic formula' (corresponding disjunction of conjunctions), i.e. is not the empty datum and does not contain the empty component, then indeed $\alpha \simeq_{RT} \beta$ iff $\alpha \simeq_{RF} \beta$. However we have e.g.: p / $\varnothing$ always has (w.r.t. every model) the same falsity-makers as q / $\varnothing$ (namely: none), but of course not always the same truthmakers.

(2) Clearly what extended $\sigma$ (from regular model $(S, \sigma)$) gives to datum $\alpha$ is a regular proposition iff $\alpha$ is not the empty datum and does not contain the empty component. ⊣

We define Regular Entailment and Regular Containment by:

$\alpha \vDash_{RE} \beta \ =_{df} \ \forall$ regular model $(S, \sigma)$: $\sigma^+(\alpha) \subseteq \sigma^+(\beta)$.

$\alpha \vDash_{RC} \beta \ =_{df} \ \forall$ regular model $(S, \sigma)$: $\sigma^-(\beta) \subseteq \sigma^-(\alpha)$.

And similarly for dual data etc.

We may also define a 'regular semantics' version of FDE consequence (which as we will see is equivalent to $\vDash_{FDE}$ in the sense of the preceding section):

$\alpha \vDash_{FDEr} \beta \ =_{df} \ \forall$ regular model $(S, \sigma)$: $\forall \Sigma \in \sigma^+(\alpha) \ \exists \Sigma' \in \sigma^+(\beta)$: $\Sigma' \subseteq \Sigma$.

Now, a bilateral proposition $\langle X, Y \rangle$ is *semiregular* if both X and Y are closed under fusions. A *semiregular model* is pair $(S, \sigma)$ where S is e.s.s. and $\sigma$ is function assigning a semiregular bilateral proposition to each propositional variable. Then σ is extended to other formulas and to data expressions etc. as above with regular models only now with Int Cl omitted. We then define ST equivalence and SF equivalence for data by:

$\alpha \simeq_{ST} \beta =_{df} \forall$ semiregular model (S, σ): $\sigma^+(\alpha) = \sigma^+(\beta)$.

$\alpha \simeq_{SF} \beta =_{df} \forall$ semiregular model (S, σ): $\sigma^-(\alpha) = \sigma^-(\beta)$.

And similarly for dual data etc. It should be clear that $\alpha \simeq_{SF} \beta$ iff $d(\alpha) \simeq_{ST} d(\beta)$, and $\delta \simeq_{SF} \varepsilon$ iff $d(\delta) \simeq_{ST} d(\varepsilon)$.

We define Semiregular Entailment and Semiregular Containment by:

$\alpha \vDash_{SE} \beta =_{df} \forall$ semiregular model (S, σ): $\sigma^+(\alpha) \subseteq \sigma^+(\beta)$.

$\alpha \vDash_{SC} \beta =_{df} \forall$ semiregular model (S, σ): $\sigma^-(\beta) \subseteq \sigma^-(\alpha)$.

And similarly for dual data etc.

– Whether in the context of regular or semiregular models, the *canonical model* is defined exactly as before with the basic semantics for FDE.

# 6. Remark on the standard of validity

FDE consequence is a 'proper fragment' of classical consequence (i.e. every case of FDE consequence is a case of classical consequence but not vice versa); and one might then expect a matching notion of 'FDE validity' (for data, or dual data, or formulas) as a corresponding proper restriction of classical validity. – But as it turns out, there does not seem to be any useful such notion of FDE validity.

Notice first that we cannot here define in usual style (say for datum $\alpha$; the situation is similar with dual data, or even just formulas): $\alpha$ is valid $=_{df} \forall$ model M: M verifies $\alpha$. For our semantics above does not really include a notion of a model verifying a formula. In this it is like the familiar possible world semantics for modal logic without designated world. – We could enrich the notion of e.s.s. so that some elementary situations are 'designated' as 'factual'. Then a state would be called factual if all its elements are factual, and a bilateral proposition (which now I suppose should mean pair $\langle X, Y \rangle$ of sets of states s.t. X contains factual state iff Y doesn't) true if some state in its first component is factual, and then finally we might say that M $= (\ldots, \sigma)$ verifies $\alpha$ if $\sigma(\alpha)$ is true. This construction would seem to yield however simply *classical* validity.

The other natural idea would be to try to define FDE validity from the already defined notion of FDE consequence along familiar lines of definition of validity from consequence. The most obvious candidates (for definition of '$\alpha$ is FDE valid') here seem to be (where '$\neg\alpha$' can be understood e.g. as the object resulting from $\alpha$ by insertion of $\neg$ before each formula [element of element] but with this resulting object now considered as a *dual* datum):

(1)  $\top \vDash_{FDE} \alpha$.

(2)  $\forall\beta: \beta \vDash_{FDE} \alpha$.

(3)  $\neg\alpha \vDash \alpha$.

But none of these seems reasonable. (1) and (2) are equivalent, and hold iff { } $\in \alpha$. This is hardly a useful 'standard of validity' – e.g. for no formula A would the corresponding datum {{A}} be valid. (And the same applies to (1′) ∅ ⊨_FDE A with the natural definitions under the more conventional understanding of ⊨ as relating set-of-formulas and formula.) (3) is less trivial but also surely unacceptable: e.g. it would make p ∨ ¬p valid, and q ∨ ¬q valid, but p ∨ ¬p .∧. q ∨ ¬q *not* valid!

– So I think we may conclude that there does not seem to be a useful distinctively relevantistic notion of validity matching the relevantistic notion of FDE consequence. This is perhaps intuitively not very surprising, since after all the idea of relevantistic requirement is that of existence of *connection* between *premiss* and *conclusion* in proper inference; and so says nothing of the 'logical truth' of an individual formula (or other formal object such as datum) per se. (Unless of course that was a formula containing a relevantistic entailment connective or the like – but that is not the case here.) (Cf. the remark in Belnap 1977b, penultimate paragraph on p. 46 of Omori and Wansing 2019.)

So at least for a language with the present resources, a relevantistic philosopher who insists that in the most proper notion of *consequence* there must be connection etc. – such a philosopher can perfectly well nevertheless accept the classical standard of validity. (Not that he has to, but if he doesn't it shouldn't be because of his relevantism.)

It is important here that the intuitive idea of validity is the idea of something that is '*always true*' – and *not* the idea of something that 'follows from the empty set of premisses'. At least, we may distinguish these two intuitive ideas; and then the relevantist as such need not disagree with the claim that a proper explication of the former in the present domain yields classical validity; although he should say that proper explication of the latter yields a trivial notion with empty or nearly empty extension (as in (1)/(1′) and (2) above).

– Essentially the same points apply also to unsatisfiability as much as to validity; to RT, ST etc. etc. as much as FDE; and to all the corresponding modal and quantificational extensions.

# 7. Note on duality

We have been speaking informally of the 'dual' of a 'system'; I hope this was not too confusing. We give now a series of definitions which allow us to say what we want to say as completely precise statements.

*Local systems.* We have defined such expressions as $\vDash_{FDE}$, $\simeq_{RT}$, etc. in connection with data (e.g. $\alpha \vDash_{FDE} \beta$), dual data ($\delta \vDash_{FDE} \varepsilon$), and datum and dual datum ($\alpha \vDash_{FDE} \delta$, or equivalently $\vDash_{FDE} \alpha \Rightarrow \delta$). So we may define first some 'local' systems – data systems, dual data systems, hybrid systems. Given a previously defined relation $\vDash_S$ (maybe actually written $\simeq_S$ or even in some other way; for our main purposes here $\vDash$ is always either $\vDash$ or $\simeq$, and S is always one of: FDE, EqFDE [and we write then $\simeq_{EqFDE}$, not $_{FDE}$], (R/S/E)(E/C/T/F)), we put:

$S[D] \ =_{df} \ \{\langle \alpha, \beta \rangle : \alpha \vDash_S \beta\}.$

$S[DD] \ =_{df} \ \{\langle \delta, \varepsilon \rangle : \delta \vDash_S \varepsilon\}.$

$S[H] \ =_{df} \ \{\langle \alpha, \delta \rangle : \alpha \vDash_S \delta\}.$

Occasionally we may want to speak also of a 'formula system':

$S[Fla] \ =_{df} \ \{\langle A, B \rangle : A \vDash_S B\}.$

(Of course, since $\vDash_S$ is by assumption defined for say data, it is in effect also defined for formulas: $\langle A, B \rangle \in S[Fla]$ iff $\langle \{\{A\}\}, \{\{B\}\} \rangle \in S[D]$.)

*Global systems.* It is then here natural to define a corresponding 'global system' as the result of putting together the 'local systems':

$S \ =_{df} \ \langle S[D], S[DD], S[H] \rangle.$

*Duality for local systems.* Given local system $S[t]$, we define its dual by:

$d(S[t]) \ =_{df} \ \{\langle X, Y \rangle \text{ of type } t^- : d(Y) \vDash_S d(X)\}$

(where t⁻ is the obvious dualization of the type: D⁻ = DD, (DD)⁻ = D, H⁻ = H [and when Fla is also considered: Fla⁻ = Fla]). We may write the relation ⊨$_S$ among formal objects of type t alternatively as ⊨$_{S[t]}$. Then we have (for all X, Y of type t⁻):

X ⊨$_{d(S[t])}$ Y  iff  d(Y) ⊨$_{S[t]}$ d(X).

*Duality for global systems*.  Then we define the dual d(S) of *global* system S by:

d(S)  =$_{df}$ ⟨d(S[DD]), d(S[D]), d(S[H])⟩.

*Main Proposition*.  – Self-duality (i.e. S = d(S)) is the normal thing: we certainly have d(CPL) = CPL, d(Q) = Q, d(S5) = S5, etc. But the systems studied in the present chapter are mostly not self-dual:

*Proposition*.  Among the systems FDE, EqFDE, (R/S/E)(E/C/T/F), the following are the (only) duality relations:

(1) d(FDE) = FDE. (Hence) d(EqFDE) = EqFDE.

(2) d((R/S/E)(E/C)) = (R/S/E)(C/E). (Hence) d((R/S/E)(T/F)) = (R/S/E)(F/T).

For *formula* systems we have also, in addition to the dualities above, the (only) further dualities:

(3) RT[Fla] = d(RT[Fla]) = (by (2)) RF[Fla] = (by (2) and the first equality here) d(RF[Fla]).

*Proof sketch*.  That there are no *other* dualities or formula-dualities is easily shown by examples, many already explicitly given here before.

(1) This follows from the soundness and completeness results for FDEDC, FDEDDC, FDEHC to be (independently) proved here later. We have e.g.: α ⊨$_{FDE}$ β iff α ⊢$_{FDEDC}$ β iff (by the duality of rules) d(β) ⊢$_{FDEDDC}$ d(α) iff d(β) ⊨$_{FDE}$ d(α).

(2) Let us show for concreteness say that d(EE[DD]) = EC[D]; the form of argument is immediately applicable also to all the other cases here. – By definition we have

$\alpha \vDash_{d(EE[DD])} \beta$  iff  $d(\beta) \vDash_{EE[DD]} d(\alpha)$.

But, writing $\neg\beta$, $\neg\alpha$ for the dual data which result from $\beta$, resp. $\alpha$, by 'driving in' the negation turning slashes and commas and $\vee$'s and $\wedge$'s into semi-colons and backslashes and $\wedge$'s and $\vee$'s until it reaches the propositional variables – i.e. we have what results from $d(\beta)$, resp. $d(\alpha)$, by inserting a negation-sign immediately before each occurrence of propositional variable –, we have:

$d(\beta) \vDash_{EE[DD]} d(\alpha)$  iff  $\neg\beta \vDash_{EE[DD]} \neg\alpha$.

For: the r.h.s. follows from the l.h.s. by general Principle of Substitution (of arbitrary formulas for propositional variables); and the l.h.s. follows from the r.h.s. by again Principle of Substitution, which permits insertion of a *further* negation before each occurrence of variable, and then the fact that we always have $\sigma(\neg\neg p) = \sigma(p)$ (for any model $(S, \sigma)$) and so may delete the double negations in question preserving the truth of the implication statement. – Finally, we have:

$\neg\beta \vDash_{EE[DD]} \neg\alpha$  iff  $\alpha \vDash_{EC[D]} \beta$.

For the l.h.s. means that (w.r.t. every model) every truthmaker of $\neg\beta$ is truthmaker of $\neg\alpha$, or equivalently every falsity-maker of $\beta$ is falsity-maker of $\alpha$, i.e. $\alpha \vDash_{EC} \beta$.

(3) We must show that $RT[Fla] = d(RT[Fla])$. Since we already know (by (2)) that $d(RT[Fla]) = RF[Fla]$, it is enough to show that $RT[Fla] = RF[Fla]$, i.e. that $\forall$ formulas A, B: A $\simeq_{RT}$ B iff A $\simeq_{RF}$ B, or equivalently $\neg$A $\simeq_{RT}$ $\neg$B. By the soundness and completeness results for the calculus AC (with an 'Angellic equivalence' connective $\leftrightarrow$) in Fine 2016, this is equivalent to $\forall$ A, B: $\vdash_{AC}$ A $\leftrightarrow$ B iff $\vdash_{AC}$ $\neg$A $\leftrightarrow$ $\neg$B. The left-to-right part of this is in effect Thm. 2 (pp. 203–204) in Fine 2016; then given this the right-to-left part follows easily (the calculus AC includes axiom A $\leftrightarrow$ $\neg\neg$A, etc.). [I would have preferred of course to argue here only in terms of our own directional deduction calculi; but I don't see how to do it.] □

*Homogeneous duality.*  This is a more recherché notion, but not really entirely unnatural.

The homogeneous dual of (local system) S[t] w.r.t. the global system S = ⟨S[D], S[DD], S[H]⟩ =$_{df}$ d(S[t⁻]).

Since the global system S is usually clear from the context, we may then speak simply of the homogeneous dual of S[t], or more briefly hd(S[t]).

Thus given global system

S = ⟨S[D], S[DD], S[H]⟩,

we have two correct descriptions of its dual:

d(S) = ⟨d(S[DD]), d(S[D]), d(S[H])⟩,

d(S) = ⟨hd(S[D]), hd(S[DD]), hd(S[H])⟩.

# 8. Completeness

We concentrate on the Data Calculi and for the most part leave to the interested reader the straightforward adaptations of the results and proofs to the other calculi.

*Proposition* (Normal-deduction completeness of FDEDC). If $\alpha \vDash_{FDE} \beta$ then $\alpha \vdash_{FDEDC, N} \beta$.

*Proof.* Suppose hypothesis. Start construction of deduction by reducing via the G-rules $\alpha$ to elementary $\alpha'$ and $\beta$ bottom-up to elementary $\beta'$. Since every G-rule is FDE equivalential, we have $\alpha \simeq_{FDE} \alpha'$ and $\beta \simeq_{FDE} \beta'$; and so $\alpha' \vDash_{FDE} \beta'$. So in particular in the canonical model every truthmaker of $\alpha'$ includes truthmaker of $\beta'$; which means in effect that every component of $\alpha'$ is superset of some component of $\beta'$ ($\forall \Gamma \in \alpha' \, \exists \Delta \in \beta'$: $\Delta \subseteq \Gamma$). So we can get from $\alpha'$ to $\beta'$ by first using ↓Del and then ↑Exp. □

*Remark*. (1) Here and in the other cases below the corresponding *soundness* results are easily verified (unless otherwise indicated).

(2) The methods for construction of deduction in the proof above and in the others below are mechanical, giving thus decision procedures for FDE consequence, RT equivalence etc., and for deducibility in the corresponding calculi. ⊣

*Proposition* (Deductive completeness of EqFDEDC). If $\alpha \simeq_{FDE} \beta$ then $\alpha \vdash_{EqFDEDC} \beta$.

*Proof.* We reduce $\alpha$ to elementary $\alpha'$ and $\beta$ bottom-up to elementary $\beta'$, as in the preceding proof. Now we have $\alpha' \simeq_{FDE} \beta'$; so every component of $\alpha'$ is superset of some component of $\beta'$, and vice versa. So with the derived Mutual Inclusion Rule we conclude the deduction. □

*Proposition* (Normal-deduction completeness of FDEDC w.r.t. isotone semantics). If $\alpha \vDash_{FDEi} \beta$ then $\alpha \vdash_{FDEDC, N} \beta$.

*Proof.* Again reduce α to elementary α′ and β to elementary β′. Clearly the G-rules are FDEi equivalential; so α′ ⊨$_{FDEi}$ β′. So in particular in the canonical model every truthmaker of α′ is truthmaker of β′, i.e. every truthmaker of component of α′ is truthmaker of component of β′. But here a truthmaker of component is a union of supersets of singletons of literals in the component, or equivalently a superset of the component. Thus every superset of component of α′ is superset of some component of β′, or equivalently: every component of α′ is superset of some component of β′. So as before we can get from α′ to β′ by ↓Del then ↑Exp, thus concluding the normal deduction. □

*Proposition.* α ⊨$_{FDE}$ β iff α ⊨$_{FDEi}$ β. And (therefore) α ≃$_{FDE}$ β iff α ≃$_{FDEi}$ β.

*Proof.* From the soundness and completeness of FDEDC w.r.t. both semantics. □

*Remark.* It follows of course that α ≃$_{FDEi}$ β iff α ⊢$_{EqFDEDC}$ β. This can also be easily shown directly, along similar lines as in the other arguments above. ⊣

We say that an elementary datum α is *maximal* if it is closed under fusions of components and intermediates between components, i.e. ∀Γ, Δ ∈ α: Γ ∪ Δ ∈ α and ∀Γ, Δ ∈ α ∀Θ with Γ ⊆ Θ ⊆ Δ: Θ ∈ α. (This condition makes sense of course for any datum; but it is the case of elementary data that will be useful here.)

In the basic semantics for FDE, any elementary datum will 'wear its canonical model truthmakers on its sleeves': for α elementary and (S, σ) the canonical model, σ$^+$(α) is literally identical with α – i.e. the components of α are the truthmakers. In the *regular* semantics for RT or FDE, by contrast, it is only *maximal* elementary data that wear their canonical model truthmakers on their sleeves, not elementary data in general. E.g. the canonical model truthmakers of the elementary datum p / p, q, r are {p}, {p, q}, {p, r}, {p, q, r}. But for α *maximal* elementary datum and (S, σ) the canonical model, we have again that σ$^+$(α) = α.

It follows that, if α and β are maximal elementary data and α ≃$_{RT}$ β, then α = β. For RT equivalence means identity of truthmakers w.r.t. all models, so in particular the canonical model; and we then have (for (S, σ) the canonical model) $\alpha = \sigma^+(\alpha) = \sigma^+(\beta) = \beta$.

Similarly, we have: If $\alpha$ and $\beta$ are maximal elementary data and $\alpha \vDash_{\text{FDEr}} \beta$, then every component of $\alpha$ is superset of some component of $\beta$. (*Here* the statement remains true if we drop the 'maximal', in view of the extensional coincidence of $\vDash_{\text{FDEr}}$ with $\vDash_{\text{FDE}}$; but we have not shown this coincidence yet.)

*Proposition* (Deductive completeness of FDEDC w.r.t. regular semantics). If $\alpha \vDash_{\text{FDEr}} \beta$ then $\alpha \vdash_{\text{FDEDC}} \beta$.

*Proof.* As before we reduce $\alpha$ to elementary $\alpha'$ and $\beta$ bottom-up to elementary $\beta'$; and we have $\alpha' \vDash_{\text{FDEr}} \beta'$. But now $\alpha'$ and $\beta'$ need not 'wear their canonical model truthmakers on their sleeves' – so we cannot automatically infer that every component of $\alpha'$ is superset of some component of $\beta'$. But we can use the equivalential derived rules ↕Fus and ↕Int (i.e. use ↓Del and ↑Exp in the corresponding ways) to transform $\alpha'$ into an equivalent *maximal* elementary datum $\alpha*$, and $\beta'$ bottom-up into an equivalent maximal elementary datum $\beta*$. And now given that $\alpha* \vDash_{\text{FDEr}} \beta*$, it does follow that every component of $\alpha*$ is superset of some component of $\beta*$; and so we can use ↓Del and then ↑Exp to get from $\alpha*$ to $\beta*$ and thus complete the deduction. □

*Remark.* The indicated deduction is not in general normal because of the additional use of inclusion and exclusion of fusions and intermediates. ⊣

*Proposition.* $\alpha \vDash_{\text{FDE}} \beta$ iff $\alpha \vDash_{\text{FDEr}} \beta$.

*Proof.* Immediate from the two preceding Propositions and the corresponding soundness results. □

*Proposition* (*Normal*-deduction completeness of FDEDC w.r.t. regular semantics). If $\alpha \vDash_{\text{FDEr}} \beta$ then $\alpha \vdash_{\text{FDEDC, N}} \beta$.

*Proof.* Immediate from the preceding Proposition and the normal-deduction completeness of FDEDC w.r.t. the basic semantics. □

*Remark.*  In our construction above of FDEDC deduction under assumption of α ⊨$_{FDEr}$ β, we observed that from α′ ⊨$_{FDEr}$ β′ we could not 'automatically infer' that every component of α′ is superset of some component of β′. Since we now know that ⊨$_{FDEr}$ coincides with ⊨$_{FDE}$, we can *now* make that inference; and thus we see that use of ↕Fus and ↕Int to reach maximal elementary data was not really necessary and just ↓Del and then ↑Exp was enough, giving a normal deduction.  ⊣

*Proposition* (Completeness of RTDC).  If α ≃$_{RT}$ β then α ⊢$_{RTDC}$ β.

*Proof.*  Reduce by G-rules α to elementary α′ and β bottom-up to elementary β′; then by ↕Fus and ↕Int transform α′ into equivalent maximal elementary datum α* and β′ into equivalent maximal elementary datum β*. Thus α* ≃$_{RT}$ β*, and since these are maximal elementary data, by previous remark we have in fact α* = β* and the deduction is complete.  □

*Proposition* (Soundness and Completeness of STDC).  α ⊢$_{STDC}$ β iff α ≃$_{ST}$ β.

*Proof.*  (⟹): The 'canonical-writing' application of G-rules is basic truthmaker invariant, hence a fortiori semiregular truthmaker invariant. And 'datum rewriting' (i.e. replacing a datum expression by another representing the same datum – i.e. changing the order and/or [positive] number of occurrences of components and/or formulas inside a component) is semiregular truthmaker invariant (indeed, incidentally, even basic truthmaker invariant except for the case of change of number of occurrences of formulas inside a component). (⟸): Reduce α top-down to elementary datum α′, and β bottom-up to elementary datum β′. So by Soundness (and hypothesis that α ≃$_{ST}$ β), α′ ≃$_{ST}$ β′. So in particular in the canonical model α′ and β′ have the same semiregular truthmakers, i.e. fusions of basic truthmakers, i.e. fusions of components. Then using ↕Fus (derived rule of STDC) we complete the deduction.  □

*Proposition* (Lemma for Completeness of ETDC).  For all elementary inferior-replica--free sms data α, β: if α ≃$_{ET}$ β then α = β.

*Proof.*  Suppose that (for α, β as indicated) α ≠ β. We show that α ≄$_{ET}$ β by showing that σ⁺(α) ≠ σ⁺(β) w.r.t. model (S, σ) defined as follows. To each variable p, let us associate a denumerable stock of dedicated new symbols $p_1$, $p_2$, $p_3$, … Then S is the set of all these new symbols, and σ is defined by: σ(p) = ⟨{{$p_1$}, {$p_3$}, {$p_5$}, …}, {{$p_2$}, {$p_4$}, {$p_6$}, …}⟩. By the *standard* truthmaker of a multiset of literals we mean the state fusing the 'initial' states in the sets of truthmakers and falsity-makers of variables as listed above, with distinct states taken for distinct occurrences of literals; thus e.g. the standard truthmaker of [p, p, p, ¬q, ¬q] is {$p_1$, $p_3$, $p_5$, $q_2$, $q_4$}. It should be clear that, if the standard truthmaker of a multiset of literals Γ is also truthmaker of multiset of literals Γ′, then Γ′ must be a superior replica of Γ. – Now since α ≠ β, there must of course be some multiset of literals Γ belonging to one of these sms data but not to the other – say (w.l.o.g., of course) Γ ∈ α, Γ ∉ β. The standard truthmaker of Γ is then truthmaker of α. If it is *not* also truthmaker of β then σ⁺(α) ≠ σ⁺(β) and we are done. If it is truthmaker of β, i.e. of some component of β, then by the above observation such a component must be a superior replica Γ′ of Γ. But now take the standard truthmaker *of Γ′*. This cannot be truthmaker of α (whence σ⁺(α) ≠ σ⁺(β) as desired): for if it were then, again, this would mean that α contains a component Γ″ which is superior replica of Γ′, whence superior replica of Γ (∈ α), contradicting the hypothesis that α is inferior-replica--free.  □

*Proposition* (Soundness and Completeness of ETDC).  ∀ sms data α, β: α ⊢$_{ETDC}$ β iff α ≃$_{ET}$ β.

*Proof.*  (⟹): Clearly the G-rules and ↕IR are truthmaker invariant. (⟸): Reduce α top-down to elementary inferior-replica--free α′, and β bottom-up to elementary inferior-replica--free β′; then α′ ≃$_{ET}$ β′, and so by preceding Proposition α′ = β′ and the deduction is complete.  □

*Proposition* (Soundness and Completeness of RFDDC).  ∀ dual data δ, ε: δ ⊢$_{RFDDC}$ ε iff δ ≃$_{RF}$ ε.

*Proof.*  δ ⊢$_{RFDDC}$ ε iff d(ε) ⊢$_{RTDC}$ d(δ). By Soundness and Completeness of RTDC, d(ε) ⊢$_{RTDC}$ d(δ) iff d(ε) ≃$_{RT}$ d(δ). And clearly d(ε) ≃$_{RT}$ d(δ) iff ε ≃$_{RF}$ δ.  □

*Proposition* (Soundness and Completeness of SFDDC).  ∀ dual data δ, ε: δ ⊢$_{SFDDC}$ ε iff δ ≃$_{SF}$ ε.

*Proof.* Exactly like the preceding proof, now using Soundness and Completeness of STDC. □

*Proposition* (Soundness and Completeness of EFDDC). ∀ sms dual data δ, ε: δ ⊢$_{EFDDC}$ ε iff δ ≃$_{EF}$ ε.

*Proof.* Again as before, now using Soundness and Completeness of ETDC. □

*Question.* Are there nice directional deduction calculi for the systems ST ∩ SF and ET ∩ EF? (There is also the same question for RT ∩ RF; but one might think that it is somewhat over-punctilious since as observed RT and RF are very nearly coincident.)

*Proposition* (Completeness of REDC). If α ⊨$_{RE}$ β then α ⊢$_{REDC}$ β.

*Proof.* Suppose α ⊨$_{RE}$ β. So (equivalently) α / β ≃$_{RT}$ β, whence by completeness of RTDC we have α / β ⊢$_{RTDC}$ β. So we have the REDC deduction:

α        ↑ Exp

α / β    [RTDC deduction]

β. □

*Proposition* (Completeness of SEDC). If α ⊨$_{SE}$ β then α ⊢$_{SEDC}$ β.

*Proof.* Exactly like the preceding proof. □

*Proposition* (Completeness of EEDC). If α ⊨$_{EE}$ β then α ⊢$_{EEDC}$ β.

*Proof.* Ditto. □

*Proposition* (Completeness of RCDDC). If δ ⊨$_{RC}$ ε then δ ⊢$_{RCDDC}$ ε.

*Proof.* Suppose δ ⊨$_{RC}$ ε. So (equivalently) δ ≃$_{RF}$ δ ; ε, whence by completeness of RFDDC we have δ ⊢$_{RFDDC}$ δ ; ε. So we can make the RCDDC deduction:

δ        [RFDDC deduction]

δ ; ε    ↓ Del

ε.  □

*Proposition* (Completeness of SCDDC).  If $\delta \vDash_{SC} \varepsilon$ then $\delta \vdash_{SCDDC} \varepsilon$.

*Proof.*  Exactly like the preceding proof.  □

*Proposition* (Completeness of ECDDC).  If $\delta \vDash_{EC} \varepsilon$ then $\delta \vdash_{ECDDC} \varepsilon$.

*Proof.*  Ditto.  □

# 9. Algebraic counterparts

We give here 'algebraic counterparts' for some of the logics considered in this chapter (in the sense in which e.g. Boolean algebras are an 'algebraic counterpart' of classical propositional logic).

(1) *Distributive lattices and Boolean algebras.*

We begin by stating standard definitions of 'lattice', 'distributive lattice', and 'Boolean algebra'. – Operations denoted by $\wedge$ and $\vee$ below are always assumed to be binary, and by $\neg$ unary.

Lattice $=_{df}$ Algebra $(U, \wedge, \vee)$ satisfying the commutative, associative and idempotent laws for $\wedge$ and $\vee$, and the absorption laws

$$p = p \wedge (p \vee q)$$

$$p = p \vee (p \wedge q).$$

Distributive lattice $=_{df}$ Lattice $(U, \wedge, \vee)$ satisfying the distributive laws

$$p \wedge (q \vee r) = (p \wedge q) \vee (p \wedge r)$$

$$p \vee (q \wedge r) = (p \vee q) \wedge (p \vee r).$$

Boolean algebra $=_{df}$ Algebra $(U, \wedge, \vee, \neg)$ where $(U, \wedge, \vee)$ is a distributive lattice, satisfying the double negation and De Morgan laws

$$\neg\neg p = p$$

$$\neg(p \wedge q) = \neg p \vee \neg q$$

$$\neg(p \vee q) = \neg p \wedge \neg q,$$

and the equations

$(p \wedge \neg p) \vee q \; = \; q$

$(p \vee \neg p) \wedge q \; = \; q.$

*Remark*. The term 'Boolean algebra' is also often used in other senses – in particular, including 'constants' 0 and 1. The above is one standard sense (see e.g. McKenkie et al. 1987, p. 17) and what will be most convenient for our purposes here. ⊣

(2) *FDE algebras*.

FDE algebra $=_{df}$ Algebra $(U, \wedge, \vee, \neg)$ where $(U, \wedge, \vee)$ is distributive lattice and the double negation and De Morgan laws are satisfied.

So this is exactly like the definition of Boolean algebra above except for the requirement that the last two ('counter-relevantistic') equations be satisfied. (FDE algebras in this sense are usually called 'De Morgan lattices'. See e.g. Anderson & Belnap 1975 § 18.3 [written by Dunn], or Font 1997, beginning of § 2. We avoid this terminology here partly because we will reserve 'lattice' for use in connection with 'lattice-like' and not 'Boolean-algebra--like' algebras – see below.)

(3) *RT lattices and RT algebras*.

RT lattice $=_{df}$ Algebra $(U, \wedge, \vee)$ satisfying the commutative, associative, idempotent and distributive laws, as well as:

$p \vee (p \wedge q \wedge r) \; = \; p \vee (p \wedge q) \vee (p \wedge q \wedge r)$

$p \wedge (p \vee q \vee r) \; = \; p \wedge (p \vee q) \wedge (p \vee q \vee r).$

*Question*. Do these conditions on intermediates already follow from the commutative, associative, idempotent and distributive laws?

*Remark*. Note that the equations

$p \vee q \; = \; p \vee q \vee (p \wedge q)$

$p \wedge q \; = \; p \wedge q \wedge (p \vee q)$

hold universally in every RT lattice: we can transform the r.h.s. into the l.h.s. by applying distribution then (three times) idempotence (also commuting and associating as needed). ⊣

RT algebra $=_{df}$ Algebra $(U, \wedge, \vee, \neg)$ where $(U, \wedge, \vee)$ is RT lattice and the double negation and De Morgan laws are satisfied.

(4) *ST / SF / ST∩SF lattices and algebras*.

ST ∩ SF lattice $=_{df}$ Algebra $(U, \wedge, \vee)$ where $\wedge$ and $\vee$ are commutative, associative and idempotent, and obey the two restricted absorption laws

$$p \vee q \ = \ p \vee q \vee (p \wedge q)$$

$$p \wedge q \ = \ p \wedge q \wedge (p \vee q),$$

and *either* obey distributivity of $\wedge$ over $\vee$ *or* obey distributivity of $\vee$ over $\wedge$.

*Questions*. Can this disjunctive distributivity condition be replaced by some equation or equations condition? Or at least can the whole definition be reformulated somehow in terms of equations only? (I.e. is the class of ST ∩ SF lattices an 'equational class'.)

ST ∩ SF algebra $=_{df}$ Algebra $(U, \wedge, \vee, \neg)$ where $(U, \wedge, \vee)$ is an ST ∩ SF lattice and the double negation and De Morgan laws are satisfied.

ST lattice $=_{df}$ ST ∩ SF lattice $(U, \wedge, \vee)$ satisfying moreover

$$p \wedge (q \vee r) \ = \ (p \wedge q) \vee (p \wedge r).$$

ST algebra $=_{df}$ ST ∩ SF algebra $(U, \wedge, \vee, \neg)$ where $(U, \wedge, \vee)$ is an ST lattice.

And similarly for SF lattice and SF algebra, now with the distributive law for $\vee$ over $\wedge$ only.

*Remarks.* (1) These four classes (ST/SF lattices/algebras) are 'equational classes', since of course once one distributive law is imposed the disjunctive condition becomes redundant and so may be omitted.

(2) Note that clearly the class of ST $\cap$ SF lattices (resp., algebras) is the union of the class of ST lattices (algebras) with the class of SF lattices (algebras). $\dashv$

(5) *ET / EF / ET$\cap$EF lattices and algebras.*

By the $\wedge$-*collapse law* we will mean the equation (given in Krämer forthcoming; I haven't seen it elsewhere)

$(p \wedge p) \vee p \ = \ p \wedge p.$

And by the $\vee$-*collapse law* we mean

$(p \vee p) \wedge p \ = \ p \vee p.$

ET $\cap$ EF lattice $=_{df}$ Algebra (U, $\wedge$, $\vee$) where $\wedge$ and $\vee$ are commutative and associative, and obey *either* $\vee$-idempotence and $\wedge$-collapse and distributivity of $\wedge$ over $\vee$ *or* $\wedge$-idempotence and $\vee$-collapse and distributivity of $\vee$ over $\wedge$.

*Questions.* Again, can this disjunctive condition be replaced by equations? Is the class 'equational'?

ET $\cap$ EF algebra $=_{df}$ Algebra (U, $\wedge$, $\vee$, $\neg$) where (U, $\wedge$, $\vee$) is ET $\cap$ EF lattice and the double negation and De Morgan laws are satisfied.

ET lattice $=_{df}$ ET $\cap$ EF lattice (U, $\wedge$, $\vee$) satisfying moreover

$p \vee p \ = \ p$

$(p \wedge p) \vee p \ = \ p \wedge p$

$p \wedge (q \vee r) \ = \ (p \wedge q) \vee (p \wedge r).$

ET algebra $=_{df}$ ET $\cap$ EF algebra (U, $\wedge$, $\vee$, $\neg$) where (U, $\wedge$, $\vee$) is ET lattice.

And similarly (dually) for EF lattice and EF algebra.

*Remarks*. (1) Again, these four classes (ET/EF lattices/algebras) are of course 'equational'.

(2) Again, clearly the class of ET ∩ EF lattices (resp., algebras) is the union of the class of ET lattices (algebras) with the class of EF lattices (algebras). ⊣

Where K is a class of algebras, by A $=_K$ B we mean that the equation A = B holds universally (i.e. for all values of the variables) in every algebra of the class K. We use obvious names for the classes of algebras defined above: DL for (the class of all) distributive lattices, BA for Boolean algebras, FDEA for FDE algebras, RTL for RT lattices, RTA for RT algebras, and so on.

*Proposition*. In each case below the three statements are equivalent (for all formulas A, B):

(1) A $\vdash_{EqDC}$ B; A $=_{BA}$ B; A $\simeq_{CPL}$ B.

(2) A $\vdash_{EqFDEDC}$ B; A $=_{FDEA}$ B; A $\simeq_{FDE}$ B.

(3) A $\vdash_{RTDC}$ B; A $=_{RTA}$ B; A $\simeq_{RT}$ B.

We have moreover:

(4) If A $\vdash_{STDC}$ B then A $=_{STA}$ B.

(5) If A $\vdash_{ETDC}$ B then A $=_{ETA}$ B.

(And similarly for SF and EF.)

*Proof sketch*. In each of the cases (1)–(3) we already know that the third statement implies the first (completeness theorem); so to complete a 'circle' of implications it will be enough to 'squeeze' the second statement between the first and the third, i.e. to show that the first statement implies the second and that the second implies the third. And with (4) and (5) we just need the 'first half' of this. – We can proceed in all these cases by the same general form of argument: –

*First statement implies second statement*, i.e. soundness of the given calculus w.r.t. the given class of algebras. – We show that each rule of the calculus corresponds (translating commas and slashes by resp. ∧ and ∨) to a schematic identity all instances of which are universally true in every algebra of the given class. We leave here detailed verifications to the reader, but give the following list of sufficient conditions for soundness of rules w.r.t. classes of algebras (leaving commutativity and associativity sometimes tacit):

Datum rewriting (i.e. replacement of a datum expression by another expression for the same datum): commutativity, associativity, idempotence.

Sms datum rewriting: commutativity, associativity, ∨-idempotence.

↕∧: automatic (same corresponding formula).

↕¬¬: double negation equation.

↕∨: distribution of ∧ over ∨.

↕¬∨: De Morgan for ¬∨.

↕¬∧: De Morgan for ¬∧, and distribution of ∧ over ∨.

↕EM: EM-equation, and distribution of ∧ over ∨.

↕NC: NC-equation, distribution of ∨ over ∧, absorption. (With this we get generalized NC-equation (p ∧ ¬p ∧ r) ∨ q = q as follows: (p ∧ ¬p ∧ r) ∨ q = (distr.) ((p ∧ ¬p) ∨ q) ∧ (r ∨ q) = (NC-equation) q ∧ (r ∨ q) = (absorption) q.)

↕SS: absorption.

↕Int: intermediates equation.

↕IR: ∧-collapse, and distribution of ∧ over ∨. (With this we get generalized ∧-collapse equation (p ∧ p ∧ q) ∨ (p ∧ q) = p ∧ p ∧ q as follows: (p ∧ p ∧ q) ∨ (p ∧ q) = (distr.) = q ∧ ((p ∧ p) ∨ p) = (∧-collapse) p ∧ p ∧ q. And this in turn implies the even more general schema S ∨ I = S where S, I stand for conjunctions of variables, with the same variables occurring in both, and each such variable having number of occurrences in S $\geq$ its number of occurrences in I.)

*Second statement implies third statement*. We show in each case that our models yield a special case of the appropriate algebras. In the case of CPL this is simply the universe {T, F} with the ordinary truth-functions $\wedge$, $\vee$, $\neg$. In the other cases the universe will consist of a set of bilateral propositions of appropriate kind (isotone, regular) closed under appropriately defined operations of conjunction, disjunction and negation. So in each case we will have that $\sigma(A \wedge B)$ is the conjunction (in the appropriate sense) of $\sigma(A)$ with $\sigma(B)$, and so on. Thus the universal truth of the equality in the class of algebras will imply the equivalence, since the values which the interpretation $\sigma$ in a model gives to the formulas, together with the appropriate operations, is an algebra of the given class. (Note here the utility of the *isotone* semantics for FDE [which is what we use for the present proof]. Here, unlike in the basic or the regular semantics, we always have literally e.g. $\sigma(p \wedge (p \vee q)) = \sigma(p)$.) □

*Questions and Remarks.* Do we have also

(4′) If A $=_{STA}$ B then A $\simeq_{ST}$ B.

(5′) If A $=_{ETA}$ B then A $\simeq_{ET}$ B.

This would mean of course that (4) and (5) in the Proposition above can be replaced by three-term equivalences just like (1)–(3).

Note that (4′) and (5′) imply, by duality:

(6) A $\vdash_{SFDDC}$ B iff A $=_{SFA}$ B iff A $\simeq_{SF}$ B.

(7) A $\vdash_{EFDDC}$ B iff A $=_{EFA}$ B iff A $\simeq_{EF}$ B.

As already pointed out above, we have

(ST $\cap$ SF)A $=$ STA $\cup$ SFA,

(ET $\cap$ EF)A $=$ ETA $\cup$ EFA.

So (if (4′) and (5′) hold) we would get here also:

(8) A $=_{(ST \cap SF)A}$ B iff A $\simeq_{ST \cap SF}$ B.

(9) A $=_{(ET \cap EF)A}$ B iff A $\simeq_{ET \cap EF}$ B.

– The argument in the proof above for showing, in the cases (1)–(3), that the second statement implies the third – this argument does not directly apply to show (4′) and (5′), because the bilateral propositions of appropriate kind (semiregular, arbitrary) constructed from the truthmaker semantics models do *not* (together with the relevant operations) form an algebra of the appropriate kind. E.g. we do not always have σ(p ∧ (q ∨ r)) = σ((p ∧ q) ∨ (p ∧ r)). What we do have is σ⁺(p ∧ (q ∨ r)) = σ⁺((p ∧ q) ∨ (p ∧ r)). One then naturally thinks of constructing instead an algebra of *uni*lateral propositions, with just the truthmaker set. But although this works fine with the operations of conjunction and disjunction (which 'depend only' on the truthmaker sets of the components [for determination of the truthmaker set of the compound]) – which indeed is why the present difficulty does not recur in the 'lattice-like' context considered below –, yet negation cannot straightforwardly be dealt with, since the truthmaker set of the negation of a proposition is not in general determined by the truthmaker set of the negatum (but by its *falsity*-maker set).

Nor can we use the procedure to show directly the left-to-right parts of (8) and (9). For the algebras of semiregular/arbitrary bilateral propositions constructed from the truthmaker models need not satisfy the disjunctive conditions in the definitions of ST∩SF/ET∩EF algebras. ⊣

*Proposition.* For any negation-free formulas A, B: A ⊨$_{FDE}$ B iff A ⊨$_{CPL}$ B; and therefore also A $\simeq_{FDE}$ B iff A $\simeq_{CPL}$ B.

*Proof.* Suppose (for the non-trivial direction) A ⊨$_{CPL}$ B. Then α ⊨$_{CPL}$ δ where α is the elementary datum resulting from applications of the ↓ G-rules starting from A, and δ the elementary dual datum resulting from 'reverse-applying' the ↑ G-rules with B at the end. Since A, B are negation-free, so are α, δ. So α criss-crosses δ and thus α ⊢$_{FDEHC}$ δ, whence A ⊢$_{FDEHC}$ B, whence A ⊨$_{FDE}$ B. □

*Proposition.* In each case below the three (or two, with (8)–(9)) statements are equivalent, for all negation-free formulas A, B:

(1) $A \vdash_{EqDC} B$; $A =_{DL} B$; $A \simeq_{CPL} B$.

(2) $A \vdash_{EqFDEDC} B$; $A =_{DL} B$; $A \simeq_{FDE} B$.

(3) $A \vdash_{RTDC} B$; $A =_{RTL} B$; $A \simeq_{RT} B$.

(4) $A \vdash_{STDC} B$; $A =_{STL} B$; $A \simeq_{ST} B$.

(5) $A \vdash_{ETDC} B$; $A =_{ETL} B$; $A \simeq_{ET} B$.

(6) $A \vdash_{SFDDC} B$; $A =_{SFL} B$; $A \simeq_{SF} B$.

(7) $A \vdash_{EFDDC} B$; $A =_{EFL} B$; $A \simeq_{EF} B$.

(8) $A =_{(ST \cap SF)L} B$; $A \simeq_{ST \cap SF} B$.

(9) $A =_{(ET \cap EF)L} B$; $A \simeq_{ET \cap EF} B$.

*Proof sketch.* (6)–(9) follow easily from (4) and (5). – For (1)–(5), we follow the same general form of argument as in the earlier Proposition with negation present. That in each case the second statement implies the third is shown as before (now using unilateral propositions for cases (4) and (5)). That in each case the first statement implies the second can also be shown essentially as before; there is only the complication that negation might appear in a deduction of negation-free formula from negation-free formula, but this we deal with as follows. In all cases (2)–(5), the only rules in the calculus involving negation are the G-rules $\updownarrow\neg\neg$, $\updownarrow\neg\wedge$, $\updownarrow\neg\vee$. But it is clear from our completeness proofs that none of these will have to be used in a deduction of a negation-free formula from a negation-free formula: our deductions (described in those completeness proofs) always proceed by G-rules top-down and bottom-up (which will be negation-free if the given formulas are negation-free), and then an essentially 'structural' middle portion. (We might use $\updownarrow$Fus as derived rule, but its derivation we saw uses only $\updownarrow\wedge$ and $\updownarrow\vee$.) So in each case we have (for negation-free A, B): $A \vdash B$ implies $A \simeq B$, which implies that in fact B is deducible from A in the restricted calculus without rules involving negation, which implies (by same reasoning as before) that A = B holds universally in every algebra of the appropriate class. – In case (1) this argument is not directly applicable, but we can argue as follows: $A \vdash_{EqDC} B$ implies $A \simeq_{CPL} B$, which implies (by preceding Proposition) $A \simeq_{FDE} B$, which implies (by (2)) $A =_{DL} B$.  □

## 10. A compactness theorem for FDE consequence

In our semantics for FDE the interpretation function $\sigma$ was extended to finite 'conjunctive sets' as for conjunctions, and to finite 'disjunctive sets' as for disjunctions. Now, it can also be extended to arbitrary – including infinite – sets in the similar way. Thus in particular for set of formulas $\Gamma$ understood conjunctively, $\sigma^+(\Gamma)$ will consist of the states fusing truthmakers of the elements of $\Gamma$; and for set of formulas $\Delta$ understood disjunctively, $\sigma^+(\Delta)$ is the set $\{\sigma^+(A) : A \in \Delta\}$ of truthmakers of elements of $\Delta$. To avoid ambiguity we may write resp. $\sigma^+(\Gamma^\wedge)$ and $\sigma^+(\Delta^\vee)$.

We can then say that $\Gamma \vDash_{FDE} \Delta$ if w.r.t. every model, every truthmaker of $\Gamma$ (understood conjunctively) includes truthmaker of $\Delta$ (understood disjunctively) – i.e. $\forall$ model $(S, \sigma)$: $\forall\Sigma \in \sigma^+(\Gamma^\wedge)$ $\exists\Sigma' \in \sigma^+(\Delta^\vee)$: $\Sigma' \subseteq \Sigma$.

We show here (principally for later use in our completeness proof for quantificational FDE) the following compactness theorem: If $\Gamma \vDash_{FDE} \Delta$ then $\exists\Gamma_0 \subseteq_{fin} \Gamma$ $\exists\Delta_0 \subseteq_{fin} \Delta$: $\Gamma_0 \vDash_{FDE} \Delta_0$. – Before proving this we define the concept of the tree for $\Gamma$.

Let $\Gamma$ be an arbitrary set of formulas; we may represent its elements in their lexicographic ordering as $A_1$, $A_2$, $A_3$, … The *tree for $\Gamma$* is the (dyadic) tree constructed as follows. At the origin of the tree we have { }, which is the only component of the datum $\top = \{\{\ \}\}$ representing the conjunction of the empty initial segment of the sequence $\langle A_1, A_2, A_3, \ldots\rangle$. – This is stage 0. Next, as single immediate successor node to { }, we have $\{A_1\}$, which is the only component of the datum $\{\{A_1\}\}$ representing the conjunction of the initial segment $\langle A_1\rangle$ of the sequence $\langle A_1, A_2, A_3, \ldots\rangle$. We then proceed to apply the $\downarrow$ G-rules in some conventionally established order until an elementary datum is obtained. The nodes of the tree correspond to components of the appropriate datum. A successor node either repeats its predecessor if it was 'inert' in the application of the G-rule, or will come from it by an elimination of double negation etc. The only cases of dyadic branching are those corresponding to applications of $\downarrow\vee$ or $\downarrow\neg\wedge$. – This is stage 1. Next, to each bottom node we adjoin as immediate successor the result of adding $A_2$ to the given set of literals. This represents the conjunction of $\langle A_1, A_2\rangle$. Again we apply $\downarrow$ G-rules until elementary datum is obtained. – This is stage 2. – And so on for all the formulas in $\Gamma$. (This is of course an infinite tree if $\Gamma$ is infinite, and finite tree if $\Gamma$ is finite.)

*Remarks*. (1) Clearly a set of literals is truthmaker of $A_1 \wedge \ldots \wedge A_n$ in the canonical model iff it is the set of literals in some branch of the n-fragment of the tree for $\Gamma$ (i.e. the subtree where we stop at the end of stage n; or equivalently the tree for $\{A_1, \ldots, A_n\}$).

(2) Likewise a set of literals is truthmaker of $\Gamma$ (i.e. fusion of resp. truthmakers of $A_1, A_2, \ldots$) in the canonical model iff it is the set of literals in some branch of the tree for $\Gamma$. ⊣

*Proposition*. If $\Gamma \vDash_{FDE} \Delta$ then $\exists \Gamma_0 \subseteq_{fin} \Gamma \; \exists \Delta_0 \subseteq_{fin} \Delta$: $\Gamma_0 \vDash_{FDE} \Delta_0$.

*Proof*. Suppose hypothesis. Consider the tree for $\Gamma$. Since in the canonical model every truthmaker of $\Gamma$ includes truthmaker of some element of $\Delta$, we have that: $\forall$ branch b of the tree for $\Gamma$, there exists some formula $B \in \Delta$ s.t. the set of literals in b includes truthmaker for B. But such truthmaker for a particular formula B is of course *finite* set of literals. Now by the *$\Delta$-subtree* of the tree for $\Gamma$ we will mean the subtree of the tree for $\Gamma$ obtained by reducing each branch b to its shortest initial segment $b_0$ satisfying the condition that the set of literals present in $b_0$ already includes truthmaker for some formula or other in $\Delta$. (Note that if $b_1$ and $b_2$ in the original tree share the initial segment $b_0$ and $b_1$ is reduced to $b_0$ then so is $b_2$.) Thus the $\Delta$-subtree is a finitely generated (specifically, dyadic) tree where every branch is finite, and so by König's lemma it is a finite tree. Let $\Gamma_0$ be the set of formulas of $\Gamma$ appearing in the $\Delta$-subtree. (We may if we wish expand the new tree to the full tree for $\Gamma_0$.) So now each branch of our new (finite) tree, where the set of literals is a truthmaker for $\Gamma_0$, includes truthmaker for some element of $\Delta$. Let then $\Delta_0$ be subset of $\Delta$ with a formula thus corresponding to each branch of the $\Gamma_0$-tree. So $\Gamma_0 \vDash_{FDE} \Delta_0$, as desired. □

## 11.  Variations in truthmaker semantics

Truthmaker semantics is a general method, subject to many variations in its concrete applications to various logical systems. In what follows I list, and make some comments on, what seem to me to be the main such variation-points.

(1) *Notion of 'frame' or 'space'.*

(1.1) *How deep do we go in the analysis of states?*  Intuitively, I at least think of a 'state' as a conjunction of asserted or negated atomic situations – or in other words, using 'elementary situation' for asserted or negated atomic situation, as a conjunction of elementary situations. So the domain of *states* is determined by the domain of *elementary situations*, which in its turn is determined by the domain of *atomic situations*, which in *its* turn is no doubt determined by the domains of 'individuals' and 'attributes of individuals':

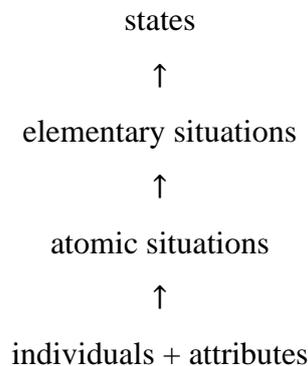

states

↑

elementary situations

↑

atomic situations

↑

individuals + attributes

So in the formulation of an abstract formal semantics of the present (truthmaker semantics) kind, there are four options concerning 'how deep we go in the analysis of states'.

At one extreme there is the option favoured by Fine where the *states* themselves are regarded as unanalyzed 'points', endowed with a part-whole relation. (We obviously have here a kind of oxymoron, or at least mixed metaphor – but that is not my fault!) Thus the basic notion of 'frame' or 'space' here is what Fine calls a *state space*, also called in other contexts a complete lattice, i.e. a partially ordered set (S, ≤) where every subset of the universe S has a least upper bound (whence also every subset of S has a greatest lower bound).

Next there is the option we have favoured here, where it is rather *elementary situations* that are regarded for present purposes as unanalyzed points. Thus now the basic notion of frame or space is what we have called an *elementary situation space* (e.s.s.), i.e. simply an arbitrary set S. From this 'states' can now be defined (as we have done) as sets of elementary situations, and the previous part-whole relation as simply the subset relation on such sets; and the earlier fusion or least upper bound is now simply union. (We think of course of the set construction here as representing conjunction, so that part-whole is then subconjunction and fusion the conjunction of the conjuncts of the given conjunctions.)

These two approaches can be compared to more familiar cases such as say the following two types of semantics for systems of modal propositional logic: (i) algebraic semantics in terms of a suitable class of algebras, where the propositions are unanalyzed points and the operations of 'meet', 'join' etc. are given inside the frame (algebra); versus (ii) the possible world semantics where a frame includes a set W thought of as the set of possible worlds, so that it is now possible worlds that are regarded as 'points' and propositions can be defined as sets of worlds, as well as the corresponding operations on propositions defined as certain operations on such sets (intersection, union, etc.). (i) is admittedly more general, and this may be useful in certain contexts; but on the other hand (ii) is somewhat simpler, and has a more solid 'feel', and in it our desired structures are naturally defined into existence rather than (as in the other case) brutally abstractly postulated – and so (ii) may well be preferred in contexts where the loss in generality is not important (and the loss may even be desirable; for surely generality like everything else can be excessive). And I would apply exactly similar remarks to the state space and elementary situation space versions of truthmaker semantics.

– Next, we might go even deeper in the analysis of states and take *atomic situations* as our unanalyzed points. Thus the basic notion of frame or space would now be what we might call an *atomic situation space* (a.s.s.), i.e. again simply an arbitrary set U, but now thought of as the set of atomic situations. We can then define an elementary situation as a pair $\langle u, 0 \rangle$ or $\langle u, 1 \rangle$ for $u \in U$ (thought of intuitively as resp. the negation of u and the assertion of u i.e. in effect u itself); and then a state as a set of elementary situations, and so on as before. – This third option can be useful sometimes (see e.g. the remarks below on FDE⁺). – One interesting point worth noting here is that we can now (we couldn't before) define a notion of *explicitly contradictory state* as a state containing both $\langle u, 0 \rangle$ and $\langle u, 1 \rangle$ for some $u \in U$.

Finally there is the fourth option, of minimum generality and maximum depth, where we go all the way down to *individuals* and *attributes of individuals* for our 'points'. – I haven't thought much about this option, but it does seem a little excessive. (If generality can be excessive, depth too.)

(1.2) *Is there 'designation' of factuality?* The 'designated world' is a familiar creature from possible world semantics: a frame (or model) with designated world is a structure $(W, w_0, \ldots)$ with $w_0 \in W$ and so on, where intuitively $W$ is thought of as the set of all possible worlds and $w_0$ as the *actual* world. Likewise, in truthmaker semantics a frame might (or might not) also include some states or situations 'designated' as *factual* (as opposed to counterfactual) states or situations.

In the Finean context of state spaces, a state space with designated ('factual') states (or corresponding model) should be a structure $(S, S_0, \leq, \ldots)$ with $S_0 \subseteq S$ and so on, with certain natural conditions imposed on the set $S_0$ of 'factual states', such as that parts of factual states are factual states, and fusions of factual states are factual states. See the definition of D-model in Fine 2014, p. 571; also Fine 2016, p. 224.

In the context of atomic situation spaces, we should define an *atomic situation space with designated facts* (a.s.s.d.f.) as an ordered pair $(U, U_0)$ where $U$ is a.s.s. (i.e. any set) and $U_0 \subseteq U$. We then define factuality of elementary situations and states in the obvious way: $\langle u, 1 \rangle$ is factual iff $u$ is (i.e. iff $u \in U_0$); $\langle u, 0 \rangle$ is factual iff $u$ isn't; state $\Sigma$ is factual iff every elementary situation in $\Sigma$ is factual. It then follows from such definitions that, w.r.t. any a.s.s.d.f.: parts of factual states are factual; fusions of factual states are factual; there are factual states (e.g. the empty state { }); no factual state is explicitly contradictory; there exists a fusion of all factual states; etc.

Finally, in the context of elementary situation spaces, I think we should define an *elementary situation space with designated facts* (e.s.s.d.f.) as an ordered pair $(S, S_0)$ where $S$ is e.s.s. (arbitrary set) and $S_0 \subseteq S$ with $card(S_0) = card(S - S_0)$. This cardinality condition may seem a little artificial; but since we are thinking of elementary situations as asserted or negated atomic situations, of course the factual elementary situations must be 'half' of the elementary situations; we couldn't have e.g. $S_0 = \varnothing$, or $S_0 = S$ (unless S itself $= \varnothing$). (With a.s.s.d.f.'s by contrast there should be no constraints: $U_0 = \varnothing$, or $U_0 = U$, or anything else in between – all such cases should be permitted.) Again, we define factual state as state all elements of which are factual elementary situations; and so on.

(Note that not every e.s.s. S can be extended to e.s.s.d.f.: for card(S) might be an odd natural number. Indeed one might prefer to define an e.s.s. as a set partitionable into two equinumerous subsets [or equivalently, whose cardinality is either infinite or finite and even]; and even include a corresponding cardinality condition on the universe of states in the definition of state space. But 'in practice' these extra conditions would hardly ever be used.)

(1.3) *Is there designation of possibility?* Similar considerations apply if we want to endow our 'frames' or 'spaces' with states or combinations of situations designated as *possible*.

Thus Fine defines a *modal state space* as a structure $(S, S^\diamond, \leq)$ where $(S, \leq)$ is state space and $S^\diamond$ ('possible states') is 'downward closed' subset of S (i.e. whenever $\Sigma \in S^\diamond$ and $\Sigma' \leq \Sigma$, also $\Sigma' \in S^\diamond$ [parts of possible states are possible states]).

In the context of a.s.s.'s, we should define an *atomic situation space with designated possibilities* (a.s.s.d.p.) as a pair $(U, U^\diamond)$ where U is a.s.s. and $U^\diamond$ is non-empty collection of sets choosing one from each pair of elementary situations $\langle u, 1 \rangle$, $\langle u, 0 \rangle$ ($u \in U$) (intuitively: these are the possible truth-value combinations for the atomic situations; there must be at least one, since at least the *factual* truth-value combination is a possible one). Then we say that an *elementary situation* is possible if it belongs to some set in $U^\diamond$; and that a *state* is possible if it is subset of some set in $U^\diamond$. It then follows that parts of possible states are possible. And also that no explicitly contradictory state is possible (although of course a state might be impossible without being explicitly contradictory).

And finally, in the context of e.s.s.'s, we define an *elementary situation space with designated possibilities* (e.s.s.d.p.) as a pair $(S, S^\diamond)$ where S is e.s.s. and, for some partition $\{S_0, S_1\}$ (with $S_0 \neq S_1$) of S with card($S_0$) = card($S_1$) and some bijection f from $S_0$ to $S_1$, $S^\diamond$ is non-empty collection of sets choosing one from each pair of elementary situations $s$ ($\in S_0$), $f(s)$. Possibility for elementary situations and states is then defined as before in the context of a.s.s.'s.

(2) *Notion of bilateral proposition*.  Once given a notion of frame or space, the next question is what are to be the possible values of the (unextended) interpretation function $\sigma$, i.e. what will be the relevant notion of 'bilateral proposition'.

(2.1) *Constraints on the truthmaker and falsity-maker sets individually*.  We might impose *no* constraints (as in the semantics for ET and the basic semantics for FDE); or upward closure, i.e. closure under superstates ('isotone' bilateral propositions, as in our isotone semantics for FDE); or closure under fusions and intermediates ('regular' bilateral propositions, as in the semantics for RT); or only closure under fusions ('semiregular' bilateral propositions, as in the semantics for ST), or only closure under intermediates; and perhaps other constraints might be considered. (Whatever such constraint is imposed, it should be imposed equally on the truthmaker set and the falsity-maker set – for of course the truthmakers of one proposition are the falsity-makers of another [e.g. its negation] and vice versa.) This is of course related to the question of the notions of truth*making* and falsity-*making* which our semantics is going to use, i.e. the question of how the interpretation-function $\sigma$ will be *extended* – and a proper truthmaker semantics must always be 'coherent' in this respect.

In particular, if we are going to use the 'basic' notions of truthmaking and falsity-making, then there should indeed be *no* constraints on the individual sets of truthmakers and falsity-makers. For let $\{\Sigma, \Sigma', \Sigma'', \ldots\}$ be *any* set of states. Then, intuitively, there is a proposition of which this is precisely the truthmaker set (in the basic sense of truthmaking), e.g. the disjunction $\vee(\Sigma, \Sigma', \Sigma'', \ldots)$; and there is a proposition of which this is precisely the falsity-maker set, e.g. $\neg\vee(\Sigma, \Sigma', \Sigma'', \ldots)$. (The situation here is similar to the situation with intensional propositions: any set of worlds $\{w, v, u, \ldots\}$ should count as a bona fide intensional proposition, since it is the set of worlds where the disjunction $\vee(w, v, u, \ldots)$ is true.)

(2.2) *Constraints on coordination*.  There is next the question of what *pairs* $\langle X, Y \rangle$ of sets of states should count as bilateral propositions. I have used here only the very minimal requirement that X ≠ Y; but although this 'works' in the sense of allowing construction of suitable semantics for this or that system, this should not make us think that it is anywhere near a proper 'intuitive' requirement on coordination. Try to find a proposition that is made true by the state of snow being white and no other state, and made false by the state of grass being green and no other state – and you will try for ever.

Nor is a logical system where this issue is felt far to seek. Take e.g. the relevantistic system which we may call FDE$^+$(a), defined as follows. We enrich $\mathcal{L}$CPL (with primitive connectives ¬, ∧, ∨) by addition of relevantistic entailment connective ⇒. A *basic entailment* is a formula of the form A ⇒ B where A, B ∈ $\mathcal{L}$CPL. Then a formula of $\mathcal{L}$FDE$^+$(a) is defined by: (i) basic entailments are formulas; (ii) if A and B are formulas then so are ¬A, A ∧ B, A ∨ B. – We would then like to use something like our truthmaker semantics for FDE to deal with this language. For w.r.t. a model we can take A ⇒ B as true or false according as it holds or not w.r.t. the model that every truthmaker of A includes truthmaker of B; and from that the truth-values of truth-functional compounds can be calculated in the usual way; so we have the notion of the truth-value of a formula of $\mathcal{L}$FDE$^+$(a) w.r.t. a model and so can define validity of such a formula as its being true w.r.t. all models.

However, take e.g. the formula

$$\text{p} \Rightarrow \text{q} \ \textbf{.}{\rightarrow}\textbf{.} \ \neg\text{q} \Rightarrow \neg\text{p}.$$

We would like this to count as valid. But it is not so by the definition above: we might have e.g. σ(p) = ⟨{Σ₁}, {Σ₂}⟩, σ(q) = ⟨{Σ₁}, {Σ₃}⟩ (with Σ₂ ⊈ Σ₃).

A radical solution is the following (say in the context of a.s.s.'s). We define the set of bilateral propositions as the least set including all pairs ⟨{{⟨u, 1⟩}}, {{⟨u, 0⟩}}⟩ for u ∈ (the a.s.s.) U and closed under the operations (defined in the obvious way) of negation and arbitrary conjunction. Then the interpretation function σ is required to assign bilateral propositions in *this* sense to the propositional variables. – It would be interesting to have more 'moderate' ways of constraining coordination so as to permit treatment of e.g. this FDE$^+$(a).

(The reader may have wondered about the '(a)' in the name FDE$^+$(a). It is that there is another, slightly more comprehensive system of 'first-degree' entailments [with entailment connective but only between $\mathcal{L}$CPL formulas], which may be called FDE$^+$(b). Clause (i) in the definition above of $\mathcal{L}$FDE$^+$(a) formula is now liberalized to (i′) basic entailments *and propositional variables* are formulas [and clause (ii) is as before]. Truthmaker semantics for FDE$^+$(a) can then be extended to truthmaker semantics for FDE$^+$(b) by including 'designation of factuality' – so that propositional variables outside the scope of entailment connective in a given formula will get a truth-value [the truth-value of the associated bilateral proposition], and so the whole formula gets a truth-value. – Exactly similar syntactical and semantical definitions and remarks apply to other systems of 'first-degree' entailment/containment/equivalence such as EqFDE$^+$(a)/(b), RT$^+$(a)/(b), RE$^+$(a)/(b), RC$^+$(a)/(b), ST$^+$(a)/(b), etc. etc.)

– For spaces with designation of factuality and/or possibility, there are some natural moderate constraints which may be imposed on coordination. E.g.: that a bilateral proposition cannot have both a factual truthmaker and a factual falsity-maker; that a bilateral proposition must have either some factual truthmaker or some factual falsity-maker; that a possible state cannot be both truthmaker and falsity-maker of the same proposition; that (more strongly) the fusion of a truthmaker and a falsity-maker of the same bilateral proposition is always impossible. And relatedly, with a.s.s.'s (where we have the notion of 'explicit contradictoriness' even without designation of anything) we might impose the strong coordination condition: the fusion of a truthmaker and a falsity-maker of the same bilateral proposition is always explicitly contradictory. – This is still not strong enough though to avoid the above problems with FDE$^+$(a): e.g. we can have $\sigma(\text{p}) = \langle\{\{\langle \text{u}, 1\rangle\}\}, \{\{\langle \text{u}, 0\rangle, \langle \text{u}', 1\rangle\}\}\rangle$, $\sigma(\text{q}) = \langle\{\{\langle \text{u}, 1\rangle\}\}, \{\{\langle \text{u}, 0\rangle, \langle \text{u}'', 1\rangle\}\}\rangle$.

(3) *Extension of the interpretation-function.* Next there is the issue of how to extend the interpretation-function to molecular formulas (and data expressions etc.). There is the style of 'basic' or 'exact' truthmaker semantics as in the main semantics we used for FDE and ET; the 'regular' style as in the semantics for RT; 'semiregular' with closure only under fusions as in the semantics for ST, or alternatively with closure only under intermediates; etc. As we have already observed before, the choice here is very much connected (in a reasonable construction) to the corresponding choice in (2).

(4) *Variety of notions of consequence, equivalence, etc.* Finally, when we come to notions of consequence, equivalence, etc., the variations from the earlier variation-points, together with a few new tweaks which we can make here, yield a large, even embarrassingly large, variety of notions. In what follows I mention only some conspicuous examples of such notions, hoping to give some idea of the relevant scope. (A more systematic classification and investigation of the classified notions, and corresponding calculi, would of course be interesting.) – To simplify the discussion we stick here with formulas; but similar distinctions apply also in the cases of data etc.

We have already considered numerous notions of the form

'(W.r.t. every model:) Every truthmaker of A is truthmaker of B',

their duals

'Every falsity-maker of B is falsity-maker of A',

and the corresponding equivalential forms

'A and B have the same truthmakers',

'A and B have the same falsity-makers',

all for various senses of truthmaking and falsity-making.

– If we have also the resources of designation of factuality and/or possibility, or even just the notion of explicit contradictoriness from atomic situation spaces, many further notions become definable. Here are just a few illustrative examples: –

'Every non--explicitly-contradictory (or, for frames with designation of possibility: possible) truthmaker of A includes truthmaker of B'.

Something like this should characterize the 'strong Kleene logic' K3. The dual system as is well known is Priest's 'Logic of Paradox' LP. (See e.g. Priest 2008 §§ 7.3 & 7.4 for definitions of these via three-valued semantics.) A direct corresponding definition of LP consequence in the present environment seems to be:

'Every non--explicitly-contradictory (or: possible) falsity-maker of B includes falsity-maker of A'.

(It seems clear that a suitable 'K3DC' would be simply DC minus EM, and LPDC dually DC minus NC; and similarly for dual data etc.)

'Every fusion of truthmaker of A with falsity-maker of B is explicitly contradictory (or: impossible)'.

This seems to be good old classical consequence.

'Every factual truthmaker of A includes truthmaker of B'.

'For every model M: If M verifies A then M verifies B'.

– One might also play around with *conjunctions* of given consequence or equivalence notions. We have already mentioned before ET ∩ EF and ST ∩ SF, which are natural notions corresponding to sameness of the full bilateral proposition of appropriate kind. And one might also consider EE ∩ EC, EE ∩ (EC)$^{-1}$, SE ∩ SC, SE ∩ (SC)$^{-1}$, etc.

We may also consider some *dis*junctions (unions). These may be interesting *notions*, but note that they will not in general characterize a *system* in a reasonable sense of the term – which, in this context of binary relations between formulas (or data etc.) should no doubt include at least the requirements of reflexivity and transitivity, the latter of which may easily fail for such disjunctions. E.g. EE ∪ EC, SE ∪ SC, RE ∪ RC are interesting disjunctive notions, but not transitive: e.g. they all hold between p ∧ q and q, and between q and q ∨ r, but not between p ∧ q and q ∨ r. – One can of course always consider the *transitive closure* of the union; I guess this might or might not be an interesting system, or a 'new' system. In the present example it seems clear that tc(EE ∪ EC) = tc(SE ∪ SC) = tc(RE ∪ RC) = (good old) FDE.

CHAPTER 5

DIRECTIONAL DEDUCTION FOR RELEVANCE LOGICS:
QUANTIFICATIONAL AND MODAL EXTENSIONS

## 1. The calculi FDEQDC etc.

(To reduce somewhat the proliferation of systems and corresponding calculi, in the present chapter we omit consideration of the non-symmetric [Entailment and Containment] systems of the regular, semiregular, and exact families.)

(1) *Extensions of FDE calculi*.

(The basic rules of) the calculi FDEQDC, FDEQDDC, FDEQHC, FDEQCC are obtained from (the basic rules of) the corresponding standard calculi QDC, QDDC, QHC, QCC simply by omission of the NC and EM rules, with the single additional difference that the Criss-Crossing rule in FDEQHC is formulated in terms of 'strict' criss-crossing rather than in terms of 'weak' criss-crossing (criss-crossing of *pure* components) as in QHC.

Likewise, the calculi FDES5DC, FDES5DDC, FDES5HC, FDES5CC are obtained from the corresponding standard calculi S5DC, S5DDC, S5HC, S5CC by omission of the NC and EM rules, and with the strict form of Criss-Crossing for FDES5HC.

And finally in modal-quantificational case too the situation is similar: the calculi FDES5QDC, FDES5QDDC, FDES5QHC, FDES5QCC are obtained from the corresponding standard calculi S5QDC, S5QDDC, S5QHC, S5QCC by omission of NC and EM, and with the strict form of Criss-Crossing for FDES5QHC.

(2) *Extensions of RT calculi*.

We now omit, in each case (w.r.t. the standard calculus), not only the NC and EM rules, but also the Del, Exp, Inst and Gen (or GenInst and GenGen in Consecution Calculus) rules; we add as new primitive rules the inverses of ↑IVU and ↓EVE (so we may now speak of the rule-pairs ↕VU and ↕VE); we use Distribution rather than Criss-Crossing in the Hybrid and Consecution calculi (as in RTHC and RTCC); and we add the ↕Fus and ↕Int rules. – However, I don't know whether this is enough (i.e. gives *complete* calculi; they are certainly sound calculi): perhaps just as in basic propositional case we need ↕Int for some limited RT-equivalential use of Del and Exp, so also in quantificational/modal cases we need extra rules for some limited RTQ/S5/S5Q-equivalential use of Inst and Gen. (I do not see e.g. how to deduce ∀xPx from [its RTQ-equivalent] ∀xPx, ∃xPx in a calculus RTQDC as described above without any supplementation.)

The following rules ↕RUI and ↕REI, for inclusion and exclusion of redundant instances of universals and existentials, are natural candidates for the role of such supplementary primitive rules. (I don't know whether they are enough for completeness. But they do provide deductions in such cases as ∀xPx, ∃xPx and ∀xPx, or □p ∧ p and □p, etc.)

$$(\text{↕RUI}) \quad \begin{array}{c} (E\underline{f}) \ (\underline{x}, y): \alpha(y), \alpha(t) \\ \updownarrow \\ (E\underline{f}) \ (\underline{x}, y): \alpha(y) \end{array} \qquad (\text{↕REI}) \quad \begin{array}{c} (E\underline{f}, g) \ (\underline{x}): \alpha(g\underline{x}) \ / \ \alpha(t) \\ \updownarrow \\ (E\underline{f}, g) \ (\underline{x}): \alpha(g\underline{x}) \end{array}$$

(where the comma in the first matrix indicates the 'conjunction', by 'distribution', of the two data). (Similarly for dual data, and in modal cases.) In QDC (or FDEQDC), ↓RUI is derivable using ↓Del, ↑RUI using ↓Inst, ↓REI using ↑Gen, and ↑REI using ↑Exp.

(3) *Extensions of ST calculi.*

As with RT above except that we do not add the ↕Fus and ↕Int rules (↕Fus is here derived rule). In particular ↕RUI and ↕REI are still sound rules and plausible candidates for primitive rules.

(4) *Extensions of ET calculi.*

Here perhaps we should have the G-rules (for sms data) supplemented by ↕IR, ↕VE, and ↕REI. Note that now the rules ↕VU and ↕RUI are *not* sound (w.r.t. the semantics to be given in the next section). – *Question*: In the appropriate notion of datum here, is it only the components of the matrix which should be multisets, or also the universal prefix (or whatever)?

## 2. Truthmaker semantics for FDEQ, RTQ, STQ, ETQ

W.r.t. e.s.s. S and domain D ('individuals'), and for n ≥ 0, we define an n-ary *bilateral propositional function* as a function assigning a bilateral proposition (w.r.t. S) to each n-tuple of elements of D. And we say that a bilateral propositional function is *regular* if all its values are regular bilateral propositions; and *semiregular* if all its values are semiregular bilateral propositions.

An *FDEQ model* is a triple (D, S, σ) where D in non-empty set ('individuals'), S is e.s.s., and σ is function assigning: individuals to individual variables; n-ary functions over D to n-ary functional variables; and n-ary bilateral propositional functions to n-ary predicate variables. The interpretation function σ can then be extended to all terms by the obvious condition, viz.

$$\sigma(f(t_1 \ldots t_n)) \;=\; \sigma(f)(\sigma(t_1) \ldots \sigma(t_n)),$$

and to all formulas as follows. For atomic formulas we put

$$\sigma(R(t_1 \ldots t_n)) \;=\; \sigma(R)(\sigma(t_1) \ldots \sigma(t_n)).$$

For negations, conjunctions and disjunctions we have the same clauses as before in our semantics for FDE. And for quantifications we put (where, for a ∈ D, σ[a/x] means the interpretation function which assigns a to x and in all other respects coincides with σ):

$$\sigma^+(\forall xA) \;=\; \{\cup C : C \text{ is collection of states selecting, for each } a \in D, \text{ a state from}$$
$\sigma[a/x]^+(A)\}$.

(Note that the relevant selection here is supposed to be made starting from the elements of D, not from the associated truthmaker sets in themselves. Two such sets might be identical and yet contribute distinct truthmakers to a fusion $\cup C$.)

$$\sigma^-(\forall xA) \;=\; \cup\{\rho^-(A) : \rho \text{ x-variant of } \sigma\}.$$

$$\sigma^+(\exists xA) \;=\; \cup\{\rho^+(A) : \rho \text{ x-variant of } \sigma\}.$$

σ⁻(∃xA) = {∪C : C is collection of states selecting, for each a ∈ D, a state from σ[a/x]⁻(A)}.

And exactly similar clauses extend σ also to data expressions etc.

We then define FDEQ consequence exactly as before with FDE:

α ⊨_FDEQ β  =_df  ∀ model (D, S, σ): ∀Σ ∈ σ⁺(α) ∃Σ′ ∈ σ⁺(β): Σ′ ⊆ Σ.

And same for dual data (δ ⊨_FDEQ ε), datum and dual datum (α ⊨_FDEQ δ), and validity of consecution (⊨_FDEQ α ⇒ δ).

And again as before, we define for sms data α, β:

α ≃_ETQ β  =_df  ∀ model (D, S, σ): σ⁺(α) = σ⁺(β).

α ≃_EFQ β  =_df  ∀ model (D, S, σ): σ⁻(α) = σ⁻(β).

– Now a *regular model* (or *RTQ model*) is an FDEQ model (D, S, σ) where σ assigns only *regular* bilateral propositional functions to the predicate variables. We then extend σ to terms and atomic formulas and negations as in the FDEQ semantics above, and to conjunctions and disjunctions and quantifications and the corresponding 'structural' devices of data expressions etc. as above only that (as in our semantics for RT) we use the prefix Int Cl Fus Cl in each clause.

We then define RT equivalence and RF equivalence for data by:

α ≃_RT β  =_df  ∀ regular model (D, S, σ): σ⁺(α) = σ⁺(β).

α ≃_RF β  =_df  ∀ regular model (D, S, σ): σ⁻(α) = σ⁻(β).

And similarly for dual data etc.

And there is the 'regular semantics' version of FDEQ consequence (which I guess coincides in extension with ⊨_FDEQ, as in propositional case):

α ⊨_FDEQr β  =_df  ∀ regular model (D, S, σ): ∀Σ ∈ σ⁺(α) ∃Σ′ ∈ σ⁺(β): Σ′ ⊆ Σ.

Now a *semiregular model* (or *STQ model*) is an FDEQ model (D, S, σ) where σ assigns only semiregular bilateral propositional functions to the predicate variables. We then extend σ to other expressions as above with regular models except that now we use the prefix Fus Cl rather than Int Cl Fus Cl. We then define STQ equivalence and SFQ equivalence for data by:

$\alpha \simeq_{STQ} \beta \ =_{df} \ \forall$ semiregular model (D, S, σ): $\sigma^+(\alpha) = \sigma^+(\beta)$.

$\alpha \simeq_{SFQ} \beta \ =_{df} \ \forall$ semiregular model (D, S, σ): $\sigma^-(\alpha) = \sigma^-(\beta)$.

And similarly for dual data etc.

### 3. Herbrand theorem for FDEQ; completeness of FDEQHC and FDEQCC

Let M = (S, σ) be an FDE model for the quantifier-free formulas of classical predicate logic with functional variables (without identity), i.e. e.s.s. together with assignment of bilateral propositions to the atomic formulas. And let N = (D, S′, ρ) be FDEQ model. We say that M and N (or N and M) *coincide* if ∀ atomic formula A: σ(A) = ρ(A) (whence ∀ quantifier-free formula B: σ(B) = ρ(B)).

*Proposition* (Coincidence Lemma). (1) ∀ FDEQ-model N ∃ FDE-model M s.t. N and M coincide.

(2) ∀ FDE-model M ∃ FDEQ-model N s.t. M and N coincide.

(3) ∀ quantifier-free Γ, Δ: Γ ⊨$_{FDEQ}$ Δ (i.e.: w.r.t. every FDEQ-model, every truthmaker of Γ includes truthmaker of some element of Δ) iff Γ ⊨$_{FDE}$ Δ.

*Proof.* (1) Immediate: take the same e.s.s. S as in the FDEQ model (D, S, ρ) and define σ by: (∀ atomic A): σ(A) = ρ(A).

(2) Let M = (S, σ). Then put N = (H, S, ρ) where: H is the Herbrand universe for the set of *all* functional variables, i.e. H is the set of all terms; and ρ is defined by:

(i) For individual variable x: ρ(x) = x.

(ii) For n-ary functional variable f: ρ(f) = {⟨⟨$t_1$ … $t_n$⟩, f($t_1$ … $t_n$)⟩ : $t_1$ … $t_n$ any terms}.

(iii) For k-ary predicate variable R: ρ(R) = {⟨⟨$t_1$ … $t_k$⟩, π⟩ : σ(R$t_1$ … $t_k$) = π}.

– Thus clearly N will coincide with M.

(3) Clearly, if M and N coincide then: it is true w.r.t. M that every truthmaker of Γ includes truthmaker of some element of Δ iff it is true w.r.t. N that every truthmaker of Γ includes truthmaker of some element of Δ. So given (1) and (2), (3) follows. □

*Proposition* (Herbrand theorem for FDEQ). For A$\underline{x}$ and B$\underline{y}$ quantifier-free and F = the set of functional variables occurring in either A$\underline{x}$ or B$\underline{y}$: $\forall \underline{x}A\underline{x} \vDash_{FDEQ} \exists \underline{y}B\underline{y}$ iff the sets $\Gamma$, $\Delta$ of their respective Herbrand instances w.r.t. F are s.t. $\Gamma \vDash_{FDE} \Delta$.

*Proof.* ($\Leftarrow$) Suppose $\Gamma \vDash_{FDE} \Delta$. So by the Coincidence Lemma, $\Gamma \vDash_{FDEQ} \Delta$. (Unfortunately we cannot now just say that we have $\forall \underline{x}A\underline{x} \vDash \wedge\Gamma \vDash \vee\Delta \vDash \exists \underline{y}B\underline{y}$ and thus $\forall \underline{x}A\underline{x} \vDash \exists \underline{y}B\underline{y}$, since there are no such formulas here as $\wedge\Gamma$ and $\vee\Delta$. But we can give a sort of corresponding semantical argument, as follows. [We could also appeal to the *compactness* of FDEQ, and use $\wedge\Gamma_0 \vDash \vee\Delta_0$; but this is not necessary – the following argument is much more direct.]) Let M = (D, S, $\sigma$) be arbitrary FDEQ model and let $\Sigma$ be truthmaker of $\forall \underline{x}A\underline{x}$ w.r.t. M (i.e. $\Sigma \in \sigma^{+}(\forall \underline{x}A\underline{x})$). Since each element of $\Gamma$ is instance of, hence FDEQ implied by, $\forall \underline{x}A\underline{x}$, $\Sigma$ includes truthmaker of each such element of $\Gamma$, and hence includes truthmaker of $\Gamma$. So, since $\Gamma \vDash_{FDEQ} \Delta$, $\Sigma$ includes truthmaker of some element C of $\Delta$. But C $\vDash_{FDEQ} \exists \underline{y}B\underline{y}$, and so $\Sigma$ includes truthmaker of $\exists \underline{y}B\underline{y}$, as desired.

($\Rightarrow$) Suppose $\Gamma \nvDash_{FDE} \Delta$, and let M = (S, $\sigma$) be a corresponding FDE countermodel. Construct N = (H$^F$, S, $\rho$) from M as in the proof of the Coincidence Lemma except that we have H$^F$, i.e. the Herbrand universe for F (rather than the full H), and correspondingly in the definition of $\rho$ we set (where x is the alphabetically first variable) $\rho(f) = \{\langle\langle t_1 \ldots t_n\rangle, x\rangle : t_1 \ldots t_n \in H^F\}$ when $f \notin F$. Thus N coincides with M, and so it is *not* the case that, w.r.t. N, every truthmaker of $\Gamma$ includes truthmaker of some element of $\Delta$. But w.r.t. N, the truthmakers of $\Gamma$ are precisely the truthmakers of $\forall \underline{x}A\underline{x}$, and state is truthmaker of some element of $\Delta$ iff it is truthmaker of $\exists \underline{y}B\underline{y}$. So w.r.t. N it is *not* the case that every truthmaker of $\forall \underline{x}A\underline{x}$ includes truthmaker of $\exists \underline{y}B\underline{y}$; and so $\forall \underline{x}A\underline{x} \nvDash_{FDEQ} \exists \underline{y}B\underline{y}$. $\square$

*Proposition* (Deductive completeness of FDEQHC). If $\alpha \vDash_{FDEQ} \delta$ then $\alpha \vdash_{FDEQHC} \delta$.

*Proof.* Given the previous results of compactness for FDE and Herbrand theorem for FDEQ, the rest of our proof in the classical case of QHC goes through essentially unchanged. $\square$

*Proposition* (Completeness of FDEQCC). If $\vDash_{FDEQ} \alpha \Rightarrow \delta$ (regular consecution) then $\vdash_{FDEQCC} \alpha \Rightarrow \delta$.

*Proof*.   Again, given compactness for FDE and Herbrand theorem for FDEQ, rest is as before.  □

*Question*.   Are FDEQDC and FDEQDDC also deductively complete? (I expect so, but don't have proof.)

# 4. Truthmaker semantics for modal systems

(1) *Semantics for FDES5 and ETS5.* W.r.t. domain W ('worlds') and e.s.s. S, we make the following definitions. A *relativized elementary situation* is a pair $\langle w, e \rangle$ where $w \in W$, $e \in S$. We usually write $\langle w, e \rangle$ more briefly as $w : e$. (Intuitively, $w : e$ should be thought of as: 'in w, e (holds)'.) We sometimes refer to elementary situations (i.e. elements of S) here as *absolute elementary situations* for emphasis. A *state* is now a set each of whose elements is either an absolute or a relativized elementary situation. A *bilateral proposition* is a pair $\langle X, Y \rangle$ where X and Y are sets of states with $X \neq Y$. – We extend the 'w :' notation to other types of object in the natural way:

(Relativized elementary situations:) $w : (v : e) =_{df} v : e$.

(States:) $w : \Sigma =_{df} \{w : e \mid e \in \Sigma\}$.

(Sets of states:) $w : X =_{df} \{w : \Sigma \mid \Sigma \in X\}$.

(Bilateral propositions:) $w : \langle X, Y \rangle =_{df} \langle w : X, w : Y \rangle$.

Now an *FDES5 model* (or *ETS5 model*, or *EFS5 model*) is a triple $(W, S, \sigma)$ where W ('worlds') is non-empty set, S is e.s.s., and $\sigma$ is function assigning: to world-variables, worlds; to n-ary functional variables, n-ary functions over W; and to propositional variables, bilateral propositions. – The function $\sigma$ can then be extended to all formulas and data expressions etc. as follows. For negations, conjunctions and disjunctions we have the same clauses as before with FDE models. For relative truth and modalities we put:

$\sigma(T^w A) = \langle \sigma(w) : \sigma^+(A), \sigma(w) : \sigma^-(A) \rangle$.

$\sigma^+(\square A) = \{\cup C \mid C$ is collection of states selecting, for each $v \in W$, a state from $\sigma[v/w]^+(T^w A)\}$.

$\sigma^-(\square A) = \cup\{\rho^-(T^w A) \mid \rho$ w-variant of $\sigma\}$.

$\sigma^+(\Diamond A) = \cup\{\rho^+(T^wA) \mid \rho \text{ w-variant of } \sigma\}$.

$\sigma^-(\Diamond A) = \{\cup C \mid C \text{ is collection of states selecting, for each } v \in W, \text{ a state from}$ $\sigma[v/w]^-(T^wA)\}$.

And similarly for the structural devices of data expressions etc. In particular for the prefixes (w̲), (Ef̲) etc. we give the natural clauses as in our truthmaker semantics for FDEQ.

We then define:

$\alpha \vDash_{FDES5} \beta =_{df}$ Either (i) $\forall$ model $(W, S, \sigma)$: $\forall \Sigma \in \sigma^+(\alpha) \exists \Sigma' \in \sigma^+(\beta)$: $\Sigma' \subseteq \Sigma$, or (ii) (where w is the first world-variable occurring in neither $\alpha$ nor $\beta$) $\forall$ model $(W, S, \sigma)$: $\forall \Sigma \in \sigma^+(T^w\alpha) \exists \Sigma' \in \sigma^+(T^w\beta)$: $\Sigma' \subseteq \Sigma$.

And similarly for dual data etc.

*Remark*. Here (i) always implies (ii), though not vice versa (take e.g. $\alpha = \Box p$, $\beta = p$). So (ii) alone would be an equivalent definiens. Still the above bipartite formulation is perhaps preferable, matching the bipartite character of the definition of FDES5DC (etc.) deductions (as already of S5DC etc. deductions). ⊣

$\alpha \simeq_{ETS5} \beta =_{df}$ Either (i) $\forall$ model $(W, S, \sigma)$: $\sigma^+(\alpha) = \sigma^+(\beta)$, or (ii) (where w is the first world-variable occurring in neither $\alpha$ nor $\beta$) $\forall$ model $(W, S, \sigma)$: $\sigma^+(T^w\alpha) = \sigma^+(T^w\beta)$.

And similarly for $\simeq_{EFS5}$: just replace $\sigma^+$ by $\sigma^-$ throughout.

I assume that the appropriate completeness results here can be obtained from the corresponding Herbrand and compactness theorems as in the case of FDEQ above, though I have not checked this in detail.

(2) *Semantics for RTS5 and STS5*. Like the above semantics for FDES5 and ETS5 except that (i) the interpretation function must assign *regular*, resp. *semiregular*, bilateral propositions to the propositional variables; and (ii) the prefix Int Cl Fus Cl, resp. Fus Cl, is inserted in the appropriate clauses for extending the interpretation function. Then we put:

α ≃$_{RTS5}$ β  =$_{df}$ Either (i) ∀ regular model (W, S, σ): σ⁺(α) = σ⁺(β), or (ii) (where w is the first world-variable occurring in neither α nor β) ∀ regular model (W, S, σ): σ⁺(Tᵂα) = σ⁺(Tᵂβ).

And similarly for dual data etc. And for ≃$_{RFS5}$ just replace σ⁺ by σ⁻ throughout.

*Remark.* Again, (i) implies (ii) but not vice versa (take e.g. α = □p ∧ p, β = □p), and so on. ⊣

≃$_{STS5}$ and ≃$_{SFS5}$ are defined similarly: just replace 'regular model' by 'semiregular model'.

(3) *Semantics for FDES5Q and ETS5Q.* We combine the ideas of previous constructions in the natural way. Thus an FDES5Q model (or ETS5Q model, or EFS5Q model) is quadruple (D, W, S, σ) where D and W are non-empty sets, S is e.s.s., and σ assigns: to individual variables, individuals; to world-variables, worlds; to functional variables, functions of appropriate type; and to n-ary predicate variables, n-ary bilateral propositional functions (where the sense of 'bilateral proposition' is as in the semantics for FDES5 above). And so on.

I haven't checked details of the completeness arguments, but think there should be no significant difficulties.

(4) *Semantics for RTS5Q and STS5Q.* Like (3) except for (i) use of *regular*, resp. *semiregular*, bilateral propositional functions and (ii) appropriate insertions of Int Cl Fus Cl, resp. Fus Cl. Then ≃$_{RTS5Q}$, ≃$_{RFS5Q}$, ≃$_{STS5Q}$, ≃$_{SFS5Q}$ are defined as in (2) above.

# SOME MAIN OPEN PROBLEMS

Among the various open questions mentioned in the course of this work (often with the ceremony of the prefix *Question(s)*, sometimes in other ways), some are more casual, others more significant. The following is a selection of some such questions which it would be particularly interesting to have answers to.

(1)  Does replacement of EM by EM* or some other reasonable modification of DC yield a *normal*-deduction complete calculus? (Pages 44–45.)

(2)  Is QDC deductively complete? (Page 79.)

(3)  Can the treatment of transition in QHC be improved? (Pages 80–81.)

(4)  Does S5Q satisfy Datum Interpolation? (Page 93.)

(5)  Prove our conjectures (or suitable modifications thereof) connecting ST and ET to corresponding classes of algebras. (Pages 135–136.)

(6)  Provide suitable completion to the formulations of RTQDC and STQDC (and ETQDC if needed). (Pages 150–152.)

# REFERENCES

A. R. Anderson & N. D. Belnap, 'Tautological entailments', *Philosophical Studies*, vol. 13, 1962, pp. 9–24. (See also Anderson & Belnap 1975, §§ 15.1–16.1.)

A. R. Anderson & N. D. Belnap, *Entailment*, vol. I, Princeton U. P., 1975.

A. R. Anderson, N. D. Belnap & J. M. Dunn, *Entailment*, vol. II, Princeton U. P., 1992.

R. B. Angell, 'Three systems of first degree entailment' (abstract), *Journal of Symbolic Logic*, vol. 42, 1977, p. 147.

R. B. Angell, 'Deducibility, entailment and analytic containment', in J. Norman & R. Sylvan (eds.), *Directions in Relevant Logic*, Dordrecht: Kluwer, 1989, pp. 119–143.

R. Batchelor, 'Grounds and consequences', *Grazer Philosophische Studien*, vol. 80, 2010, pp. 65–77.

R. Batchelor, 'Complexes and their constituents', *Theoria*, vol. 79, 2013, pp. 326–352.

N. D. Belnap, 'Tautological entailments' (abstract), *Journal of Symbolic Logic*, vol. 24, 1959, p. 316.

N. D. Belnap, 'A useful four-valued logic', in J. M. Dunn & G. Epstein (eds.), *Modern Uses of Multiple-Valued Logic*, Dordrecht: Reidel, 1977a, pp. 8–37. (Repr. in Omori & Wansing 2019. See also Anderson, Belnap & Dunn 1992, § 81.)

N. D. Belnap, 'How a computer should think', in G. Ryle (ed.), *Contemporary Aspects of Philosophy*, Stocksfield: Oriel Press, 1977b, pp. 30–55. (Repr. in Omori & Wansing 2019. See also Anderson, Belnap & Dunn 1992, § 81.)

R. Carnap, *Formalization of Logic*, Harvard U. P., 1943.

F. Correia, 'Semantics for analytic containment', *Studia Logica*, vol. 77, 2004, pp. 87–104.

F. Correia, 'Grounding and truth-functions', *Logique et Analyse*, vol. 53, 2010, pp. 251–279.

F. Correia, 'On the logic of factual equivalence', *Review of Symbolic Logic*, vol. 9, 2016, pp. 103–122.

W. Craig, 'Linear reasoning: A new form of the Herbrand-Gentzen theorem', *Journal of Symbolic Logic*, vol. 22, 1957, pp. 250–268.

J. M. Dunn, 'Intuitive semantics for first-degree entailments and "coupled trees"', *Philosophical Studies*, vol. 29, 1976, pp. 149–168. (Repr. in Anderson, Belnap & Dunn 1992 [§ 50] and in Omori & Wansing 2019.)

K. Fine, 'Failures of the interpolation lemma in quantified modal logic', *Journal of Symbolic Logic*, vol. 44, 1979, pp. 201–206.

K. Fine, 'Analytic implication', *Notre Dame Journal of Formal Logic*, vol. 27, 1986, pp. 169–179.

K. Fine, 'Guide to ground', in F. Correia and B. Schnieder (eds.), *Metaphysical Grounding*, Cambridge U. P., 2012, pp. 37–80.

K. Fine, 'Truth-maker semantics for intuitionistic logic', *Journal of Philosophical Logic*, vol. 43, 2014, pp. 549–577.

K. Fine, 'Unified foundations for essence and ground', *Journal of the American Philosophical Association*, vol. 1, 2015, pp. 296–311.

K. Fine, 'Angellic content', *Journal of Philosophical Logic*, vol. 45, 2016, pp. 199–226.

K. Fine, 'A theory of truthmaker content I: Conjunction, disjunction and negation', *Journal of Philosophical Logic*, vol. 46, 2017a, pp. 625–674.

K. Fine, 'Truthmaker semantics', in B. Hale, C. Wright & A. Miller (eds.), *A Companion to the Philosophy of Language*, 2nd ed., Vol. II, Oxford: Wiley-Blackwell, 2017b, pp. 556–577.

K. Fine, 'Some remarks on Bolzano on ground', in S. Roski and B. Schnieder (eds.), *Bolzano's Philosophy of Grounding*, Oxford U. P., 2022, pp. 276–300.

K. Fine & M. Jago, 'Logic for exact entailment', *Review of Symbolic Logic*, vol. 12, 2019, pp. 536–556.

J. M. Font, 'Belnap's four-valued logic and De Morgan lattices', *Logic Journal of the IGPL*, vol. 5, 1997, pp. 413–440.

R. C. Jeffrey, *Formal Logic*, 1st ed., New York: McGraw-Hill, 1967.

R. C. Jeffrey, *Formal Logic*, 2nd ed., New York: McGraw-Hill, 1981.

W. Kneale, 'The province of logic', in H. D. Lewis (ed.), *Contemporary British Philosophy*, Third Series, London: Allen & Unwin, 1956, pp. 235–261.

W. Kneale & M. Kneale, *The Development of Logic*, paperback ed., Oxford U. P., 1984.

S. Krämer, 'Truthmaker equivalence', to appear in S. Leuenberger & A. Rieger (eds.), *Themes from Alan Weir*, Cham: Springer.

R. N. McKenzie, G. F. McNulty & W. F. Taylor, *Algebras, Lattices, Varieties*, Vol. I, Monterey, California: Wadsworth & Brooks/Cole, 1987.

H. Omori & H. Wansing (eds.), *New Essays on Belnap-Dunn Logic*, Cham: Springer, 2019.

G. Priest, *An Introduction to Non-Classical Logic*, 2nd ed., Cambridge U. P., 2008.

I. Rumfitt, '"Yes" and "No"', *Mind*, vol. 109, 2000, pp. 781–823.

D. J. Shoesmith & T. J. Smiley, *Multiple-Conclusion Logic*, Cambridge U. P., 2nd printing, 1980.

T. J. Smiley, 'Rejection', *Analysis*, vol. 56, 1996, pp. 1–9.

R. M. Smullyan, 'A unifying principle in quantification theory', *Proceedings of the National Academy of Sciences* (USA), vol. 49, 1963, pp. 828–832.

R. M. Smullyan, *First-Order Logic*, Berlin: Springer, 1968 (repr. Dover, 1995).

B. van Fraassen, 'Facts and tautological entailments', *Journal of Philosophy*, vol. 66, 1969, pp. 477–487. (Repr. in Anderson & Belnap 1975 [§ 20.3].)